\documentclass[11pt,twoside]{article}
\usepackage{amsmath,amsfonts,amssymb,cases,mathdots,color,colordvi,url,tikz}
\usepackage{hyperref}
\usepackage{bm}
\usepackage{bbm}
\usepackage{pgfpages}
\usepackage{tabularx}
\usepackage{graphics}%\textwidth 16.cm
\usepackage{subfigure}
\usepackage{epsfig}
\usepackage{algorithm}
\usepackage{algorithmic}
\usepackage{physics}
\usepackage{epstopdf}
\usepackage{parskip}
\usepackage{footmisc}
\usepackage{tablefootnote}
\usepackage[justification=centering]{caption}

\pdfoutput=1
\newtheorem{proposition}{Proposition}%[section]
\newtheorem{theorem}[proposition]{Theorem}
\newtheorem{lemma}[proposition]{Lemma}
\newtheorem{corollary}[proposition]{Corollary}
\newtheorem{definition}[proposition]{Definition}

\newtheorem{Remark}{Remark}
\newtheorem{assumption}{Assumption}
\newtheorem{example}{Example}

\newcommand{\be}{\begin{equation}}
	\newcommand{\ee}{\end{equation}}
\newcommand{\ba}{\begin{eqnarray}}
	\newcommand{\ea}{\end{eqnarray}}
\newcommand{\bas}{\begin{eqnarray*}}
	\newcommand{\eas}{\end{eqnarray*}}

\graphicspath{ {./Figures/} }
\def\bfq{{\bf q}}

\def\bfa{{\bf a}}
\def\bfb{{\bf b}}

\def\bfd{{\bf d}}
\def\bfe{{\bf e}}

\def\bfv{{\bf v}}
\def\bfx{{\bf x}}

\def\bfz{{\bf z}}
\def\bfw{{\bf w}}
\def\bfu{{\bf u}}

\def\H{{\cal H}}

\def\M{{\cal M}}
\def\N{{\cal N}}

\def\R{{\cal R}}
\def\T{{\cal T}}
\def\K{{\mathcal{K}}}
\def\F{{\mathcal{F}}}
\def\J{{\mathcal{J}}}
\def\Re{\mathbb{R}}

\def\bfxi{{\boldsymbol{\xi}}}

\def\bfone{{\bf 1}}
\def\bfzero{\bf 0  }
\def\Prox{{\rm Prox}}
\def\Proj{{\rm Proj}}

\def\bfone{{\bf 1}}
\def\bfzero{\bf 0  }

\def\oG{{\overline{\Gamma}}}

\def\proj{{\rm Proj}_{\mathbb{S}}}
\def\prox{{\rm Prox}_{\beta \lambda \|(\cdot)_+\|_0}}
\def\oT{\overline{T}} 
\def\mI{\mathcal{I}^*_-}
\def\moI{\overline{\mathcal{I}}^*_-}
\def\bp{{\bf Proof.}}

\def\nwg{\nabla_{\bfw} g}
\def\nyxi{\nabla_{\bfxi} g}
\def\tbu{ \widetilde{\bfu} }
\def\hbw{ \widehat{\bfw} }
\def\hbu{ \widehat{\bfu} }

\def\hbxi{ \widehat{\bfxi} }
\def\hmI{ \widehat{\mathcal{I}} }
\def\hT{\widehat{T}}

\begin{document}

\title{Sparse SVM with Hard-Margin Loss: a Newton-Augmented Lagrangian Method in Reduced Dimensions}

%\author{\name Penghe Zhang \email 19118011@bjtu.edu.cn
%	\AND
%       \name Naihua Xiu \email nhxiu@bjtu.edu.cn \\
%       \addr School of Mathematics and Statistics \\
%             Beijing Jiaotong University \\
%             Beijing, China
%       \AND
%       \name Hou-Duo Qi \email houduo.qi@polyu.edu.hk \\
%       \addr Department of Applied Mathematics\\
%       The Hong Kong Polytecnic University\\
%       Kowloon,  Hong Kong}
   
   \author{Penghe Zhang\thanks{School of Mathematics and Statistics, Beijing Jiaotong University, Beijing 100044, PR China, E-mail: {19118011@bjtu.edu.cn} },
   	\ \ Naihua Xiu\thanks{School of Mathematics and Statistics, Beijing Jiaotong University, Beijing 100044, PR China, E-mail: {nhxiu@bjtu.edu.cn} } \ \ and \ \
   	Hou-Duo Qi\thanks{Department of Applied Mathematics, The Hong Kong Polytechnic University, Hung Hom, Hong Kong, E-mail: {houduo.qi@polyu.ac.hk} }
   }

%\editor{ }

\maketitle

\begin{abstract}
The hard margin loss function has been at the core of the support vector machine (SVM)
research from the very beginning due to its generalization capability.
On the other hand, the cardinality constraint has been widely used for feature selection, leading to sparse solutions. This paper studies the sparse SVM with the hard-margin loss (SSVM-HM) that integrates the virtues of both worlds.
However, SSVM-HM is one of the most challenging models to solve. In this paper,
we cast the problem as a composite optimization with the cardinality constraint.
We characterize its local minimizers in terms of {\rm P}-stationarity that well captures the combinatorial 
structure of the problem. 
We then propose an inexact proximal augmented Lagrangian method (iPAL).
The different parts of the inexactness measurements from the {\rm P}-stationarity are controlled at
different scales in a way that the generated sequence converges both globally and at a linear rate.
This matches the best convergence theory for composite optimization. 
To make iPAL practically efficient, we propose a gradient-Newton method in a subspace for the iPAL subproblem. 
This is accomplished by detecting active samples and features with the help of the proximal operator of the hard margin loss and the projection of cardinality constraint. 
Extensive numerical results on both simulated and real datasets demonstrate that the proposed method is
fast, produces sparse solution of high accuracy, and can lead to effective reduction on active samples and features  when compared with several leading solvers.
\end{abstract}


%\vskip 1cm
%\noindent
%{\bf Running Title:} Convex Optimization of Low Dimensional Euclidean Distances
%%%%%%%%%%%%%%%%%%%%%%%%%%%%%%%%%%%%%%%%%%%%%%%%%%%%%%%%%%%%%%%%%%%%%%%%%%%%%%%%%%%%%%%%%%%%%%%%%%%
%\newpage
\section{Introduction}

This paper is concerned with one of the most challenging formulations in the study of support vector machines
(SVM):
\be\label{HM-SSVM}
\min_{\bfw\in \Re^n, b \in \Re } \; \frac 12 \|\bfw\|^2 + \frac 12 b^2 
     + \lambda \sum_{i=1}^m h \Big( 1 - y_i (\langle \bfx_i, \bfw \rangle +b) \Big), \ \
     \mbox{s.t.} \ \ \mathbb{S} := \left\{ \bfw \in \Re^n \ | \  \| \bfw\|_0 \le s \right\},
\ee 
where $\{(\bfx_i, y_i) \}_{i=1}^m$ are the sample data with $\bfx_i \in \Re^n$ and $y_i \in \{1, -1\}$ being its label.
The separating hyperplane is $\langle \bfx, \bfw \rangle = b$ and the loss function is the hard-margin loss:
\[
  h(t) = \left\{
  \begin{array}{ll}
  	 1  & \mbox{if} \ t > 0 \\
  	 0  & \mbox{if} \ t \le 0 .
  \end{array} 
  \right .
\]
Furthermore, the model aims to seek a hyperplane of sparse features selected by the $\ell_0$-norm
$\| \cdot\|_0$
with a user-specified sparsity level $s \ge 1$ and $\mathbb{S}$ is known as the $s$-sparse set. 
\cite{Vapnik1998} discussed the hard-margin loss (also known as the $0/1$-loss), which is to
construct the hyperplane that makes the smallest number of separating errors. However, the optimization of it
is NP-complete.
The use of $\ell_0$-norm is getting popular in selecting sparse features.
The first two terms in the objective is to maximize the separation gap in the $(\bfw, b)$ space rather than
in the feature space of $\bfw$. 
This objective has been promoted by Mangasarian and his collaborators (see, \cite{mangasarian2001lagrangian, fung2001proximal, lee2001ssvm}). Due to its strong convexity in both $\bfw$ and $b$, Newton's method 
has been the core of those studies for the ridge/hinge-loss function.
The purpose of this paper is to extend Newton's method to the sparse SVM with hard-margin loss under the framework of
augmented Lagrangian method with proved convergence. 
This section is organized as follows. We will first conduct a literature review, followed by an explanation of our 
numerical approach.

%%%%%%%%%%%%%%%%%%%%%%%%%%%
\subsection{Related work}

There exists extensive research on SVMs. We refer to \cite{Vapnik1998, cristianini2000introduction, smola2004tutorial, steinwart2008support, chang2011libsvm} for many of the models and the solvers.
We restrict our review to the sparse SVM with the hard-margin loss and the related numerical methods.
We split the papers into three groups. The first is the mixed-integer programming (MIP) approach.
The second group is to treat (\ref{HM-SSVM}) as a composite optimization and the augmented Lagrangian method is
a natural choice. The last group is on Newton's method for such composite optimization.

\noindent
{\bf (A) MIP and its convex relaxation.}
The advantage of simultaneously addressing the $0/1$-loss and the $\ell_0$-norm for feature selection
was thoroughly justified by \cite{ustun2016supersparse} for a medical scoring problem. 
In this application, both the solution accuracy (controlled by the $0/1$-loss) and 
solution sparsity (controlled by the $\ell_0$-norm) are crucial to  yield a reliable medical score.
The solution method is to reformulate the problem as a mixed integer programming (MIP) by using the old trick:
Big-M constraint on both the $0/1$-loss and the $\ell_0$-norm. We refer to 
\cite{ liittschwager1978integer, bajgier1982experimental, brooks2011support} for earlier works along this line.
Another trick for MIP reformulation is based on the following fact:
\be \label{NCP}
  \mbox{(complementarity reformulation)} \quad  \| \bfw\|_0   = \min_{\bfz \in \mathbb{R}^n} \sum_{i=1}^n z_i, \ \ s.t. \ w_i(1-z_i) =0, \ z_i \in [0, 1],
\ee
see \cite{feng2013complementarity, kanzow2022sparse}.
One potential drawback for the smooth approach is the drastic increase in the dimensionality, especially when 
Newton's method is applied, see Section 7.4 of \cite{kanzow2022sparse} for a numerical example. 
One can imagine that this drawback would get worse when the $0/1$-loss is also represented by the
complementarity reformulation. 
It is worth pointing out that exciting progress has been made in a recent MIP approach (e.g., via Big-M constraint) by \cite{dedieu2021learning}, who cleverly combines a continuous approach and MIP to develop
a fast algorithm for an $\ell_0$-norm minimization problem. 
It remains to be seen how the approach would be adapted to Problem (\ref{HM-SSVM}), which involves both $\ell_0$-norm
and the $0/1$-loss.

Extensive work has been done in relaxing the $\ell_0$-norm by its convex surrogate $\ell_1$-norm
see, e.g., \cite{zhu20031, fung2004feature, shao2019joint, yuan2010comparison, dedieu2022solving}.
Although the approximation models are easier to tackle, they may not exactly recover the solution to the original $\ell_0$-based model. 
For example, comparison studies on linear regression and convex quantile regression show that $\ell_0$-norm has better performance than $\ell_1$-norm on feature selection, see \cite{johnson2015risk,dai2023variable}.
Therefore, for applications that require higher solution accuracy, solving Problem (\ref{HM-SSVM}) directly
seems necessary as done in 
\cite{ustun2016supersparse}.
However, MIP approach has drawbacks on scalability and computational speed for Problem (\ref{HM-SSVM}).
%This may be the reason why no numerical experiments were reported on large problems in \cite{ustun2016supersparse}.

\noindent
{\bf (B) Augmented Lagrangian methods for nonconvex problems.}
From the perspective of constrained optimization, it is natural to consider the augmented Lagrangian method (ALM) of 
\cite{hestenes1969multiplier, powell1969method} for Problem (\ref{HM-SSVM}). 
ALMs have become standard textbook material (see, e.g., \cite{bertsekas2014constrained, nocedal2006numerical, birgin2014practical}).
However, direct application is not possible due to the problem being a type of nonsmooth, nonconvex, and composite 
optimization.
Despite this, significant progress has been recently made for this type of problems by \cite{bolte2018nonconvex}:
\be \label{COP}
 \min \ f(\bfx) + \theta(F(\bfx)),
\ee 
where $f: \Re^n \mapsto \Re$ is continuously differentiable ($C^1$ class),
$F: \Re^n \mapsto \Re^m$ ($m \le n$) is also $C^1$, and $\theta: \Re^m \mapsto (-\infty, + \infty]$ is a proper and
lower-semicontinuous (lsc) function.
A key message delivered in \cite{bolte2018nonconvex} was that adaptive Lagrangian-based multiplier methods
can be developed with guaranteed convergence properties.
An essential requirement is that the primal iterates are kept close to the so-called
information zone, where certain regularity conditions are assumed.
This requirement is often met when the subproblems are solved exactly. 
Other developments also appear in \cite{li2015global,wang2018con,boct2020proximal} for unconstrained composite optimization.

Another possible solution method for (\ref{HM-SSVM}) is to follow the framework of the augmented Lagrangian method
of \cite{kanzow2021augmented, de2023constrained, jia2023augmented} for composite optimization covering the
cardinality constraint (i.e., $\ell_0$-norm constraint).
One of the techniques used is to represent the cardinality constraint as a smooth complementarity system 
in the spirit of (\ref{NCP}). Similarly, the hard-margin loss can also be represented by a
system of complementarity. 
This would drastically increase the dimensions of the resulting formulation.

Our problem (\ref{HM-SSVM}) can be put in the framework of (\ref{COP}) by making use of the indicator function on the
sparse constraint. The number of smooth functions in $F(\cdot)$ would be $(n+m)$, violating the requirement of
$m \le n$ in \cite{bolte2018nonconvex}. It is also not clear how the primal iterates would be kept close to the problem
information zone as we are simultaneously dealing with both the hard-margin loss and the $\ell_0$-norm constraint. 
Furthermore, Mangasarian's original proposal for introducing the quadratic objective in the $(\bfw, b)$ space is for
Newton's method to be used as its Hessian matrix is diagonal (i.e., sparse). 
Therefore, our proposal in this paper is to
develop an augmented Lagrangian method sharing similar convergence properties as in \cite{bolte2018nonconvex}
while allowing Newton's method to be used. Furthermore, we allow its
subproblems to be solved approximately.

\noindent
{\bf (C) Newton's method for composite optimization.}
We briefly discuss our own work on this aspect. For the application of compressed sensing with cardinality constraint,
we developed a Newton-based hard-thresholding method in \cite{zhou2021global}, which is also proved to be globally
convergent. For the hard-margin loss, we were only able to prove its local quadratic convergence in \cite{zhou2021quadratic}. 
Our recent attempt of \cite{zhang2023inalm} studies an ALM for a hard-margin loss composite optimization without any constraints.
The current paper can be seen as an extension to the constrained case with
the cardinality constraint.
%%%%%%%%%%%%%%%%%%%%%%%%%
%The difference between the cardinality constraint and the hard-margin loss is that the former is
%unconstrained (counting nonzeros on the whole real line) while the latter is constrained in the sense that the counting of nonzeros is restricted in the positive part of the real line.
%Furthermore, such counting is composited with a linear mapping. 
%Consequently, it is the case of unconstrained optimization vs constrained optimization.
%%%%%%%%%%%%%%%%%%%%%%%%%%%%%%%
Extension of optimization methods from unconstrained optimization to constrained counterpart is sometime very challenging. 
%Fortunately, the new framework of \cite{bolte2018nonconvex} opens the possibility of such extension 
%through the vehicle of ALMs. 
The difficulty lies with the challenge of simultaneously handling both the sparse set and the hard-margin loss, both of which are of combinatorial nature.
This paper successfully resolved this difficulty in
the venue of SVMs.

%%%%%%%%%%%%%%%%%%%%%%%%%%%%%%%%%%%%%%%%%%%%%%%%%%%
\subsection{Main contributions}

The review above establishes that SSVM-HM (\ref{HM-SSVM}) is a very useful yet
challenging model to solve.
There lacks efficient numerical methods for it especially for large data sets. 
Since we are not following the MIP approach, we are contented with being capable of computing 
a local minimizer. 
Our first contribution is on the characterization of local minimizers of (\ref{HM-SSVM}). This
is explained below with other innovative contributions.
% for finding a local minimizer.

{\bf (i) On the concept of stationarity.} Since Problem (\ref{HM-SSVM}) is essentially a nonconvex composite
optimization with a cardinality constraint. It has various formulations (e.g., via the complementarity 
systems as we review above). Stationary points can then be characterized for those reformulations. One good example to
follow is the recent paper of \cite{cui2022nonconvex}. We choose to define the stationarity through two 
proximal mapping involving the hard-margin loss function and the $s$-sparse set. We hence call it the 
{\rm P}-stationarity. This extends the previous stationarity concepts of \cite{beck2013sparsity, pan2015solutions, zhou2021global}
on sparse optimization to the hard-margin case.
Moreover, we establish one-to-one correspondence between {\rm P}-stationary points and local minimizers of (\ref{HM-SSVM}). This shows that {\rm P}-stationarity is adequate for Problem (\ref{HM-SSVM}).
% and there is no need to look for other concepts to enhance/replace it.

{\bf (ii) Inexact framework of proximal augmented Lagrangian method.}
To make the proposed ALM implementable, we solve its subproblem inexactly in a way that the generated iterates 
should enjoy the best known convergence properties, namely global convergence to 
a stationary point with a linear rate.
%(see \cite{bolte2018nonconvex}). 
It turns out that the accuracy of different parts of the stationarity measurement of the iterates should
satisfy certain relationship between them. In other words, a new set of computable stopping criteria for
solving each subproblem of ALM is developed.
Unlike the case where each subproblem is solved exactly in terms of satisfying its optimality condition, the
inexactness of the approximate solution creates some unavoidable obstacles in applying the traditional convergence analysis tools. 
A new Lyapunov function is constructed by adding a proximal term to the standard augmented Lagrangian to prove the global convergence as well as the linear rate of
convergence under certain regularity conditions often met by data with $n \gg m$.

{\bf (iii) Optimization methods in reduced dimensions.}
Since Problem (\ref{HM-SSVM}) is highly combinatorial defined by the both hard-margin loss and the sparse set,
a (local) solution should stay in a subspace when the iterates are close to it. This raises the question whether 
we can develop a subspace-based optimization method for each of the subproblems in the ALM framework. 
Intuitively, it is possible. However, for thus generated sequences to have good convergence properties as stated in
(ii) above requires delicate tracking of the true underlying space. 
We achieved this tracking by making use of a sharp observation that the optimal solution should satisfy some
complementarity conditions. Those conditions naturally define a subspace at each iteration. 
We then apply a gradient descent method  in this subspace to get a sufficient decrease in the Lyapunov function. 
To speed up the convergence, we further update the iterate by Newton's method in the same subspace. 
The generated iterate is guaranteed to meet the stopping criteria discussed in (ii). 
The Newton method enjoys the quadratic convergence under the assumption of strict complementarity condition.

The resulting algorithm is highly efficient and is benchmarked against several leading SVM solvers
on both simulated and real datasets. The proposed method is capable of computing a sparse solution with high classification accuracy and a smaller number of support vectors. 
And it is fast due to the fact that subproblems were
often solved in a much smaller subspace than the full space.

%%%%%%%%%%%%%%%%%%%%%%%%%%%%%%%%%%%%%%%%%%%%%%%%%%%%
\subsection{Organization}

In next section, we explain the notations used in the paper and present the basic properties of
the projection operator to the $s$-sparse set and the positive hard-thresholding operator for the hard-margin
loss function.
Section \ref{Section-Optimality} introduces the stationarity and characterizes it in terms of
the local minimizers of Problem (\ref{HM-SSVM}).
Section \ref{Section-PAL} develops the inexact framework of the proximal augmented Lagrangian method (iPAL) and
conducts its convergence analysis. 
In Section \ref{Section-Newton}, we propose an efficient numerical strategy to solve the subproblem 
in iPAL in a subspace.
The strategy consists of two parts: first apply a gradient descent to guarantee a sufficient decrease, followed 
by a Newton step. Both are computed in a well defined subspace. 
We also conduct convergence analysis of this numerical strategy. 
We report extensive numerical experiments in Section \ref{Section-Numerical}.

The new algorithmic framework does not rely on any external optimization solvers for its subproblems.
The design of the algorithm is constructive and is active-set based. It requires a new set of convergence analysis.
We provide all the detailed proofs in Appendix.

%%%%%%%%%%%%%%%%%%%%%%%%%%%%%%%%%%%%%%%%%%%%%%%%%%%%
\section{Preliminaries and Positive Hard-Thresholding Operator} \label{Section-Preliminaries}

%%%%
\subsection{Notation and Definitions} 
We use boldfaced lowercase letters to denote vectors. For example, $\bfw \in \mathbb{R}^n$ is a column vector of size $n$ and $\bfw^\top$ is its transpose.
Let $w_i$ or $[\bfw]_i$ denote the $i$th element of $\bfw$. The norm $\| \bfw\|$ denotes the Euclidean norm of $\bfw$ and for
a matrix $A$, $\|A\|$ is the induced norm by the Euclidean norm so that we always have
$
\| A\bfw \| \le \| A \| \| \bfw\|.
$
For two column vectors $\bfw$ and $\bfxi$, we use the Matlab notation $[\bfw; \bfxi]$ to 
denote the new column vector concatenating $\bfw$ and $\bfxi$.
The neighborhood of $\bfw^* \in \mathbb{R}^n$ with radius $\delta > 0$ is denoted by $\mathcal{N}(\bfw^*, \delta) := \{ \bfw \in \mathbb{R}^n \ | \ \| \bfw - \bfw^* \| \le \delta \}$, where ``$:=$'' means ``define''.
We let
$
I
$ denote the identity matrix of appropriate dimension. $\mathbb{N}$ (resp. $\mathbb{N}^+$) denotes
 the set of all natural (resp. positive natural) numbers. 
 For convenience, we sometimes use the shorthand symbol $\bfu := [\bfw; \bfxi]$ (similarly, $\bfu^*:= [\bfw^*; \bfxi^*]$).

Let $[n]$ denote the set of indices $\{1, \ldots, n\}$.
For a subset $T \subset [n]$,  
$|T|$ denotes the number of elements in $T$ (cardinality of $T$) and
 $\bfw_{T}$ denotes the subvector of
$\bfw$ indexed by $T$. 
We also let $\overline{T}$ denote the set of indices not in $T$ (i.e., $\overline{T} = [n] \setminus T$).
Given $\Gamma \in [m]$ and $A \in \mathbb{R}^{m \times n}$, $A_{\Gamma, T}$ denotes a submatrix of $A$ with row and column indexed by $\Gamma$ and $T$ respectively. Particularly, $A_{\Gamma:}$ (resp. $A_{:T}$) is the submatrix with full column (resp. row) index.

We recall from \cite[Definition 1.22]{rockafellar1976augmented} that the Moreau envelop
for a proper and lower semi-continuous function $f: \Re^n \mapsto \Re$ with $\lambda > 0$ is defined as
\begin{align*}
	&\Phi_{\lambda f(\cdot)} (\bfxi) : = \min_{ \bfq \in \Re^n } \  f(\bfq) + \frac{1}{2\lambda} \| \bfq - \bfxi \|^2 . %\\
%	&{\rm Prox}_{\lambda f(\cdot)} (\bfxi) := \argmin_{ \bfq } \frac{1}{2\beta} \| \bfq - \bfxi \|^2 + f(\bfq).
\end{align*}
The set of the solutions achieving the value $\Phi_{\lambda f(\cdot)} (\bfxi)$ is denoted by
${\rm Prox}_{\lambda f(\cdot)} (\bfxi)$ (the proximal operator of $f$).
Throughout the paper, we only deal with functions whose Moreau envelop is always achieved.

%%%%%%%%%%%%%%%
\subsection{Projection onto the $s$-sparse set}

The orthogonal projection onto the $s$-sparse set $\mathbb{S}$ is known, see \cite[Sect.~6.8.3]{beck2017first}.
We use a different (but equivalent) description below. For a given $\bfw \in \Re^n$, 
let $|\bfw|$ be the vector whose element is the absolute value of the corresponding element in $\bfw$.
Let $|\bfw|_{(i)}$ denote the $i$th largest value in $|\bfw|$. Define $\T_s(\bfw)$ to be the collection of all sets,
each consisting the $s$ indices which give rise to the largest $s$  elements in $|\bfw|$:
\be \label{Set-T}
 \T_s(\bfw) := \left\{
  \{ i_1, \ldots, i_s\} \ | \
  |w_{i_1}| = |\bfw|_{(1)}, \ldots, |w_{i_s}| = |\bfw|_{(s)}
 \right\}
\ee 
For example, for $\bfw = [10; 20; 10]$, we have
$
  \T_2(\bfw) = \left\{
  \{ 2, 1\}, \ \{2, 3\}
  \right\}.
$
%
%\[
%  \I_s^> (\bfw) := \left\{
%   i \in [n] \ | \ |w_i| > |\bfw|_{(s)}
%  \right\} \quad \mbox{and} \quad
%  \I_s^= (\bfw) := \left\{
%  i \in [n] \ | \ |w_i| = |\bfw|_{(s)}
%  \right\} .
%\]
The orthogonal projection onto $\mathbb{S}$ is given by
\[
 \Proj_{\mathbb{S}} (\bfw) =
 \left\{
   \bfq \in \Re^n \ | \ 
   \bfq = \sum_{i \in T} w_i \bfe_i, \ \ T \in \T_s(\bfw)
 \right\} ,
\]
where $\bfe_i$ is the $i$th standard unit vector in $\Re^n$.
An easy consequence of this description is the following result, see also \cite[Table 1]{pan2015solutions}. 

\begin{lemma} \label{Lemma-Proj}
	(Fixed-point characterization of the $s$-sparse set)
	Given vectors $\bfw, \bfq \in \mathbb{R}^n$ and $\alpha > 0$, we have
\begin{align} \label{proj_explicit}
	\bfw \in \Proj_{\mathbb{S} } (\bfw - \alpha \bfq) \Longleftrightarrow \left\{ 
	\begin{aligned}
		& q_i = 0, \ & \ \mbox{if} \ w_i \neq 0, \\
		& | q_i | \leq |\bfw|_{(s)}/\alpha, \ & \ \mbox{if} \ w_i = 0, 
	\end{aligned} 
	\right.
\end{align}	
Moreover, for such pair $(\bfw, \bfq)$,  the complementarity condition holds:
\[
  w_i \times q_i = 0, \quad \forall \ i \in [n].
\]
\end{lemma}

%%%%
\subsection{Positive hard-thresholding operator}

For the ease of description, we define the $(0,1)$-norm of $\bfxi \in \Re^m$ by $J: \Re^m \mapsto \Re$:
\[
  J(\bfxi) := \| \bfxi\|_{(0,1)} = \sum_{i=1}^n h(\xi_i),
\]
where $h(t)$ is the $0/1$-loss. It is not a real norm.
The proximal operator of the $h(t)$ has a simple characterization (it can be computed directly through its definition) for $\beta>0$:
\be \label{H-holding}
  \Prox_{\beta h(\cdot)} (t) = \H_{\sqrt{2\beta}} (t), \quad
  \mbox{where} \quad
  \H_{\nu}(t) := \left\{
   \begin{array}{ll}
   	 \left\{ \min \{0, \; t \} \right\},  & \mbox{if} \ t < \nu, \\
   	 \{ t \} ,               & \mbox{if} \ t > \nu, \\
   	 \{0, \; t\} ,     & \mbox{if} \ t = \nu . \\
   \end{array} 
  \right .
\ee
The operator $\H_\nu(t)$ treats small positive values $t$ as zero and is very similar to the well-known hard-thresholding operator that treats small absolute values of $t$ as zero, see \cite[Example 6.10]{beck2017first}. 
We call $\H_\nu(t)$ the 
{\em positive hard-thresholding operator}.
Consequently, the proximal operator of $J(\cdot)$ is given by
\be \label{Prox-J}
 \Prox_{\beta J(\cdot)} (\bfxi) =
 \H_{\sqrt{2\beta}} (\xi_1) \times \cdots \times \H_{\sqrt{2\beta}} (\xi_n) .
\ee 
It follows from (\ref{H-holding}) that $t \not\in \H_\nu(t)$ whenever $t \in (0, \nu)$.
Consequently, we have
\[
  \bfxi \not\in \Prox_{\beta J(\cdot)} (\bfxi) 
  \quad \mbox{if and only if there exists an index} \
  i \in [m] \ \mbox{such that} \ \xi_i \in (0, \sqrt{2 \beta}).
\]
Equivalently, we have
\be \label{Pxi}
 \bfxi \in \Prox_{\beta J(\cdot)} (\bfxi)  \quad 
 \mbox{if and only if} \ \ \xi_i \in (-\infty, 0] \cup [\sqrt{2 \beta}, \; \infty )
 \ \mbox{for all} \ i \in [m] .
\ee
We extend this result to a more general situation and it will be used in characterizing the stationary point
of our problem (\ref{HM-SSVM}) .

\begin{lemma} \label{Lemma-Prox}
	(Fix-point characterization of the hard-margin loss)
Suppose $\beta, \lambda$ are two positive constants. 
Let
$\bfxi, \bfv \in \Re^m$ be given. It holds that
\[
\bfxi \in  {\rm Prox}_{\beta \lambda J(\cdot)} (\bfxi + \beta \bfv )
\]
if and only if
\[
  \bfxi \in {\rm Prox}_{\beta \lambda J(\cdot)} (\bfxi)
  \quad \mbox{and} \quad
  \left\{
   \begin{array}{ll}
   	v_i = 0 , & \mbox{if} \ \xi_i \in (-\infty, 0) \cup [\sqrt{2 \beta \lambda}, \; \infty ) \\  [0.2ex]
   	v_i \in [0, \; \sqrt{2\lambda / \beta } ],   & \mbox{if} \ \xi_i = 0.
   \end{array} 
  \right .
\]
Consequently, the complementarity condition holds for such pair $(\bfxi, \bfv)$:
\[
  \xi_i \times v_i = 0, \ \forall \ i \in [m] .
\]
\end{lemma}

%%%%%%%%%%%%%%%%%%%%%%%%%%
\section{Stationarity Characterization of Local Minimizers} \label{Section-Optimality}

%%%%%%%%%%%%%%%%%%%%
%\subsection{Problem reformulation and its stationarity}

For the sake of simplicity, it is without loss of generality that we merge the variable $b$ into $\bfw$ in (\ref{HM-SSVM}): $\bfw := [\bfw; b]$ (Matlab notation). Define the corresponding matrix $A$ with its $i$th row being
$A_{i:} = -y_i [\bfx_i^\top, 1]$, $i=1, \ldots, m$. 
We still treat thus defined vector $\bfw$ as $n$-dimensional vector (to save us from using $(n+1)$) 
and $A$ is $m \times n$ data matrix. $\bfone$ is a vector with appropriate dimension and all entries being one. Problem (\ref{HM-SSVM}) then becomes
\be \label{HM-J}
  \min f(\bfw) := \frac 12 \| \bfw\|^2 + \lambda J(A\bfw + \bfone) , \quad \mbox{s.t.} \ \
  \| \bfw \|_0 \le s .
\ee
By introducing the auxiliary variable $\bfxi \in \Re^m$, we consider the following 
reformulation:
\be \label{Constrained-HM}
\min \frac 12 \| \bfw\|^2 + \lambda J( \bfxi) + \delta_{\mathbb{S} } (\bfw) , \quad \mbox{s.t.} \ \
A\bfw + \bfone = \bfxi,
\ee 
where $\delta_{\mathbb{S}}(\cdot)$ is the indicator function of the set $\mathbb{S}$.
The augmented Lagrangian function of (\ref{Constrained-HM}) is 
\[
\mathcal{L}_\rho(\bfw, \bfxi,\bfz) := 
\frac{1}{2}\| \bfw \|^2  + \langle \bfz, A\bfw + \bfone - \bfxi \rangle + \frac{\rho}{2} \| A\bfw + \bfone - \bfxi \|^2 + \lambda J( \bfxi) + \delta_{\mathbb{S}}(\bfw),
\]
where $\bfz \in \Re^m$ is the Lagrange multiplier and $\rho>0$ is a penalty parameter.
We will interchangeably refer to (\ref{HM-J}) and (\ref{Constrained-HM}) depending on the situation whether 
$\bfxi$ is needed or not.

\begin{definition}
	A point $\bfw^*$ is called a {\rm P}-stationary point of Problem (\ref{HM-J}) if
	there exist a Lagrange multiplier $\bfz^*$ and two positive constants $\alpha >0$ and
	$\beta >0$ such that
	\begin{equation} \label{P-stat}
		\left\{ \begin{aligned}
			& \bfw^* \in \Proj_{\mathbb{S}} (\bfw^* - \alpha ( \bfw^* + A^\top \bfz^* )),\\
			& \bfxi^* \in \Prox_{\beta \lambda J(\cdot)  }(\bfxi^* + \beta \bfz^*), \\
			& A\bfw^* + \bfone - \bfxi^* = 0,
		\end{aligned}  \right . .
	\end{equation} 
We also say that $(\bfw^*, \bfxi^* = A \bfw^* + \bfone_m)$ is a {\rm P}-stationary point of Problem (\ref{Constrained-HM}) with the
Lagrange multiplier $\bfz^*$.
\end{definition}

\begin{Remark} \label{rem_Psta}
The notation of ${\rm P}$-stationarity has its reference to the projection and proximal operators used in its 
definition.
The first inclusion relationship in (\ref{P-stat}) characterizes the stationarity with regarding to the $s$-sparse set $\mathbb{S}$. The projection operator is actually the proximal operator of the indicator function $\delta_{\mathbb{S}}(\cdot)$.
The second inclusion relationship is about the hard-margin loss function. Proximal operators have been used to
characterize stationary points in sparse optimization, see \cite{beck2013sparsity, zhou2021global}.
We also note that 
if the ${\rm P}$-stationary condition (\ref{P-stat}) is satisfied for some $\alpha = \alpha_0$ and
$\beta = \beta_0$, then it is also satisfied with any $\alpha \le \alpha_0$ and $\beta \le \beta_0$.
This follows from the fixed-point characterizations in Lemmas~\ref{Lemma-Proj} and \ref{Lemma-Prox}.
Therefore, the stationarity can be searched over an interval $\alpha \in (0, \alpha_0]$ and $\beta \in (0, \beta_0]$
even $\alpha_0$ and $\beta_0$ are often unknown in practice. 
%Our algorithm will automatically find such constants using the problem data.
In fact, the ${\rm P}$-stationarity is quite strong.
As we show below, the ${\rm P}$-stationary point and the strict local minimizer of Problem (\ref{HM-J}) has
one-to-one correspondence.

If we denote $\mathcal{S}^*:= \{ i \in [n]: w^*_i \neq 0 \}$ and $\Gamma^*:= \{ i \in [m]: \xi^*_i \neq 0 \}$, then by using Lemmas \ref{Lemma-Proj} and \ref{Lemma-Prox}, we can derive $\bfz^*_{\Gamma^*} = 0$ and $(\bfw^*+ A^\top \bfz^*)_{\mathcal{S}^*} = 0$ from \eqref{P-stat}, which further leads to
\begin{equation*}
	\bfw^*_{\mathcal{S}^*} = -A^\top_{\overline{\Gamma}^*,\mathcal{S}^*} \bfz^*_{\overline{\Gamma}^*} \quad \mbox{and} \quad \bfw^*_{\overline{\mathcal{S}}^*} = 0.
\end{equation*}
This means that $\overline{\Gamma}^*$ actually includes all the support vectors of $\bfz^*$.
\end{Remark} 

\begin{theorem} \label{Thm-Stationarity}
	(Stationarity characterization of local minimizers)
Suppose $\bfw^*$ is a local minimizer of (\ref{HM-J}). Then $\bfw^*$ is a  {\rm P}-stationary point of
(\ref{HM-J}). 
Conversely, if $\bfw^*$ is a {\rm P}-stationary point of
(\ref{HM-J}), then it must be a strict local minimizer of (\ref{HM-J}).
Moreover, there exist constants $c_*$ and $\epsilon_*>0$ such that
\be \label{Quadratic-Growth}
  f(\bfw) 
  \ge f(\bfw^*) + c_* \Big(
   \| \bfw - \bfw^*\|^2 + \| A( \bfw - \bfw^*) \|^2
  \Big), \quad \forall \ \bfw \in \N(\bfw^*, \epsilon_*) \cap \mathbb{S} .
\ee
\end{theorem} 

\begin{Remark}
In optimization, the inequality (\ref{Quadratic-Growth}) is known as the quadratic growth condition.
The objection function $f(\bfw)$ involves the hard-margin loss, which has a combinatorial nature. 
This complicates the proof.
In Appendix, we will characterize the {\rm P}-stationary point in terms of a smooth optimization problem
and eventually establish this quadratic growth condition.
It plays a very important role in convergence analysis.
We note that it is satisfied over the $s$-sparse set $\mathbb{S}$. 
Our algorithm will guarantee that all iterates $\bfw^k$ will stay in $\mathbb{S}$.
\end{Remark}

\begin{Remark}
This characterization also justifies the proposal of Mangasarian for separating data in
the $(\bfw, b)$ (feature-intercept) space. 
In this space, the distance between separating planes is strongly convex in 
$(\bfw, b)$.
Without this strong convexity, we would need extra conditions for the quadratic growth condition.
Consequently, we would not be able to establish the one-to-one correspondence between
{\rm P}-stationary points and local minimizers.
\end{Remark}

%%%%%%%%%%%%%%%%%%%%%%%%%%%%%%%%%%%%%%%%%%%%%%%%%%%%%%%%%%%%%
\section{Inexact Proximal Augmented Lagrangian Method} \label{Section-PAL}

As mentioned in Introduction, Problem (\ref{HM-SSVM}) can be put in the framework of composite optimization.
Therefore, general principle for developing augmented Lagrangian methods (ALM) set in \cite{bolte2018nonconvex}
serves a guidance for us. In this part, we develop an implementable ALM, which is based on the following important innovations.

\begin{itemize}
	\item[(i)] The subproblems of our ALM are solved inexactly.
	Computable stopping criteria are designed and are sufficient for the generated sequence to have both
	global and local linear convergence rate. This is the most challenging part of our method.
	
	\item[(ii)] In general, ALM generates infeasible iterates. Our problem has two constraints:
	\[
	    \bfw \in \mathbb{S} \quad \mbox{and} \quad A \bfw + \bfone = \bfxi .
	\]
	We treat the first constraint as ``hard'' constraint, which mus be met. In other words, we will
	generate feasible iterates $\bfw^k \in \mathbb{S}$. However, we allow the second constraint be only approximately
	satisfied. 
	This gives us much freedom to control the quality of the iterates that satisfy some decrease condition.
	
	\item[(iii)] We take the advantage of the combinatorial nature of the hard-margin loss function to define a subspace sufficiently big enough to contain a local minimizer of Problem (\ref{Constrained-HM}). 
	This subspace is potentially much smaller than the full space at each iteration.
	The benefit is that the ALM subproblems can be efficiently solved by Newton's method.  
	
\end{itemize}

The consideration above results in a new ALM.
We first describe the framework of the ALM and then state its convergence properties.

%%%%%%%%%%%%%%%%%%%%%%%%%%%%%%%%%%%%%%%%%%%%%%%%%%%%%%%%
\subsection{Framework of iPAL.} \label{frame-iPAL}

Throughout, we denote $\bfu := [\bfw; \bfxi] \in \Re^{n+m}$ and $\bfu^k := [\bfw^k; \bfxi^k]$ for each iterate.
We further define the Lyapunov function $\mathcal{M}_{\rho,\mu}: \mathbb{R}^{n+m} \times \mathbb{R}^n \times \mathbb{R}^{n} \to \mathbb{R} $ by
\begin{align*}
& \mathcal{M}_{\rho,\mu} (\bfu, \bfz, \bfv) :=\mathcal{L}_\rho (\bfu, \bfz) + \frac{\mu}{2} \| \bfw - \bfv \|^2 \\ 
	=& \underbrace{\frac{1}{2}\| \bfw \|^2 + \langle \bfz, A\bfw + \bfone - \bfxi \rangle + \frac{\rho}{2} \| A \bfw + \bfone - \bfxi \|^2 + \frac{\mu}{2} \| \bfw - \bfv \|^2}_{=: g(\bfu, \bfz, \bfv)} 
	+ \delta_{\mathbb{S}}(\bfw) + \lambda J( \bfxi),
\end{align*}  
where $\bfz$ represents the Lagrangian multiplier and $\bfv$ is a point that acts as a proximal to $\bfw$.
The function $g(\bfu, \bfz, \bfv)$ is the smooth part of the Lyapunov function.

Suppose the current iterate is $(\bfu^k, \bfz^k)$. We obtain $\bfu^{k+1}$ by
\begin{align} \label{subproblem}
	\bfu^{k+1} \approx \mathop{\arg\min}\limits_{\bfu} \mathcal{M}_{\rho,\mu} (\bfu, \bfz^k, \bfw^k) = \mathop{\arg\min}\limits_{\bfu} \underbrace{g(\bfu, \bfz^k, \bfw^k)}_{ := g_k(\bfu) }
	 + \delta_{\mathbb{S}}(\bfw) + \lambda J( \bfxi),
\end{align} 
and the Lagrange multiplier is updated according to the usual rule.
The question now is how accurate $\bfu^{k+1}$ should be calculated. We must come up with a reasonable and computable 
criterion for it. Suppose Problem (\ref{subproblem}) were to be solved exactly and let $\widehat{\bfu}^{k+1}$
denote its solution. 
Then it must satisfy the following first-order optimality condition for some $\alpha>0$ and $\beta>0$:
\be \label{P-stat-subk}
	\left\{ 
	\begin{array}{l}
		 \hbw^{k+1} \in \proj (\hbw^{k+1} - \alpha \nabla_\bfw g_k (\widehat{\bfu}^{k+1}) )) \\
		 \widehat{\bfxi}^{k+1} \in \Prox_{\beta\lambda J(\cdot) } (\widehat{\bfxi}^{k+1} - \beta \nabla_\bfxi g_k (\widehat{\bfu}^{k+1})) 
	\end{array}  
	\right.
\ee
Both the projection and the proximal operators in (\ref{P-stat-subk}) have been well studied in
Lemmas \ref{Lemma-Proj} and \ref{Lemma-Prox}, where the complementarity relationships show the different magnitudes of
the quantities involved. Let us expand those quantities in order to derive a good approximation to
(\ref{P-stat-subk}). 

Given a point $\bfu$, let us define its gradient step by
\[
  \widetilde{\bfw}^k(\bfu) :=  \bfw - \alpha \nabla_\bfw g_k (\bfu)  \quad \mbox{and} \quad
  \widetilde{\bfxi}^k(\bfu) := \bfxi - \beta \nabla_\bfxi g_k (\bfu) .
\]
Pick the index sets $T_{\bfu}$ and $\Gamma_{\bfu}$ respectively by
\[
  T_{\bfu} \in \T_s (\widetilde{\bfw}^k(\bfu)) \quad \mbox{and} \quad
  \Gamma_{\bfu} = \{
   i \in [m] \ | \ [\widetilde{\bfxi}^k(\bfu) ]_i \in (-\infty, 0] \cup [\sqrt{2 \beta\lambda }, \infty)
  \} ,
\]
where $\T_s$ is defined in (\ref{Set-T}). 
We simply use $T$ and $\Gamma$ instead of $T_{\bfu}$ and $\Gamma_{\bfu}$ when no confusion is caused.
Using Lemmas \ref{Lemma-Proj} and \ref{Lemma-Prox}, we see that
(\ref{P-stat-subk}) holds if and only if
\[
 \R_1(\bfu^{k+1}) = 0, \quad \R_2(\bfu^{k+1}) = 0, \quad \mbox{and} \quad \R_3(\bfu^{k+1}) = 0,
\]
where
\begin{align*}
	\left\{ 
	\begin{aligned}
		& \R_1(\bfu):= \| [\nabla_{T} g_k (\bfu); \bfw_{\overline{T}}]  \| \\
		& \R_2(\bfu):= \| [\nabla_{\Gamma} g_k (\bfu); \bfxi_{\overline{\Gamma}}] \| \\
		& \R_3 (\bfu) := (\beta/2) \| \nabla_\bfxi g_k (\bfu) \|^2 + \lambda J(\bfxi) - \Phi_{\beta\lambda J(\cdot) } ( \bfxi - \beta \nabla_\bfxi g_k(\bfu) ) 
	\end{aligned}\right.
\end{align*}
where $\nabla_{T} g_k (\bfu) := [ \nabla_\bfw g_k(\bfu) ]_T$ and $\nabla_{\Gamma} g_k (\bfu) := [ \nabla_\bfxi g_k(\bfu) ]_\Gamma$. We note that the residual $\R_3$ involves the Maureau envelop of the hard-margin loss $J(\bfxi)$ and plays an important role in our analysis. We now present our inexact ALM in Alg.~\ref{Alg-iPAL}.

\begin{algorithm}[H]
	\caption{(iPAL: inexact Proximal Augmented Lagrangian Method)} \label{Alg-iPAL}
	\begin{algorithmic}
		\STATE{\textbf{Initialization:}} Given positive constants $ c_1,  c_2$
		and initial point $(\bfu^0, \; \bfz^0)$. Select a positive sequence $\{ \vartheta_k \}_{k \in \mathbb{N}}$ converging to zero.
		
		\FOR{$k = 0,1, \cdots$}
		\STATE{\textbf{1. Primal step:}} 
		Starting with $(\bfu^k, \bfz^k)$, solve the
		subproblem \eqref{subproblem} for
		$\bfu^{k+1}$, which satisfies
		the following criteria:
		\be \label{error-metric}
		\left\{
		\begin{aligned}
			& \M_{\rho,\mu} (\bfu^{k+1}, \bfz^{k}, \bfw^k ) \leq \M_{\rho,\mu} (\bfu^{k}, \bfz^{k}, \bfw^k ) \ \mbox{and} \ \| \bfw^{k+1} \|_0 \leq s  \\
			&\R_1( \bfu^{k+1}  ) \le  c_1 \| \bfw^{k+1} - \bfw^k \|, \\ 
			&\R_2( \bfu^{k+1}  ) \le  c_2 \| \bfw^{k+1} - \bfw^k \|^2, \\ 
			&\R_3( \bfu^{k+1}  ) \le \vartheta_k.
		\end{aligned}
		\right.
		\ee 
		
		\STATE{\textbf{2. Multiplier step:}}
		\begin{equation} \label{Multiplier-update}
			\bfz^{k+1} = \bfz^k + \rho(A \bfw^{k+1} + \bfone - \bfxi^{k+1}).
		\end{equation}
		
		\ENDFOR
	\end{algorithmic}
\end{algorithm}

\begin{Remark}
The algorithm iPAL follows the standard framework of ALM having both the primal and the multiplier steps.
The only difference is that the subproblem was solved inexactly, but increasingly accurate.
In particular, the residual $\R_2$ is one order more accurate than $\R_1$ is. 
This requirement is crucial in ensuring the generated sequence to converge linearly. 
We will design Newton's method for the subproblem in the next section to meet those criteria.
For now, we present the convergence results.
\end{Remark} 

%%%%%%%%%%%%%%%%%%%%%%%%%%%%%%%%%%%%%%%%%%%%%%%%%%
\subsection{Convergence of iPAL} \label{Subsection-Global}

As rightly emphasized in \cite{bolte2018nonconvex}, certain regularity is needed on the constraints in composite optimization for global convergence of ALMs. We need the following regularity assumption. 
Let $r := \lfloor s/2 \rfloor$ and define $\Theta:= \{ T \subseteq [n]: |T| = r \}$, where $\lfloor \cdot \rfloor$ is the floor function. 

\begin{assumption}\label{assum}
	For any $T \in \Theta$, $A_{:T}$ has full row rank. 
	Consequently, there exists $\gamma > 0$ satisfying $\gamma^2 = \min_{T \in \Theta} \lambda_{\min} (A_{:,T}A_{:,T}^\top)$.
\end{assumption}

The assumption is particularly useful when the sample data is small (i.e., $m \ll n$). This has been confirmed
in our numerical experiments for such data. The assumption can be weakened to only those rows of $A$ indexed by
$\oG_k$ in Alg.~\ref{Alg-iPAL}. A further result (see \eqref{y+0-iden}) indicates that $\oG_k$ can be seen as an approximation of the support vector index set $\oG^*$ (defined in Remark \ref{rem_Psta}), which is usually much smaller than $m$. This increases the
chance for the assumption to hold. The general assumption significantly simplifies our analysis.\\

\noindent 
\textbf{Parameter Setup:} Let $c_1$ and $c_2$ be two constants used in Alg.~\ref{Alg-iPAL}. Given $\mu > 0$, 
set $\rho$ and $\eta$ as follows:
\begin{align} \label{para_set}	\rho \geq \left\{ \frac{2}{\gamma^2}, \frac{8( c_3^2 +  c_4^2)}{\mu} \right\}, \quad
	 \eta = \frac{4 c_4^2}{\rho},
\end{align}
where $ c_3:= ( 2c_1 + \mu + 2 )/\gamma$ and $ c_4:= ( 2c_1 + \mu )/\gamma$.

%%%%%%%%%%%%%%%%%%%%%%%%%%%%%%%%%%%%
Our first result states that
Alg.~\ref{Alg-iPAL} leads to a sufficient decrease in the function value of the Lyapunov function
$\M_{\rho,\eta}(\cdot)$. Let
\[
  \M_{k+1} := \M_{\rho,\eta} (\bfu^{k+1}, \bfz^{k+1},  \bfw^k), \ \ \mbox{for} \ k=0, 1, \ldots, .
\]

\begin{proposition} \label{Prop-Sufficient-Decrease} 
	Suppose that Assumption \ref{assum} holds and parameters are chosen as in \eqref{para_set}. If $\{ (\bfu^{k} ; \bfz^k ) \}_{k \in \mathbb{N}}$ is a sequence generated by iPAL. The following hold.
	
	\begin{itemize} 
		\item[(i)] (Sufficient decrease) The sequence $\{ \M_k\}$ is nonincreasing and
	\begin{align} \label{lya-des}
		\M_k - \M_{k+1} \geq \frac{\mu}{4} \| \bfu^{k+1} - \bfu^k \|^2.
	\end{align}

       \item[(ii)] (Sequence boundedness) The sequence $\{ (\bfu^{k} ; \bfz^k ) \}$ is bounded. Moreover
       	\begin{align} \label{suc_change}
       	\lim_{k \to \infty} \| \bfu^{k+1} - \bfu^k \| = 0\quad \mbox{and} \quad \lim_{k \to \infty} \| \bfz^{k+1} - \bfz^k \| = 0.
       \end{align}
   
   \end{itemize} 
\end{proposition}

\begin{Remark}
	If the Lyapunov function $\M_{\rho,\eta}(\bfu, \bfz, \bfv)$ is bounded from below by a constant $M_\infty$, then (\ref{lya-des}) would 
	imply
	\[
	  \frac{\mu}{4} \sum_{k} \| \bfu^{k+1} - \bfu^k \|^2 \le \sum_{k} (\M_k - \M_{k+1}) \le  \M_1 - M_\infty \le \infty.
	\]
	Then \eqref{suc_change} would be a direct consequence.
%	This means that the sequence $\{ \bfu^k\}$ has a finite length property and hence is Cauchy sequence, which means
%	that $\{ \bfu^k\}$ converges, see \cite[Thm.1]{bolte2014proximal}. The boundedness in (ii) would be a direct consequence. However, 
%	we do not know if the Lyapunov function is bounded from below on the generated sequence. Therefore, we
%	need a different proof for the boundedness result in (ii).	
\end{Remark}

Those results ensures the global convergence as well as linear convergence rate of iPAL.

\begin{theorem} (Global Convergence) \label{Thm-Global}
	Suppose that Assumption \ref{assum} holds and parameters are chosen as \eqref{para_set}. Let $\{ (\bfu^{k} ; \bfz^k ) \}_{k \in \mathbb{N}}$ be a sequence generated by iPAL.
	Then the whole sequence converges to a {\rm P}-stationary pair $(\bfu^*, \bfz^*)$ of (\ref{Constrained-HM}). 
	Furthermore, $\bfu^*$ is a strict minimizer of (\ref{Constrained-HM}).
\end{theorem}

Since the whole sequence $\{ \bfu^k; \bfz^k\}$ converges and the Lyapunov sequence $\{ \M_k\}$ is nonincreasing,
there must exist a limit, denoted by $\M_*$. Actually, we can prove $\M_* = \M_{\rho,\eta}(\bfu^*, \bfz^*, \bfw^*)$. For more details, please refer to Corollary \ref{col_convergence} in Appendix.
%:
%\[
%  \M_* := \lim_{k \rightarrow \infty} \M_k =\lim_{k \rightarrow \infty} \M_{\rho,\eta}(\bfu^k, \bfz^k, \bfv^k) 
%  \ge \M_{\rho,\eta}(\bfu^*, \bfz^*, \bfv^*) ,
%\]
%where the last inequality used the lower-semicontinuity of $\M_{\rho,\eta}(\bfu, \bfz, \bfv)$.

\begin{theorem} \label{Thm-Convergence-Rate}
	(Linear rate of convergence)
Under the premise in Theorem \ref{Thm-Global}, the following estimations hold with a constant $q \in (0,1)$.
	
	\begin{itemize}
	 \item[(i)] (Linear convergence in Lyapunov function) 
	 There exists a positive constant $c_m$ and a sufficiently large index $k^*$ such that 
	\begin{align} \label{r-linear-v}
		\M_k - \M_* \leq  c_m q^k, \quad \forall k \geq k^*.
	\end{align}

     \item[(ii)] (Linear convergence in iterative sequence)
     There exist a sufficiently large index $k^*$ and positive constants $ c_w $, $c_\xi$ and $ c_z$ such that for any $k \geq k^*$, it holds
     \begin{align} \label{r-linear-wz}
     	\| \bfw^{k} - \bfw^* \| \leq  c_w \sqrt{q}^k, \ \ \| \bfxi^{k} - \bfxi^* \| \leq  c_\xi \sqrt{q}^k ,\ \ \mbox{and}\ \ \ \| \bfz^{k} - \bfz^* \| \leq  c_z \sqrt{q}^k .
     \end{align}

\end{itemize}
\end{theorem}

%%%%%%%%%%%%%%%%%%%%%%%%%%%%%%%%%%%%%%%%%%%%%%%%%%%%%%%%
\section{Projected Gradient-Newton Method for Subproblems} \label{Section-Newton}

The algorithmic framework of iPAL looks promising in terms of its global and linear convergence.
To make it practically effective, we need to address how the subproblem (\ref{subproblem}) can be
efficiently solved so as to meet the stopping criteria (\ref{error-metric}).
As mentioned earlier, our ultimate purpose is to apply Newton's method in reduced dimensions.
However, it is widely known that Newton's method is a local method. This motivates us to 
use a gradient descent method to initialize the computation.
We put those considerations in precise formulation.

First, the subproblem (\ref{subproblem}) takes the following form:
\begin{align} \label{inner-sub}
	\min_{\bfu:= (\bfw,\bfxi)} G(\bfu) := g(\bfu) + \delta_{\mathbb{S}} (\bfw) + \lambda J( \bfxi),
\end{align} 
where we dropped the dependence of $g$ on the iterate $k$. The main purpose is to solve
(\ref{inner-sub}). It is very important to note that
(i) the gradient $\nabla g(\bfu)$ is Lipschitzian continuous with  constant $\ell_g$:
\[
  \| \nabla g(\bfu) - \nabla g(\bfv) \| \le \ell_g \| \bfu - \bfv\| \quad \forall \ \bfu , \bfv \in \Re^n
\]
and (ii) $g(\bfu)$ is strongly convex with  constant $\sigma_g$:
\[
  g(\bfu) \ge g(\bfv) + \langle \nabla g(\bfv), \; \bfu -\bfv \rangle + \frac{\sigma_g}2 \| \bfu - \bfv\|^2
  \quad \forall \ \bfu , \bfv \in \Re^n .
\]

Now suppose $\bfu^j = (\bfw^j, \bfxi^j)$ be the current iterate. For given two constants $\alpha >0$ and $\beta>0$
(they serve as stepsizes respectively for $\bfw$ and $\bfxi$), the new iterate by the gradient step is given by
\be \label{what}
  \widehat{\bfw}^j := \bfw^j - \alpha \nabla_{\bfw} g(\bfu^j)
  \quad \mbox{and} \quad
  \widehat{\bfxi}^j := \bfw^j - \beta \nabla_{\bfxi} g(\bfu^j) .
\ee
We then project $ \widehat{\bfw}^j$ to the $s$-sparse set $\mathbb{S}$ and compute the hard-margin proximal of
$\widehat{\bfxi}^j$ and denote them by $\bfu^{j+1/2} = (\bfw^{j+1/2}, \bfxi^{j+1/2})$
\be\label{w-half}
  \bfw^{j+1/2} = \Proj_{\mathbb{S}} ( \widehat{\bfw}^j ) \quad \mbox{and} \quad
  \bfxi^{j+1/2} \in \Prox_{\lambda\beta J(\cdot)} ( \widehat{\bfxi}^j ) .
\ee
We only consider those indices where $\bfw^{j+1/2} $ and $\bfxi^{j+1/2} $ are not zero:
\be \label{Tj}
 T_j \in \T_s ( \widehat{\bfw}^j ) \quad \mbox{and} \quad
 \Gamma_j = \left\{
  i \in [m] \ | \ [ \widehat{\bfxi}^j ]_i \in (-\infty, 0) \cup (\sqrt{2 \lambda\beta}, \infty) 
 \right\} .
\ee
Consequently, when restricting to the subspace:
\[
  \left\{
  \bfu= (\bfw, \bfxi) \in \Re^n \times \Re^m \ | \ \bfw_{\overline{T}_j } = 0, \ \
                                              \bfxi_{\overline{\Gamma}_j } = 0
  \right\},
\]
the objective function $G(\bfu)$ is locally twice continuously differentiable. Newton's method is well defined
over this subspace. The resulting algorithm is called the projected gradient-Newton method, which
is detailed in Alg.~\ref{Alg-PGN}

\begin{algorithm}[H] 
	\caption{(PGN: Projected Gradient-Newton Method)} \label{Alg-PGN}
	\begin{algorithmic}
		
		\STATE{Initialization: Set $\alpha, \beta \in (0,1/\ell_g)$, take initial point $\bfu^0 := (\bfw^0, \bfxi^0) \in \mathbb{R}^{n+m}$} with $\| \bfw^0 \|_0 \leq s$.
		\FOR{$j=0,1,\cdots$}
		
		\STATE{\textbf{1. Identification step:} } 
		Compute $\widehat{\bfw}^j$ and $\widehat{\bfxi}^j$ by (\ref{what}). 
		Select $T_{j}$ and $\Gamma_j$ by (\ref{Tj}) 
		
		\STATE{\textbf{2. Gradient step: } } Compute $\bfu^{j+1/2} = (\bfw^{j+1/2}, \bfxi^{j+1/2})$ by
		(\ref{w-half}).
		
		\STATE{\textbf{3. Newton step:} } Denote $\Upsilon_j:= T_j \cup \Gamma_j$ and compute $\tbu^{j+1} := (\widetilde{\bfw}^{j+1}, \widetilde{\bfxi}^{j+1})$ by solving the following
		reduced Newton equation in $\bfu: = (\bfw,\bfxi)$
		\begin{align}\label{subspace-newton}
			\left\{ \begin{aligned}
				& H^{j+1/2} ( \bfu - \bfu^{j+1/2} )_{\Upsilon_j} = -\nabla_{\Upsilon_j} g ( \bfu^{j+1/2} ) \\
				& \bfw_{\oT_j} = 0, \ \bfxi_{\oG_j} = 0, 
			\end{aligned} \right.
		\end{align}
		where $H^{j+1/2}:= [ \nabla^2 g(\bfu^{j+1/2}) ]_{\Upsilon_j, \Upsilon_j}$.
		\STATE{\textbf{4. Update step:} }	Update $\bfu^j$ either by the Newton step or
		the gradient step as follows:
		\be \label{Newton-Condition}
		\bfu^{j+1} = \left\{
		\begin{array}{ll}
			\tbu^{j+1} , & \ \mbox{if} \ G ( \bfu^{j+1/2} ) - G ( \tbu^{j+1}) \geq (\sigma_g/4) \|  \bfu^{j+1/2} - \tbu^{j+1} \|^2 \\ [1ex]
			\bfu^{j+1/2} , & \ \mbox{otherwise}
		\end{array} 
		\right .
		\ee 
		
		\ENDFOR
	\end{algorithmic}
\end{algorithm}
%%%%%%%%%%%%%%%%%%%%%%%%%%%%%%%%%%%%%%%%%%%%%%%%%%%%%%%%%
\begin{Remark}
	(i) Computational complexity of the gradient step.
	Assuming the gradient of $g(\bfu)$ is available, 
	the complexity of selecting $T_j$ and $\Gamma_j$ is $O(ns)$. According to Lemmas \ref{Lemma-Proj} and
	\ref{Lemma-Prox}, the gradient update is computed by
			\begin{align} \label{gradient_xy}
					\begin{aligned}
							&\bfw^{j+1/2}_i = \left\{ \begin{aligned}
									& [ \bfw^j - \alpha \nwg (\bfu^j) ]_i, &&\ \mbox{if} \ i \in T_j \\
									& 0, &&\ {\rm otherwise}.
								\end{aligned}  \right.  \\ &\bfxi^{j+1/2}_i = \left\{ \begin{aligned}
									& [ \bfxi^j - \beta \nyxi (\bfu^j) ]_i, &&\ \mbox{if} \ i \in \Gamma_j \\
									& 0, &&\ {\rm otherwise}.
								\end{aligned} 
							\right.
						\end{aligned}
				\end{align}
	Therefore, the overall complexity for computing $\bfu^{j+1/2}$ is $O(ns)$.
	
 (ii) Computational complexity of the Newton step.
  We expand the Newton equation \eqref{subspace-newton} as follows:
 \begin{align*}
 	\left[\begin{array}{cc}
 		[(\mu + 1) I + \rho A^\top A ]_{T_j, T_j} & - \rho(A_{\Gamma_j,T_j})^\top \\
 		- \rho A_{\Gamma_j,T_j} & \rho I
 	\end{array}\right] \left[
 \begin{array}{c}
 	\bfd_w \\ \bfd_\xi
 \end{array} 
 \right] 
  = - \nabla_{\Upsilon_j} g_k(\bfu^{j+1/2}) = \left[
  \begin{array}{c}
  	\bfb_w \\ \bfb_\xi
  \end{array} 
  \right] 
 \end{align*}
 with variable $\bfd = [\bfd_w; \bfd_\xi] \in \Re^{|\Upsilon_j| }$ to be computed. 
 By using Schur complement theorem, it is equivalent to 
% \begin{align}\label{linear_sub}
\begin{equation} \label{lin_eq}
	\left\{
	\begin{array}{ll}
		\left( (\mu + 1) I + \rho A^\top_{ \oG_j,T_j } A_{ \oG_j,T_j } \right) \bfd_w 
		= \bfb_w + (A_{\Gamma_j,T_j})^\top \bfb_\xi
		,\\ [1ex]
		\bfd_\xi = \frac{1}{\rho}  \bfb_\xi   + A_{\Gamma_j,T_j} \bfd_w.
	\end{array}  \right.
\end{equation}
% \end{align}
The computational complexity for solving this linear system is $O( | \oG_j | | T_j |^2 )$. We can also apply Sherman-Morrison-Woodbury formula to  this linear equation when $ | \oG_j | \ll | T_j | $ and the corresponding computational complexity will be $O( | \oG_j |^2 | T_j | )$.
\end{Remark}

\begin{theorem}[Global Convergence of PGN] \label{Thm-PGN-Global}
	Let $\{ \bfu^j \}_{j \in \mathbb{N}}$ be the sequence produced by PGN. The following statements hold.
	
	\begin{itemize}
		\item[(i)] (Sufficient decrease)  We have
		\begin{align} \label{des-G}
			G(\bfu^j) - G(\bfu^{j+1}) \geq \zeta \| \bfu^{j+1/2} - \bfu^j \|^2 + (\sigma_g/4) \| \bfu^{j+1} - \bfu^{j+1/2} \|^2,
		\end{align} 
		where $\zeta:= \min \{ ( 1/\alpha - \ell_g )/2, ( 1/\beta - \ell_g )/2 \}$.
		 This further leads to
		\begin{align}\label{inn-suc-chan}
			\lim_{j\to \infty} \| \bfu^{j+1} - \bfu^j \| = 0\quad  \mbox{and} \quad \lim_{j \to \infty} \| \bfu^{j+1/2} - \bfu^j \| = 0
		\end{align}
		\item[(ii)] (Convergence to stationary point) 
		The sequence $\{ \bfu^j \}_{j \in \mathbb{N}}$ converges to a {\rm P}-stationary point $\widehat{\bfu}:= ( \hbw, \hbxi )$ satisfying
		\begin{align} \label{P-stat-sub}
			\left\{\begin{aligned}
				& \hbw \in \proj (\hbw - \alpha \nabla_\bfw g (\widehat{\bfu}) )) \\
				& \widehat{\bfxi} \in \Prox_{\beta\lambda J(\cdot)}(\widehat{\bfxi} - \beta \nabla_\bfxi g (\widehat{\bfu}))  
			\end{aligned}\right.
		\end{align}
	
		\item[(iii)] (iPAL is well defined) If $\bfw^0 \neq \hbw$, then there exists a sufficiently large index $j_k$ such that $\bfu^{j_k}$ satisfies the stopping
		criteria (\ref{error-metric}).
		
%		It holds
%		\begin{align} \label{y+-iden}
%			\lim_{j \to \infty} J(\bfxi^{j+1}) = \lim_{j \to \infty} J(\bfxi^{j+1/2}) = J(\bfxi^*) \ \mbox{and} \ \lim_{j \to \infty} G(\bfu^j) = G(\bfu^*).
%		\end{align}
	\end{itemize}
\end{theorem}

The global convergence theorem states that as long as the optimization method for the subproblem yields the sufficient
decrease in terms of (\ref{des-G}), then the generated sequence must converge to a {\rm P}-stationary point.
Moreover, the stopping criteria (\ref{error-metric}) is met as soon as $j \ge j_k$. Note that we use $j_k$ for 
the iterate index because we apply Alg.~\ref{Alg-PGN} to the subproblem at the $k$th iteration of iPAL. 
Given the linear convergence rate of iPAL, the smaller $j_k$ is at each iteration, the more efficient iPAL would be.
Therefore, we study when the Newton iteration takes place and whether it has a quadratic convergence. 
We consider the situation near the stationary point $\widehat{\bfu}$ in (\ref{P-stat-sub}).
It follows from Lemma~\ref{Lemma-Prox} that $\widehat{\bfxi}$ and $\nabla_\bfxi g (\widehat{\bfu})$ must satisfy
the complementarity condition. We assume further that they satisfy the strict complementarity condition:
\begin{align}\label{scc}
	\hbxi_i + [ \nabla_\bfxi g ( \hbu ) ]_i \neq 0, \quad  \forall i \in [m].
\end{align}
Under this assumption, we can prove that Newton's step is always accepted when $j \ge j_k$ and hence
PGN is quadratically convergent.

\begin{theorem}[Local Quadratic Convergence of PGN] \label{qua_convergence}
	Let $\{ \bfu^j \}_{j \in \mathbb{N}}$ be a sequence converging to a {\rm P}-stationary point $\hbu$ of \eqref{subproblem}. Suppose that $\hbxi$ and $\nabla_{\bfxi} g (\hbu)$ satisfy strictly complementary condition \eqref{scc}, then there exists sufficiently large integer $j_k$ such that %the following statements hold.
	Newton's step will always be accepted for all iterations $j \geq j_k$. Moreover,
	we have
	\[  %\]\begin{align}\label{qua-rate}
		\| \bfu^{j+1} - \hbu \| \leq O(\| \bfu^j - \hbu \|^2) \quad \mbox{for} \ \ j \ge j_k .
    \] %	\end{align}
	
\end{theorem}

This may be the best result one may hope for when Newton's method is used. 
The question now is whether the Newton equation can be efficiently solved.
Our numerical results demonstrate that it is the case for many types of data.

%%%%%%%%%%%%%%%%%%%%%%%%%%%%%%%%%%%%%%%%%%%%%%%%%%%%%%%%%
\section{Numerical Experiments} \label{Section-Numerical}

In this section, extensive numerical experiments will be conducted by using Matlab 2022a on a laptop with 32GB memory and Intel CORE i7 2.6 GHz CPU.

%%%%%%%%%%%%%%%%%%%%%%%%%%%%%%%%%%%%%%%%%%%%%%%%%%%%%%%
\subsection{Benchmark Methods and Experimental Setting} \label{bmes}

To implement iPAL, we need to set up two types of parameters. One type called model parameters of \eqref{Constrained-HM} contains $\lambda$, $\rho$, $\mu$ and $s$. To simplify the parameter tuning, we will set $\lambda = \rho$. The best choices are often dependent on data, and thus we will give more details about the selection in the subsequent experiments. Another type of parameters appearing in Alg. \ref{Alg-iPAL} and Alg. \ref{Alg-PGN} is called algorithmic parameters. We set 
\begin{align}
	c_1 =  c_2 = 0.1,~~~ \gamma = 0.1\min\{ \| \bfa_i \| | i \in [m] \},~~~\epsilon_k = \lambda/k
\end{align}
and $\eta$ is taken as \eqref{para_set}.
We adopt $(\bfw^0, \bfxi^0,\bfz^0) = \bfzero$ as initial point and iPAL will stop if the following criterion holds
\begin{align} \notag
	\frac{\| \bfw^k - \bfw^{k-1} \| + \| \bfxi^k - \bfxi^{k-1} \| + \| \bfz^k - \bfz^{k-1} \|}{\| \bfw^k \| + \| \bfxi^k \| + \| \bfz^k \|} < 10^{-3}
\end{align}

We also select four efficient SVM solvers for numerical comparison. Together with iPAL, the five algorithms designed for solving different SVM models are summarized in Table \ref{tb1}
\begin{table}[htbp]
	\centering
	\caption{Benchmark Algorithms and Their Models} \label{tb1}
	\resizebox{\textwidth}{14mm}
	{\begin{tabular}{ccccc}
		\hline
		Algorithm    & Reference   & Loss Function                     & Regularizer & Constraint \\ \hline
		iPAL  & This work  & Hard margin                    & $\ell_2$  & $\ell_0$      \\
		ADMM0/1  &\cite{wang2021support} & Hard margin                    & $\ell_2$  & --          \\
		LISVM    &\cite{yuan2010comparison}   & Hinge                         & $\ell_1$  & --          \\
		NLPSVM  &\cite{fung2004feature}  & Hinge                         & $\ell_1$  & --          \\
		PDLSVM  &\cite{shao2019joint}   & Least square primal and dual  & $\ell_1$  & Linear      \\ \hline
	\end{tabular}}
\end{table}

Four metrics are used for evaluating performance of the algorithms. They are classification accuracy: $\texttt{Acc}:= 1 - J(A\bfw) / m $, CPU time (\texttt{Time}), the number of support vectors (\texttt{nSV}), and the number of nonzero elements $\texttt{nnz}:= \| \bfw \|_0$. As LISVM adopts a coordinate descent method without introducing dual variables, this solver does not provide a dual solution and thus we do not record the \texttt{nSV} for it.

\subsection{Experiments on Simulated Data}
In this subsection, we will test all the solvers on datasets generated by the following example.

\begin{example} \label{ex1}
	Samples with positive (resp. negative) labels are drawn from 
	the normal distribution
	$N(\mu_1,\Sigma_1)$ (resp. $N(\mu_2,\Sigma_2)$), where the parameters
	$\mu_1 \in \mathbb{R}^n$ (resp. $\mu_2$) are mean vectors, and $\Sigma_1 \in \mathbb{R}^{n\times n}$ (resp. $\Sigma_2$) are diagonal covariance matrices. 
	We then flip $r$ percentage (noise ratio) of those samples, making them be marked with reverse labels.
\end{example}

\subsubsection{Convergence Test} In this part, we will observe how the model parameters ($\lambda$, $\rho$, $\mu$ and $s$) influence the convergence of iPAL. We will use the following metric to judge the violation of first-order optimality condition of \eqref{Constrained-HM} for an iterate
\begin{align} \notag
	\texttt{ VFC} := \max\{  {\rm dist}_p, \ \ {\rm dist}_d, \ \ {\rm dist}_c  \},
%	\max\{ {\rm dist} (\bfw^k , {\rm Proj}_\mathbb{S} ( \bfw^k - \alpha ( \bfw^{k} + A^\top \bfz^k ) ) ) , 
%	{\rm dist}(\bfxi^k,  {\rm Prox_{\alpha\lambda h(\cdot)}} ( \bfxi^k + \alpha \bfz^k )), \| A\bfw^k + \bfb - \bfxi^k \| \}.
\end{align}
where
\begin{eqnarray*}
{\rm dist}_p &:=& \| \bfw^k - {\rm Proj}_\mathbb{S} ( \bfw^k - \alpha ( \bfw^{k} + A^\top \bfz^k ) ) \| , \\
{\rm dist}_d &:=& \| \bfxi^k -  {\rm Prox_{\alpha\lambda h(\cdot)}} ( \bfxi^k + \alpha \bfz^k ) \|, \\
{\rm dist}_c &:=& \| A\bfw^k + \bfone - \bfxi^k \| .
\end{eqnarray*}
A simulated dataset with $m = 1000$ and $n = 2000$ is generated as the way described in Ex. \ref{ex1}. 
As mentioned at the beginning of Subsection \ref{bmes}, we will set $\lambda = \rho$ with the model parameters selected from the following sets:
\begin{equation} \notag
	\Omega_\rho = \{ 10^{-3}, 10^{-2}, \cdots, 10^{3} \},~\Omega_\mu := 10^{-2}\times\{2^0,2^1,\cdots,2^{10}  \},~\Omega_s = \{20,40,\cdots,200  \}
\end{equation}
%Changes of \texttt{VFC} and \texttt{Time} along with iteration are recorded as follows.
We have the following comments.

\begin{figure}[htbp]	
	\subfigure{
		\begin{minipage}[t]{0.5\linewidth}
			\centering
			\includegraphics[width=2.9in]{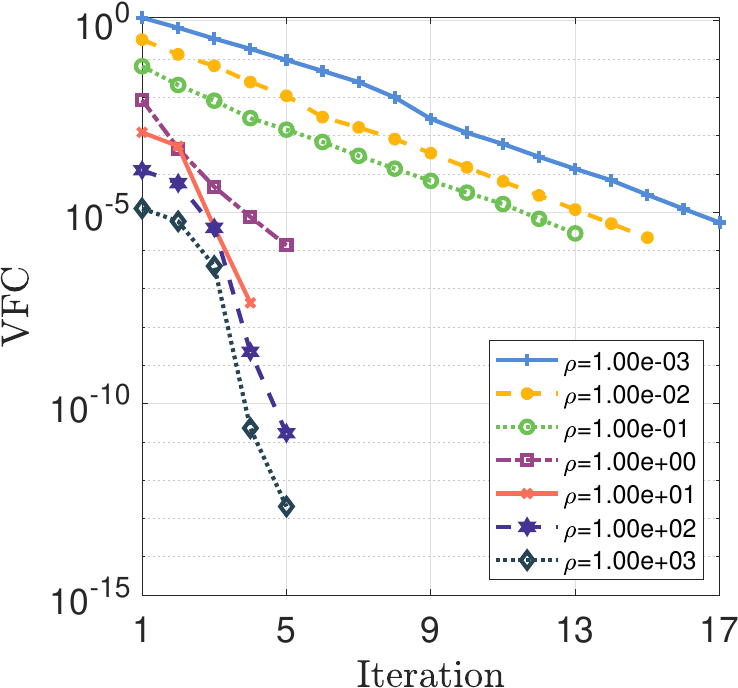}
			%\caption{fig1}
		\end{minipage}%
	}%
	\subfigure{
		\begin{minipage}[t]{0.5\linewidth}
			\centering
			\includegraphics[width=2.9in]{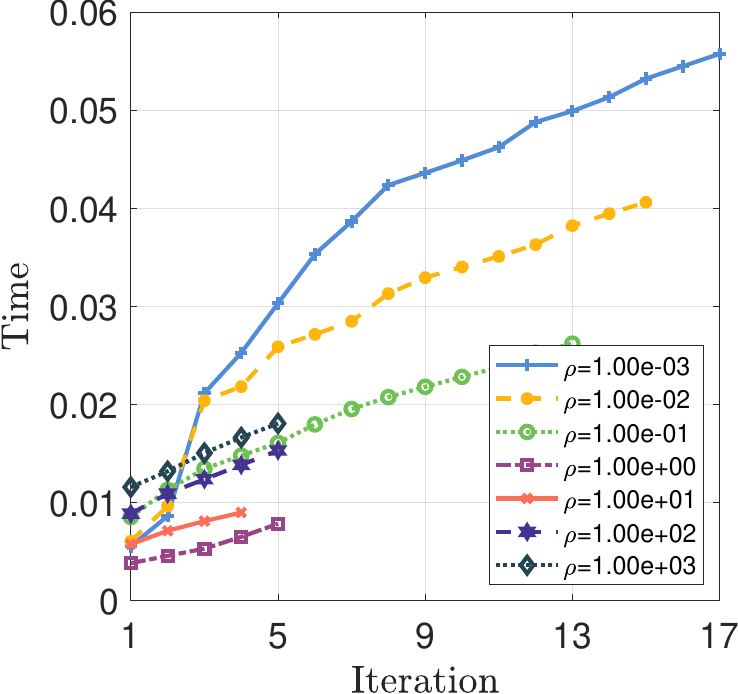}
			%\caption{fig2}
		\end{minipage}%
	}%
	\centering
	\caption{\texttt{VFC} and \texttt{Time} of iPAL along with iteration when $\mu = 10^{-2}$, $s = 20$ and $\rho = \lambda \in \Omega_\rho$.}
	\label{fig1}
\end{figure}

\begin{figure}[htbp]	
	\subfigure{
		\begin{minipage}[t]{0.5\linewidth}
			\centering
			\includegraphics[width=2.9in]{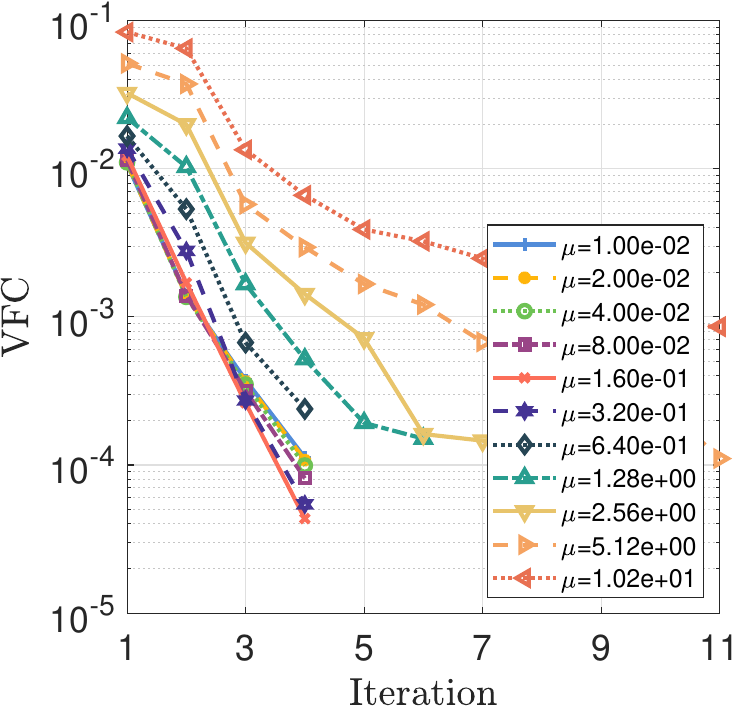}
			%\caption{fig1}
		\end{minipage}%
	}%
	\subfigure{
		\begin{minipage}[t]{0.5\linewidth}
			\centering
			\includegraphics[width=2.9in]{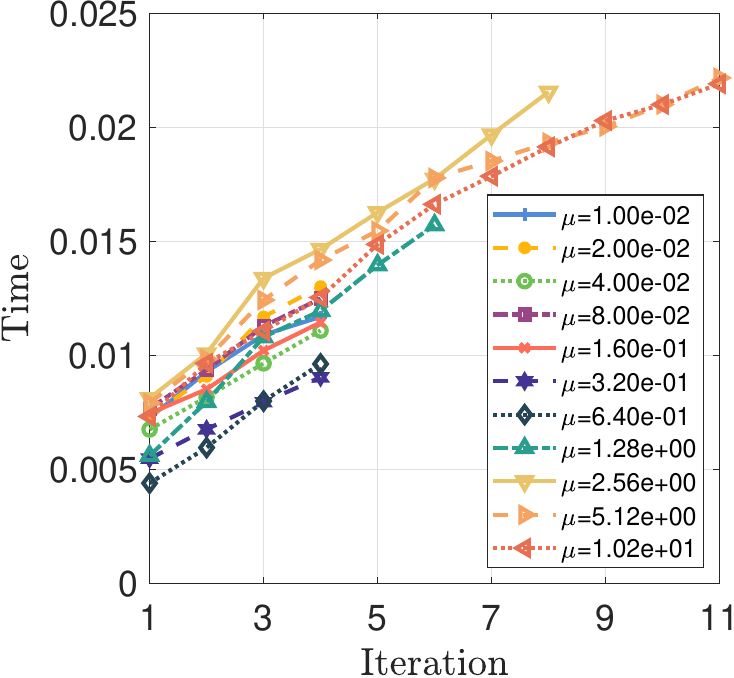}
			%\caption{fig2}
		\end{minipage}%
	}%
	\centering
	\caption{\texttt{VFC} and \texttt{Time} of iPAL along with iteration when $\rho = \lambda = 1$, $s = 20$ and $\mu \in \Omega_\mu$.}
	\label{fig2}
\end{figure}

\begin{figure}[htbp]	
	\subfigure{
		\begin{minipage}[t]{0.5\linewidth}
			\centering
			\includegraphics[width=2.9in]{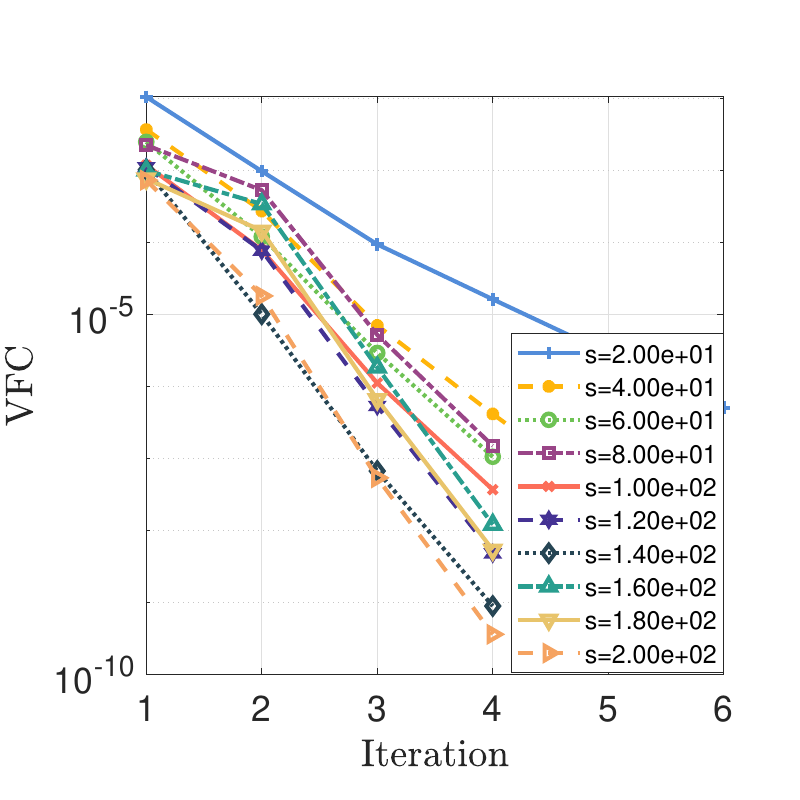}
			%\caption{fig1}
		\end{minipage}%
	}%
	\subfigure{
		\begin{minipage}[t]{0.5\linewidth}
			\centering
			\includegraphics[width=2.9in]{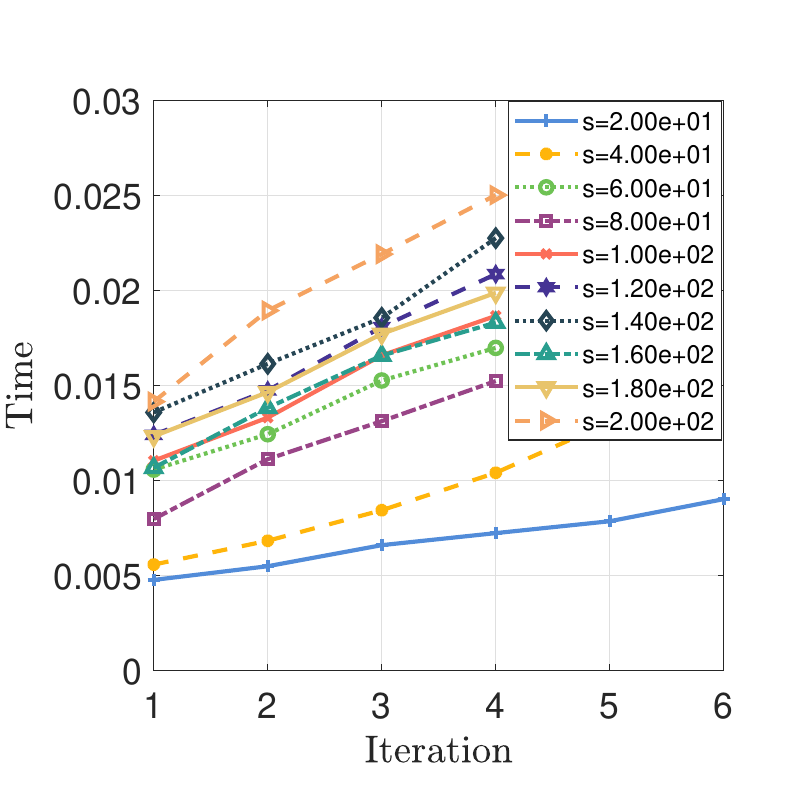}
			%\caption{fig2}
		\end{minipage}%
	}%
	\centering
	\caption{\texttt{VFC} and \texttt{Time} of iPAL along with iteration when $\mu = 10^{-2}$, $\rho = \lambda = 1$ and $s \in \Omega_s$.}
	\label{fig3}
\end{figure}

\begin{itemize}
\item[(i)] From Fig. \ref{fig1}, we can observe that \texttt{VFC} decreases faster when $\rho$ grows. However, the \texttt{Time} v.s. Iteration graph in Fig \ref{fig1} shows that a large $\rho$ does not always leads to a smaller \texttt{Time}. In fact, when $\rho$ increases, the conditional number of linear equation \eqref{lin_eq}
becomes bigger and thus it takes more time to solve.

\item[(ii)] We illustrate how the change of $\mu$ influence the convergence of iPAL. 
As shown in Fig \ref{fig2}, iPAL with larger $\mu$ tends to converge slower.
But it might spend less \texttt{Time} because the linear system \eqref{lin_eq} admits smaller conditional number. 
For example, a medium value $\mu = 0.32$ leads to the least \texttt{Time} in this simulation.

\item[(iii)]  We can see from Fig. \ref{fig3} that the convergence rate shows a faster decreasing trend when the $s$ grows. 
This is because the matrix dimension in linear system \eqref{lin_eq} is $s \times s$.
A smaller $s$ will lead to a significant reduction in dimension and computation. 
That is why iPAL with $s = 20$ (the smallest value of $s$) runs much faster than other cases 
(see \texttt{Time} v.s. Iteration in Fig. \ref{fig3} ).

\end{itemize} 

\subsubsection{Numerical Comparison } In this part, we will generate datasets with various $m$, $n$ and $r$ (noise rate) by the method in Ex. \ref{ex1}.
The performance of all the five algorithms will be compared. 
Half of the samples will be chosen as training set, and the rest of the samples are used for testing. In the following three tests, for iPAL, we set $\lambda = 1$, $\rho = 1$, $\mu = 10^{-2}$ and $s = 20$. Other algorithms used their default parameter settings. 

\begin{figure}[htbp]	
	\subfigure{
		\begin{minipage}[t]{0.5\linewidth}
			\centering
			\includegraphics[width=2.9in]{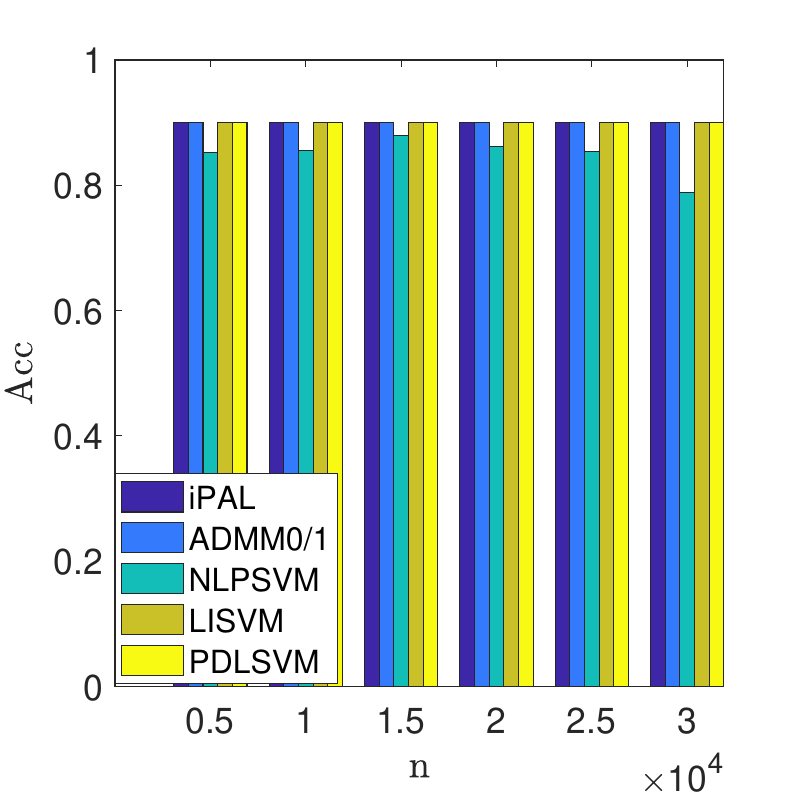}
		\end{minipage}%
	}%
	\subfigure{
		\begin{minipage}[t]{0.5\linewidth}
			\centering
			\includegraphics[width=2.9in]{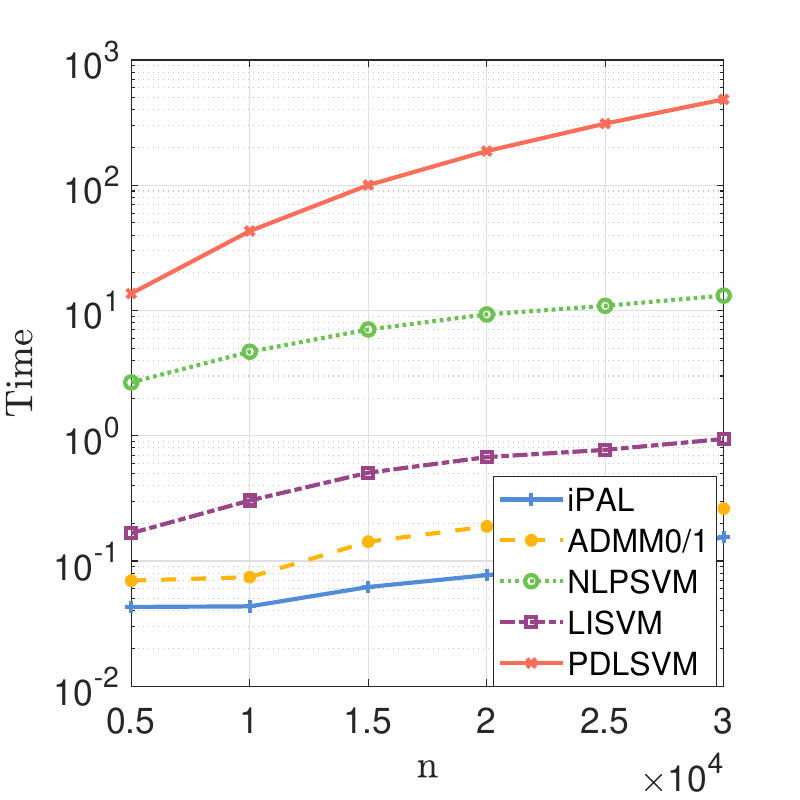}
		\end{minipage}%
	}%
	
	\subfigure{
		\begin{minipage}[t]{0.5\linewidth}
			\centering
			\includegraphics[width=2.9in]{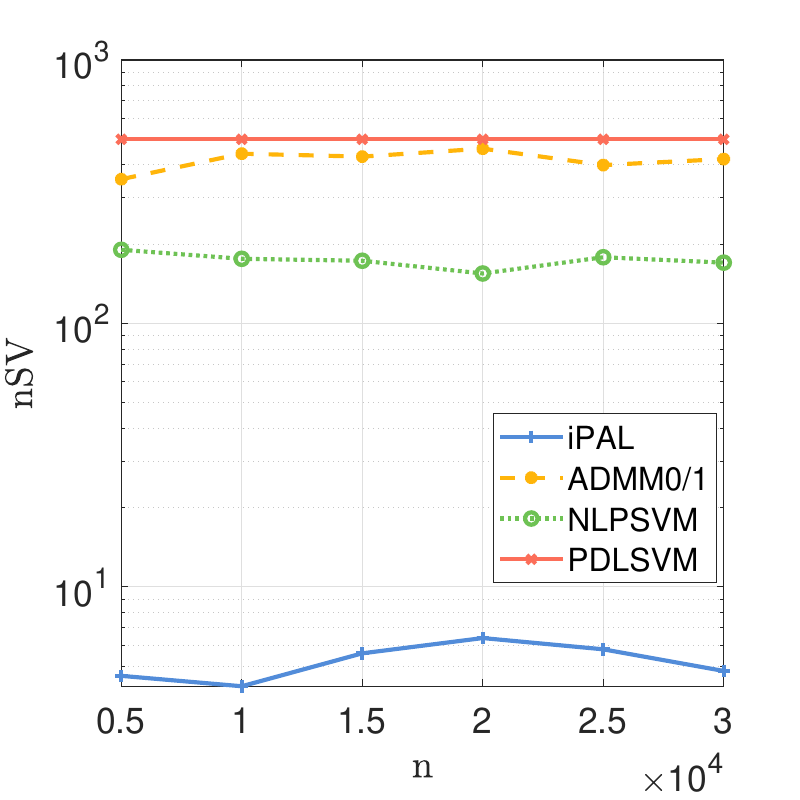}
		\end{minipage}%
	}%
	\subfigure{
		\begin{minipage}[t]{0.5\linewidth}
			\centering
			\includegraphics[width=2.9in]{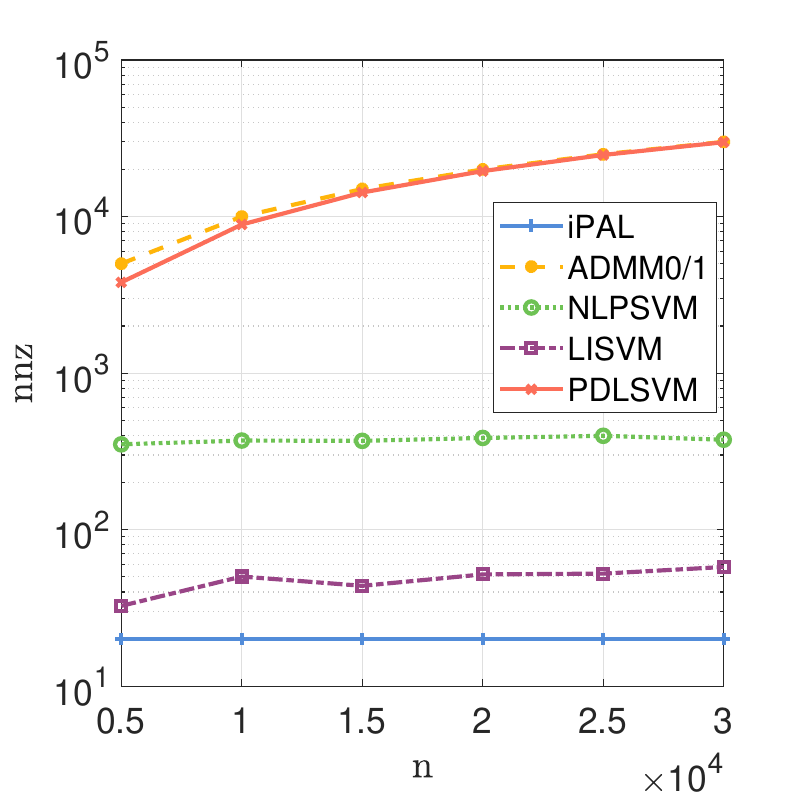}
		\end{minipage}%
	}%
	\centering
	\caption{Comparison results on simulated dataset with $m = 1000$, $r = 0.1$ and $n \in \{ 5000,10000,\cdots,30000 \}$.}
	\label{fig4}
\end{figure}

\textbf{Test I.} We fix $m = 1000$, $r = 0.1$ and vary $n \in \{ 5000,10000,\cdots,30000 \}$. 
In this test, we can see from Fig.~\ref{fig4} that except NLPSVM, all the other solvers achieve the best \texttt{Acc}. Particularly, iPAL spends the least amount of \texttt{Time} with the fewest \texttt{nSV} and \texttt{nnz}. ADMM0/1 is the second fastest solver in this test, but its \texttt{nnz} is much larger and increases as $n$ grows.
This is because this algorithm is designed for a SVM problem without a sparsity constraint on its solutions. 
PDLSVM also shows a significant increase on \texttt{nnz} when $n$ rises, whereas \texttt{nnz} of the other three solvers remain stable. 
As the number of samples $m$ is fixed, the numbers for \texttt{nSV} of all the algorithms are steady.

\begin{figure}[htbp]	
	\subfigure{
		\begin{minipage}[t]{0.5\linewidth}
			\centering
			\includegraphics[width=2.9in]{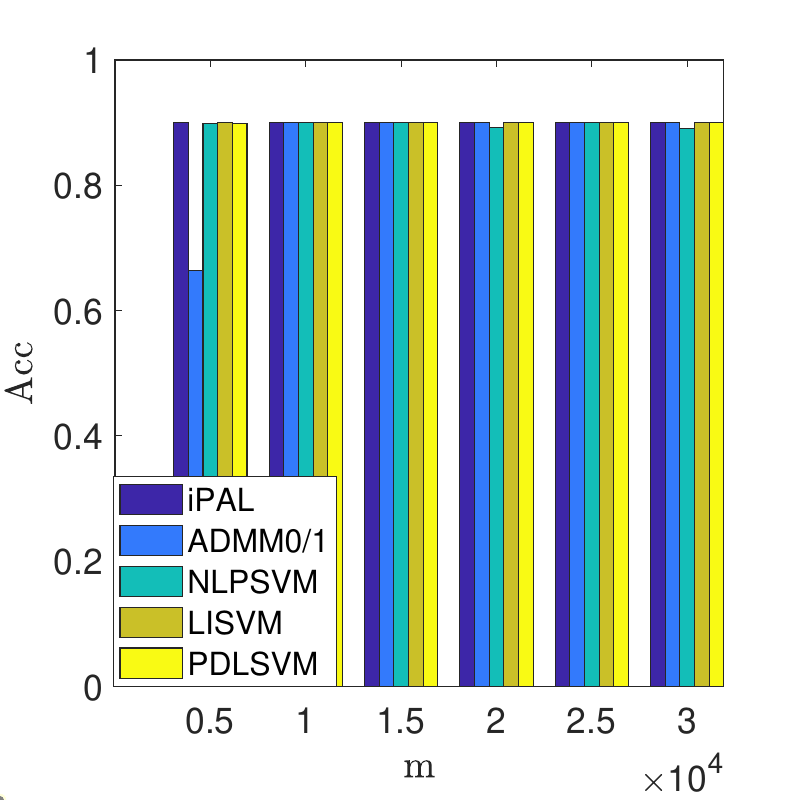}
		\end{minipage}%
	}%
	\subfigure{
		\begin{minipage}[t]{0.5\linewidth}
			\centering
			\includegraphics[width=2.9in]{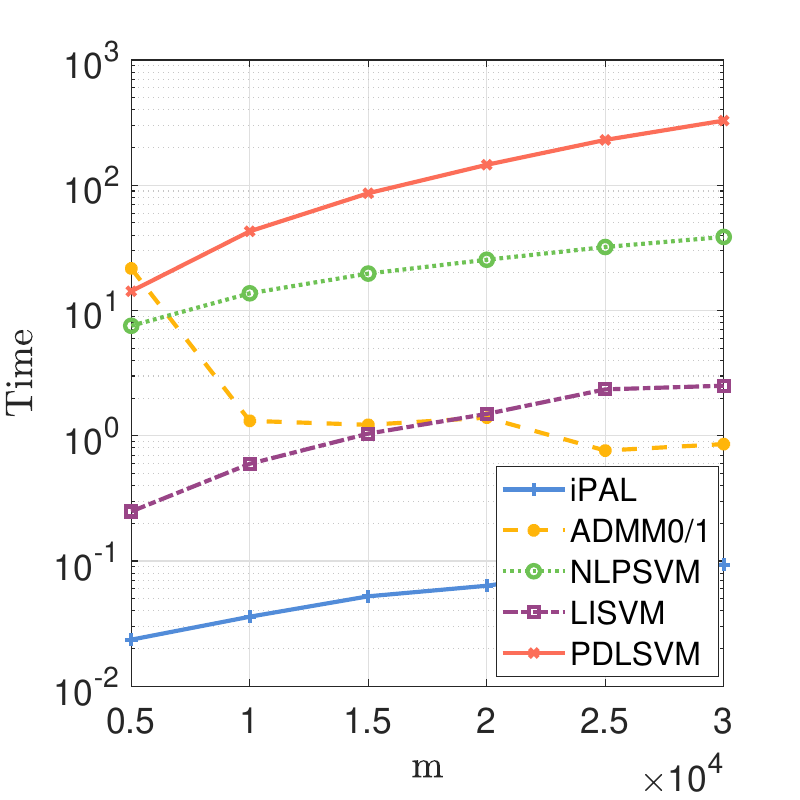}
		\end{minipage}%
	}%
	
	\subfigure{
		\begin{minipage}[t]{0.5\linewidth}
			\centering
			\includegraphics[width=2.9in]{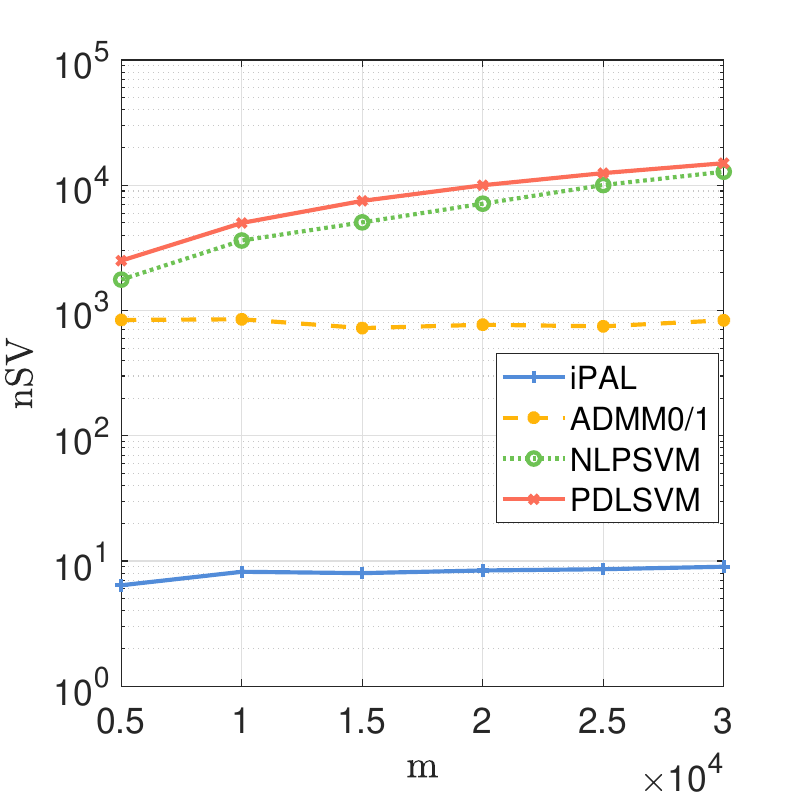}
		\end{minipage}%
	}%
	\subfigure{
		\begin{minipage}[t]{0.5\linewidth}
			\centering
			\includegraphics[width=2.9in]{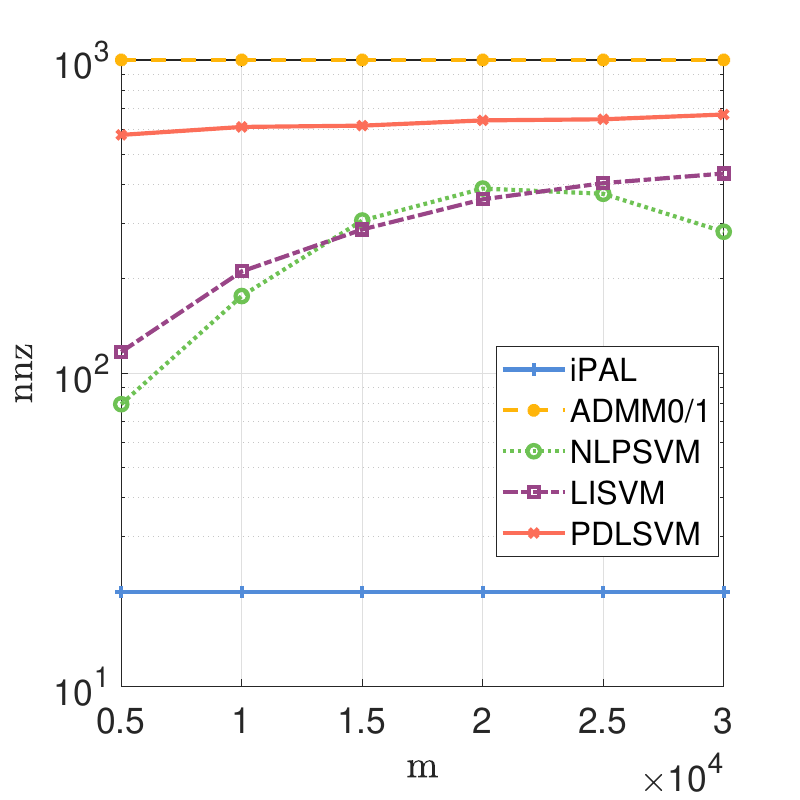}
		\end{minipage}%
	}%
	\centering
	\caption{Comparison results on simulated dataset with $n = 1000$, $r = 0.1$ and $m \in \{ 5000,10000,\cdots,30000 \}$.}
	\label{fig5}
\end{figure}

\textbf{Test II.} We fix $n = 1000$, $r = 0.1$ and alter $m \in \{ 5000,10000,\cdots,30000 \}$. 
Please refer to Fig.~\ref{fig5} for the discussion below.
Again, iPAL performs best on all the evaluating metrics. It has much smaller \texttt{nnz} and \texttt{nSV} than other solvers.
This
significantly reduces the dimension of data matrix and thus contributes to the lower computational cost. In particular, \texttt{Time} of iPAL is almost one order faster than that of ADMM0/1 and LISVM. 
When $m$ becomes larger, there are significant increases on \texttt{nSV} of NLPSVM and PDLSVM, as well as on
\texttt{nnz} of LISVM and NLPSVM.

\begin{figure}[htbp]	
	\subfigure{
		\begin{minipage}[t]{0.5\linewidth}
			\centering
			\includegraphics[width=2.9in]{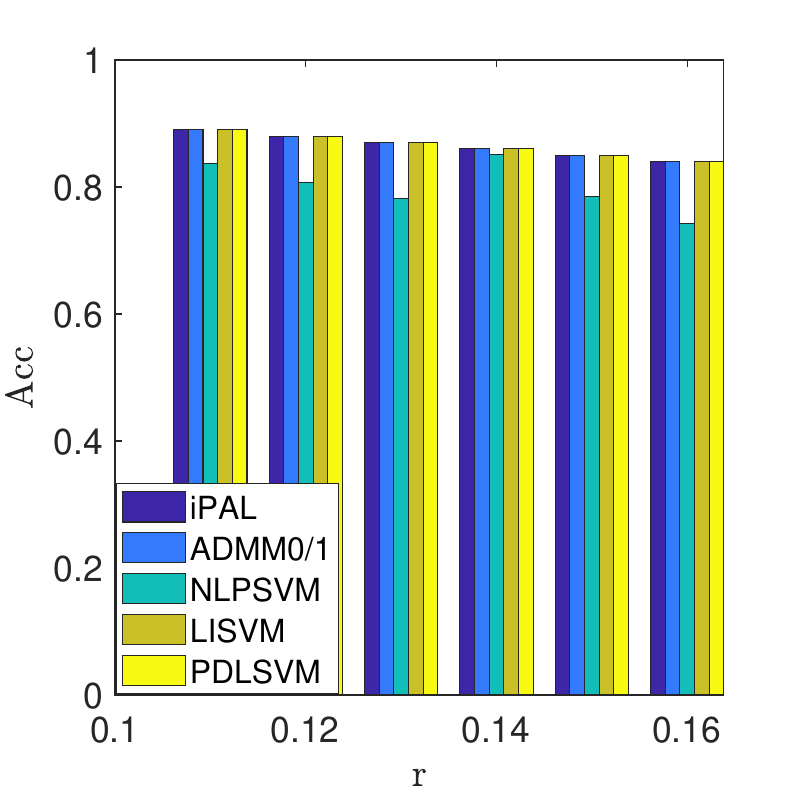}
		\end{minipage}%
	}%
	\subfigure{
		\begin{minipage}[t]{0.5\linewidth}
			\centering
			\includegraphics[width=2.9in]{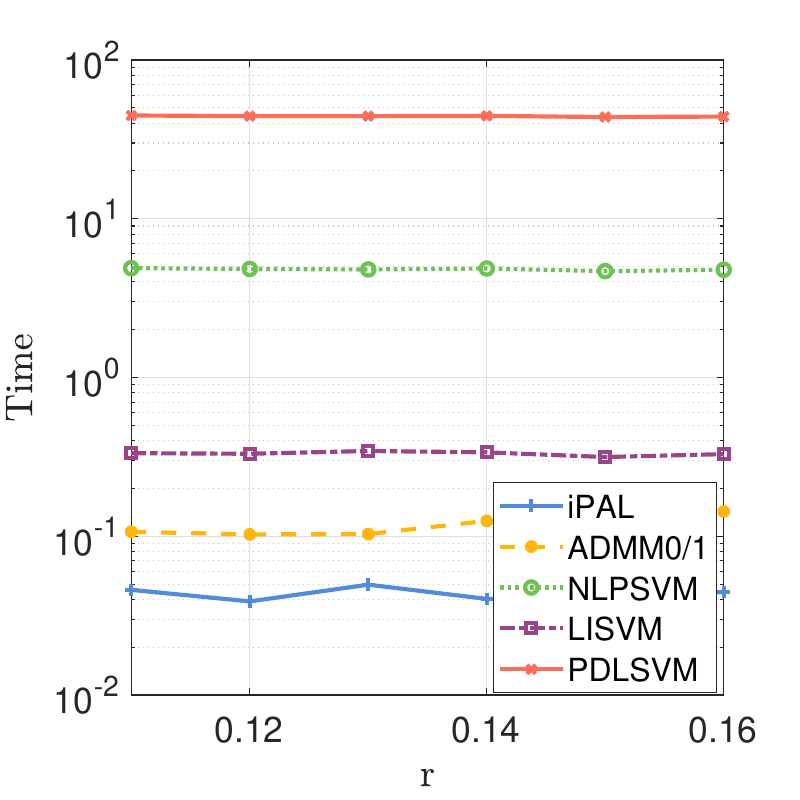}
		\end{minipage}%
	}%
	
	\subfigure{
		\begin{minipage}[t]{0.5\linewidth}
			\centering
			\includegraphics[width=2.9in]{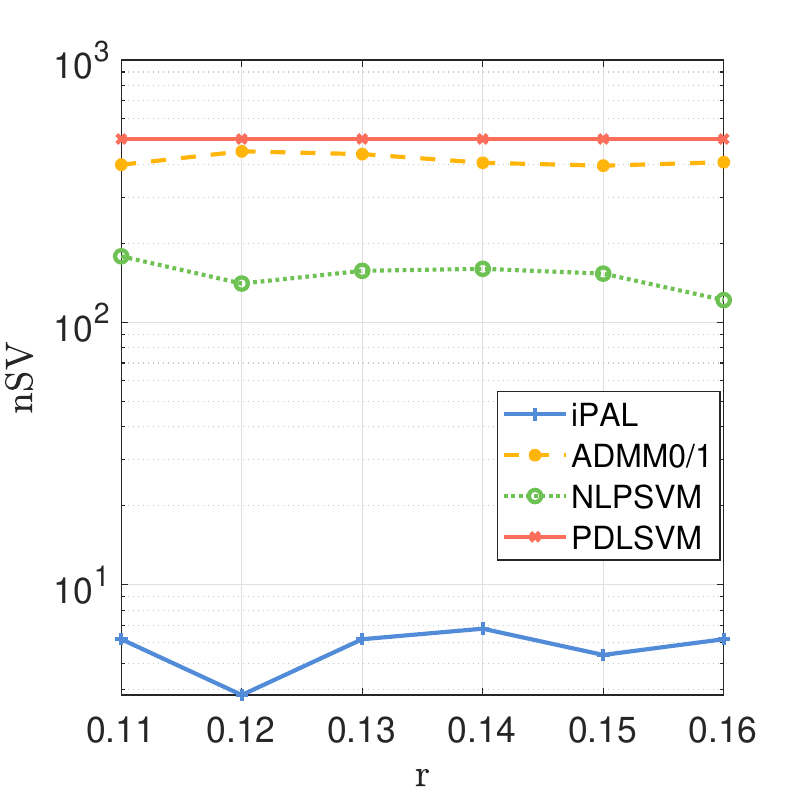}
		\end{minipage}%
	}%
	\subfigure{
		\begin{minipage}[t]{0.5\linewidth}
			\centering
			\includegraphics[width=2.9in]{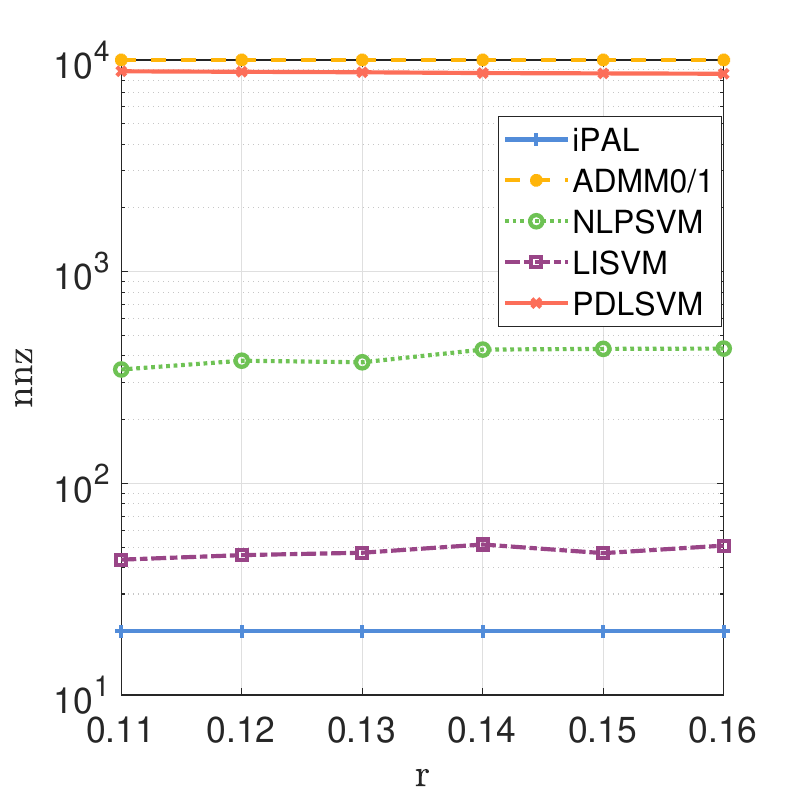}
		\end{minipage}%
	}%
	\centering
	\caption{Comparison results on simulated dataset with $m = 1000$, $n = 10000$ and $r \in \{ 0.11,0.12,\cdots,0.16 \}$.}
	\label{fig6}
\end{figure}

\textbf{Test III.} We fix $m = 1000$, $n = 10000$ and vary $r \in \{ 0.11,0.12,\cdots,0.16 \}$. 
The numerical results are illustrated in Fig.~\ref{fig6}.
It can be observed that with the increase of noise rate, 
the \texttt{Acc} of all the algorithms drops. 
Particularly, the \texttt{Acc} of NLPSVM is more sensitive to noise rate than any other solver. 
The \texttt{nSV}, \texttt{nnz} and \texttt{Time} of all the algorithms are relatively stable with the change of $r$. In this test, iPAL has the best results on all the evaluating metrics. 

The numerical experiments on the simulated data seem to suggest that iPAL is very competitive in terms of 
the four evaluating metrics. 
Similar behaviour of iPAL has also been consistently observed with the real data as we report below. 

\subsection{Experiments on Real Data}

In this section, we will conduct numerical comparison on the real datasets listed in Table \ref{tab2}. 

\begin{example} \label{ex2}
	We select the datasets in Tables \ref{tab2} and \ref{real-data-n<m} with large number of features. Apart from \texttt{gli} and \texttt{dex}, all other datasets are preprocessed by feature-wise scaling to [-1,1]. 
	\begin{table}[htbp]
		\centering
		\caption{Real binary classification datasets with $n > m$ } \label{tab2}
		\begin{tabular}{ccccc}
			\hline
			ID           & Dataset               & Source                     & number of features & number of instances \\ \hline
			\texttt{all} & ALLAML                & feature selection database$\tablefootnote{https://jundongl.github.io/scikit-feature/\label{featureselection}}$ & 7129               & 72                  \\
			\texttt{col} & Colon                 & feature selection database & 2000               & 62                  \\
			\texttt{gli} & GLI85                 & feature selection database & 22283              & 85                  \\
			\texttt{pro} & Prostate              & feature selection database & 5966               & 102                 \\
			\texttt{smk} & SMK187                & feature selection database & 19993              & 187                 \\
			\texttt{dor} & Dorothea              & uci$\tablefootnote{https://archive-beta.ics.uci.edu/datasets}$                        & 100000             & 1950                \\
			\texttt{dex} & Dexter                & uci                        & 20000              & 2600                \\
			\texttt{dbb} & Dbworld\_bodies        & uci                        & 3721               & 64                  \\
			\texttt{dbs} & Dbworld\_subjects      & uci                        & 3721               & 64                  \\
			\texttt{abu} & AP\_Breast\_Uterus      & openML$\tablefootnote{https://www.openml.org/}$                     & 10936              & 468                 \\
			\texttt{alk} & AP\_Lung\_Kidney        & openML                     & 10936              & 368                 \\
			\texttt{aou} & AP\_Ovary\_Uterus       & openML                     & 10936              & 322                 \\
			\texttt{ove} & OVA\_Endometrium       & openML                     & 10936              & 1545                \\
			\texttt{ovk} & OVA\_Kidney            & openML                     & 10936              & 1545                \\
			\texttt{ovo} & OVA\_Ovary             & openML                     & 10936              & 1545                \\
			\texttt{bre} & Breast                & openML                     & 24482            & 97               \\
			\texttt{ova} & Ovarian           & openML                     & 15155              & 253               \\
			\texttt{dbc} & Duke\_breast\_cancer    & openML                    & 7129               & 44                  \\ 
			\texttt{ge1} & Gse2880			& refine.bio $\tablefootnote{https://www.refine.bio/}$ 					&11868					&27    				\\
		\texttt{ge2}	& Gse7670			& refine.bio					&11868				&54				\\
		\texttt{ge3}	& Gse25099			& refine.bio					&16738				&79					\\
		\texttt{ge4}	& Gse27612			& refine.bio					&11868				&195				\\		
			\hline
		\end{tabular}
	\end{table}

\begin{table}[]
	\centering
	\caption{Real binary classification datasets with $n < m$ } \label{real-data-n<m}
	\begin{tabular}{ccccc}
		\hline
		ID           & Dataset                 & Source & number of features & number of instances \\ \hline
		\texttt{chr} & Christine               & openML & 1636               & 5418                \\
		\texttt{jas} & Jasmine                 & openML & 144                & 2984                \\
		\texttt{mad} & Madeline                & openML & 259                & 3140                \\
		\texttt{phi} & Philippine              & openML & 308                & 5382                \\
		\texttt{hiv} & Hiva\_Agnostic           & openML & 1617               & 4229                \\
		\texttt{gui} & Guillermo               & openML & 4296               & 20000               \\
		\texttt{evi} & Evita                   & openML & 3000               & 20000               \\
		\texttt{bio} & Bioresponse             & openML & 1776               & 3751                \\
		\texttt{sw1} & Swarm\_Aligned           & uci    & 2400               & 24016               \\
		\texttt{sw2} & Swarm\_Flocking          & uci    & 2400               & 24016               \\
		\texttt{sw3} & Swarm\_Grouped           & uci    & 2400               & 24016               \\
		\texttt{int} & Internet-Advertisements & uci    & 1559               & 3279                \\
		\texttt{qar} & QSAR\_aquatic\_receptor   & uci    & 1024               & 8992                \\
		\texttt{qot} & QSAR\_oral\_toxity        & uci    & 1024               & 1687                \\ \hline
	\end{tabular}
\end{table}

\end{example}

In the following experiments, for iPAL, we set $\lambda = 1$, $\rho = 1$, $\mu = 10^{-2}$ and $s$ was chosen 
from $\{\lceil0.001p\rceil,\lceil 0.002p  \rceil,\cdots,\lceil 0.01p \rceil ,\lceil 0.02p \rceil, \cdots, \lceil 0.1p \rceil,\lceil 0.2p  \rceil \cdots,p \}$, where $\lceil \cdot \rceil$ is the ceil function. For all other algorithms, the parameters used for trade-off between regularizer and loss function were chosen from $\{ 10^{-5}, 10^{-4}, \cdots, 10^{5} \}$. We conduct five-fold cross-validation on all the datasets in Tables \ref{tab2} and \ref{real-data-n<m},  and the average result are summarized into Table \ref{tab3} and \ref{real_res_n<m}. The PDLSVM fails to give a solution within 2 hours when solving \texttt{dor} and LISVM with $\ell_1$ regularizer does not provide a way to compute \texttt{nSV}, so the corresponding results are indicated by ``--".  
The final results are reported in Tables \ref{tab3} and \ref{real_res_n<m}. 
We have the following observations.

\begin{itemize}
	\item[(i)] {\bf On \texttt{Acc}}. iPAL has the best \texttt{Acc} on most of datasets. When compared with ADMM0/1, iPAL achieves higher \texttt{Acc} with much smaller \texttt{nnz}. This means that the cardinality constraint is beneficial to improving performance of SVM when the number of features is large.
	
	\item[(ii)] {\bf On \texttt{Time}}. iPAL also shows competitive \texttt{Time} in this comparison. For example, on the \texttt{all}, \texttt{gli} and \texttt{ge4}, \texttt{Time} of iPAL is less than 1/3 that of LISVM and even 1\% of the time by ADMM0/1. The high speed of iPAL mainly benefits from the reduction on \texttt{nSV} and \texttt{nnz}, which are obtained by the proximal operator of hard margin loss function and projection of cardinality constraint.
	
	\item[(iii)] {\bf On \texttt{nSV}}. Both iPAL and NLPSVM have small \texttt{nSV}. However, NLPSVM tends to be aggressive on reducing \texttt{nSV} and causes the low \texttt{Acc}, see, for instance, \texttt{all}, \texttt{dbb} and \texttt{bre}. We can also observe that iPAL has smaller \texttt{nSV} than that of ADMM0/1 on almost all the datasets. A possible explanation is that when the redundant features are eliminated, it is easier for a classifier to identify support vectors. 
	
	\item[(iv)] {\bf On \texttt{nnz}}. iPAL, NLPSVM and LISVM show significant reduction on \texttt{nnz} of solution.  In the case of $n \gg m$, feature selection is particularly effective because iPAL tends to have much smaller \texttt{nnz} with higher \texttt{Acc}. NLPSVM or LISVM has smaller \texttt{nnz} than that of iPAL in some cases such as \texttt{chr}, \texttt{gui} and \texttt{qar}, but iPAL has better \texttt{Acc} in those cases. 
	
\end{itemize}

%shows that iPAL can achieve best \texttt{Acc} and \texttt{Time} on most datasets in Ex. \ref{ex2}. Compared with ADMM0/1 without sparsity restriction on solutions, iPAL spends less \texttt{Time} computing solutions with much smaller \texttt{nSV} and \texttt{nnz} without sacrificing \texttt{Acc}. NLPSVM and LISVM are also competitive on reducing \texttt{nSV} and \texttt{nnz}, but the \texttt{Acc} will be lower than that of iPAL, for example \texttt{bre}, \texttt{ova} and \texttt{dbc}.

\begin{table}[htbp]
	\caption{Experiment results on real datasets with $n > m$} \label{tab3}
	\resizebox{\textwidth}{100mm}
	{\begin{tabular}{c|cccccccccc}
			\hline
			& \multicolumn{5}{c|}{\texttt{Acc}~(\%)}                                                                      & \multicolumn{5}{c}{\texttt{Time} (sec)}                                                \\ \cline{2-11} 
			& \texttt{iPAL} & \texttt{ADMM0/1} & \texttt{NLPSVM} & \texttt{LISVM} & \multicolumn{1}{c|}{\texttt{PDLSVM}} & \texttt{iPAL} & \texttt{ADMM0/1} & \texttt{NLPSVM} & \texttt{LISVM} & \texttt{PDLSVM} \\ \hline
			\texttt{all} & 98.57          & 96.07            & 91.79           & 91.61          & \multicolumn{1}{c|}{98.57}           & 9.788e-3       & 2.168e-1         & 1.826e-1        & 3.498e-2       & 1.287e+1        \\
			\texttt{col} & 90.00          & 87.38            & 80.71           & 80.95          & \multicolumn{1}{c|}{85.71}           & 3.817e-3       & 2.407e+0         & 6.580e-2        & 3.344e-3       & 8.048e-1        \\
			\texttt{gli} & 88.24          & 88.24            & 82.35           & 88.24          & \multicolumn{1}{c|}{88.24}           & 2.027e-1       & 9.5501e-1        & 3.502e-1        & 7.980e-1       & 2.147e+2        \\
			\texttt{pro} & 95.00          & 93.00            & 93.00           & 94.00          & \multicolumn{1}{c|}{90.18}           & 6.852e-2       & 2.171e+0         & 1.176e-1        & 5.537e-2       & 9.259e+0        \\
			\texttt{smk} & 77.46          & 74.79            & 72.11           & 73.68          & \multicolumn{1}{c|}{75.93}           & 8.418e-1       & 1.146e+1         & 3.974e-1        & 8.095e+0       & 1.621e+2        \\
			\texttt{dor} & 93.39          & 92.52            & 80.00           & 93.30          & \multicolumn{1}{c|}{--}              & 1.113e+0       & 1.918e+0         & 3.405e+0        & 2.305e-1       & --              \\
			\texttt{dex} & 95.00          & 94.67            & 70.17           & 91.33          & \multicolumn{1}{c|}{94.33}           & 1.759e-2       & 4.538e-1         & 8.497e-1        & 1.375e-1       & 1.808e+2        \\
			\texttt{dbb} & 90.42          & 89.17            & 78.75           & 86.25          & \multicolumn{1}{c|}{90.83}           & 2.326e-2       & 1.163e+0         & 1.867e-1        & 2.325e-2       & 5.512e+0        \\
			\texttt{dbs} & 88.75          & 88.75            & 84.17           & 87.08          & \multicolumn{1}{c|}{77.92}           & 1.458e-3       & 2.811e-3         & 4.442e-3        & 1.951e-3       & 1.714e-2        \\
			\texttt{abu} & 96.38          & 95.94            & 86.77           & 95.93          & \multicolumn{1}{c|}{94.88}           & 5.454e-1       & 1.064e+1         & 3.232e+0        & 6.317e-1       & 4.638e+1        \\
			\texttt{alk} & 97.92          & 97.40            & 90.66           & 97.92          & \multicolumn{1}{c|}{96.88}           & 4.640e-1       & 1.707e+0         & 2.578e+0        & 8.298e-1       & 4.403e+1        \\
			\texttt{aou} & 90.06          & 88.48            & 84.21           & 89.73          & \multicolumn{1}{c|}{85.10}           & 8.096e-1       & 3.056e+0         & 2.006e+0        & 1.629e+0       & 4.216e+1        \\
			\texttt{ove} & 96.63          & 96.50            & 96.05           & 96.38          & \multicolumn{1}{c|}{96.05}           & 2.057e+0       & 1.218e+1         & 1.997e+1        & 1.247e+0       & 7.398e+1        \\
			\texttt{ovk} & 98.71          & 98.71            & 88.22           & 98.58          & \multicolumn{1}{c|}{97.86}           & 2.738e+0       & 1.028e+1         & 1.927e+1        & 2.750e+0       & 7.328e+1        \\
			\texttt{ovo} & 92.49          & 91.78            & 87.18           & 92.36          & \multicolumn{1}{c|}{89.26}           & 3.259e+0       & 1.596e+1         & 1.916e+1        & 1.985e+1       & 7.285e+1        \\
			\texttt{bre} & 79.35          & 70.93            & 75.04           & 76.19          & \multicolumn{1}{c|}{75.34}           & 3.984e-2       & 4.084e-2         & 2.025e-1        & 9.198e-2       & 2.634e+2        \\
			\texttt{ova} & 100.0          & 100.0            & 98.80           & 100.0          & \multicolumn{1}{c|}{99.20}           & 5.618e-2       & 1.148e+0         & 5.429e-1        & 1.978e-1       & 8.554e+1        \\
			\texttt{dbc} & 90.83          & 90.83            & 81.67           & 88.33          & \multicolumn{1}{c|}{86.67}           & 2.517e-2       & 3.153e-1         & 1.502e-1        & 2.063e-2       & 1.340e+1        \\
			\texttt{ge1} & 85.14          & 81.14            & 76.00           & 82.29          & \multicolumn{1}{c|}{81.14}           & 7.778e-3       & 1.978+0          & 1.116e-1        & 8.129e-2       & 4.198e+1        \\
			\texttt{ge2} & 98.00          & 98.00            & 90.57           & 96.00          & \multicolumn{1}{c|}{96.00}           & 9.653e-2       & 1.506e+0         & 3.762e-1        & 8.858e-2       & 4.275e+1        \\
			\texttt{ge3} & 100.0          & 100.0            & 96.56           & 100.0          & \multicolumn{1}{c|}{100.0}           & 1.489e-2       & 8.836e-1         & 7.921e-1        & 9.266e-2       & 1.005e+2        \\
			\texttt{ge4} & 100.0          & 100.0            & 90.26           & 100.0          & \multicolumn{1}{c|}{100.0}           & 2.440e-2       & 1.184e+0         & 1.318e+0        & 9.131e-2       & 4.588e+1        \\
			 \hline
			& \multicolumn{5}{c}{\texttt{nSV}}                                                                            & \multicolumn{5}{c}{\texttt{nnz}}                                                       \\ \cline{2-11} 
			& \texttt{iPAL} & \texttt{ADMM0/1} & \texttt{NLPSVM} & \texttt{LISVM} & \multicolumn{1}{c|}{\texttt{PDLSVM}} & \texttt{iPAL} & \texttt{ADMM0/1} & \texttt{NLPSVM} & \texttt{LISVM} & \texttt{PDLSVM} \\ \hline
			\texttt{all} & 22           & 49             & 11            & --             & \multicolumn{1}{c|}{20}            & 43             & 7130             & 35            & 63           & 4206          \\
			\texttt{col} & 11           & 36             & 16            & --             & \multicolumn{1}{c|}{43}            & 7              & 2001             & 42            & 8            & 1485            \\
			\texttt{gli} & 18           & 37             & 20            & --             & \multicolumn{1}{c|}{13}              & 45             & 22284            & 56              & 1012         & 10846         \\
			\texttt{pro} & 14           & 59             & 41            & --             & \multicolumn{1}{c|}{28}            & 36             & 5967             & 52              & 159          & 3503            \\
			\texttt{smk} & 57           & 127            & 38            & --             & \multicolumn{1}{c|}{145}           & 800            & 19994            & 175           & 14732        & 13804         \\
			\texttt{dor} & 121            & 810            & 275           & --             & \multicolumn{1}{c|}{--}              & 301            & 85488          & 2525          & 343          & --              \\
			\texttt{dex} & 362          & 409              & 189           & --             & \multicolumn{1}{c|}{480}             & 1200           & 9244           & 634           & 1291         & 19999           \\
			\texttt{dbb} & 33             & 49             & 12            & --             & \multicolumn{1}{c|}{51}            & 189            & 3971           & 747             & 135            & 2340          \\
			\texttt{dbs} & 39           & 45             & 19            & --             & \multicolumn{1}{c|}{26}            & 73             & 193            & 44            & 44           & 138           \\
			\texttt{abu} & 65           & 124              & 167           & --             & \multicolumn{1}{c|}{97}            & 219            & 10937            & 209           & 355          & 4874          \\
			\texttt{alk} & 61           & 87             & 130           & --             & \multicolumn{1}{c|}{72}              & 329            & 10937            & 154           & 1010         & 5208          \\
			\texttt{aou} & 122          & 157            & 120           & --             & \multicolumn{1}{c|}{72}            & 547            & 10937            & 201           & 1759         & 5728          \\
			\texttt{ove} & 111          & 162              & 749           & --             & \multicolumn{1}{c|}{353}           & 438            & 10937            & 571           & 129            & 3601          \\
			\texttt{ovk} & 107            & 159            & 543           & --             & \multicolumn{1}{c|}{381}           & 766            & 10937            & 439           & 371          & 4331          \\
			\texttt{ovo} & 227          & 272            & 479           & --             & \multicolumn{1}{c|}{200}           & 657            & 10937            & 500           & 3244         & 4539            \\
			\texttt{bre} & 43           & 63             & 9             & --             & \multicolumn{1}{c|}{13}            & 74             & 24482            & 10            & 25           & 18099         \\
			\texttt{ova} & 27           & 50             & 12              & --             & \multicolumn{1}{c|}{46}            & 46             & 15155            & 12            & 10             & 7908          \\
			\texttt{dbc} & 17           & 30             & 8             & --             & \multicolumn{1}{c|}{8}             & 58             & 7130             & 34            & 32           & 4159          \\
	\texttt{ge1} & 10             & 20               & 13              & --             & \multicolumn{1}{c|}{21}              & 36             & 11869            & 22              & 412            & 10982           \\
	\texttt{ge2} & 13              & 26               & 20              & --             & \multicolumn{1}{c|}{10}              & 60             & 11869            & 251             & 416            & 525             \\
	\texttt{ge3} & 12             & 26               & 17              & --             & \multicolumn{1}{c|}{15}              & 51             & 16739            & 33              & 40             & 7053            \\
	\texttt{ge4} & 25             & 38               & 8               & --             & \multicolumn{1}{c|}{54}              & 24             & 11869            & 157             & 10             & 3501            \\ \hline
	\end{tabular}}{ Note: the PDLSVM fails to give a solution within 2 hours when solving \texttt{dor}. LISVM with $\ell_1$ regularizer does not provide a way to compute \texttt{nSV}}
\end{table}

\begin{table}[]
	\caption{Experiment results on real datasets with $n < m$} \label{real_res_n<m}
	\resizebox{\textwidth}{70mm}
	{\begin{tabular}{c|cccccccccc}
			\hline
			& \multicolumn{5}{c|}{\texttt{Acc}~(\%)}                                                                     & \multicolumn{5}{c}{\texttt{Time} (sec)}                                               \\ \cline{2-11} 
			& \texttt{iPAL} & \texttt{ADMM0/1} & \texttt{NLPSVM} & \texttt{LISVM} & \multicolumn{1}{c|}{\texttt{PDLSVM}} & \texttt{iPAL} & \texttt{ADMM0/1} & \texttt{NLPSVM} & \texttt{LISVM} & \texttt{PDLSVM} \\ \hline
			\texttt{chr} & 73.29         & 54.32            & 68.83           & 72.90          & \multicolumn{1}{c|}{70.27}           & 4.227e+0      & 5.572e+1         & 2.936e+1        & 3.388e+0       & 4.324e+1        \\
			\texttt{jas} & 79.79         & 77.72            & 77.01           & 77.88          & \multicolumn{1}{c|}{77.75}           & 1.367e-1      & 4.826e-1         & 2.242e-1        & 2.231e-1       & 7.480e+0        \\
			\texttt{mad} & 61.97         & 61.88            & 59.20           & 61.88          & \multicolumn{1}{c|}{56.56}           & 2.287e-1      & 2.415e-1         & 1.282e+0        & 9.529e-2       & 9.029e+0        \\
			\texttt{phi} & 71.16         & 70.47            & 72.31           & 72.46          & \multicolumn{1}{c|}{70.35}           & 5.671e-1      & 4.977e-1         & 2.921e+0        & 2.795e+0       & 2.844e+1        \\
			\texttt{hiv} & 96.48         & 93.26            & 96.48           & 96.69          & \multicolumn{1}{c|}{96.50}           & 8.283e-2      & 1.545e+2         & 8.204e+0        & 5.345e-1       & 2.903e+1        \\
			\texttt{gui} & 72.91         & 70.08            & 60.16           & 72.00          & \multicolumn{1}{c|}{70.02}           & 2.483e+0      & 1.892e+2         & 4.979e+2        & 2.388e+0       & 5.304e+2        \\
			\texttt{evi} & 96.59         & 96.59            & 96.70           & 96.59          & \multicolumn{1}{c|}{96.80}           & 3.457e-1      & 1.129e+0         & 9.467e+0        & 1.287e+0       & 4.608e+2        \\
			\texttt{bio} & 73.69         & 61.72            & 52.07           & 76.73          & \multicolumn{1}{c|}{74.94}           & 5.817e-2      & 1.683e+1         & 2.196e+1        & 4.263e-1       & 2.547e+1        \\
			\texttt{sw1} & 100.0         & 99.99            & 72.56           & 100.0          & \multicolumn{1}{c|}{100.0}           & 1.261e+0      & 3.073e+2         & 1.867e+2        & 3.132e+0       & 6.132e+2        \\
			\texttt{sw2} & 99.97         & 99.98            & 72.20           & 99.99          & \multicolumn{1}{c|}{99.94}           & 1.232e+0      & 4.101e+2         & 1.963e+2        & 3.884e+0       & 6.085e+2        \\
			\texttt{sw3} & 100.0         & 99.93            & 72.12           & 100.0          & \multicolumn{1}{c|}{99.94}           & 1.347e+0      & 3.767e+2         & 1.873e+2        & 2.871e+0       & 6.232e+2        \\
			\texttt{int} & 97.10         & 97.07            & 86.00           & 95.88          & \multicolumn{1}{c|}{90.30}           & 6.209e-1      & 1.090e+2         & 2.427e+0        & 2.151e-2       & 1.983e+1        \\
			\texttt{qar} & 89.03         & 77.77            & 88.21           & 86.72          & \multicolumn{1}{c|}{89.75}           & 1.527e-1      & 1.660e+1         & 6.141e-1        & 4.209e-1       & 5.874e+0        \\
			\texttt{qot} & 92.39         & 92.19            & 92.38           & 91.16          & \multicolumn{1}{c|}{92.26}           & 8.777e-1      & 7.737e+0         & 9.038e+0        & 3.165e+0       & 8.099e+1        \\ \cline{2-11} 
			& \multicolumn{5}{c}{\texttt{nSV}}                                                                           & \multicolumn{5}{c}{\texttt{nnz}}                                                      \\ \cline{2-11} 
			& \texttt{iPAL} & \texttt{ADMM0/1} & \texttt{NLPSVM} & \texttt{LISVM} & \multicolumn{1}{c|}{\texttt{PDLSVM}} & \texttt{iPAL} & \texttt{ADMM0/1} & \texttt{NLPSVM} & \texttt{LISVM} & \texttt{PDLSVM} \\ \hline
			\texttt{chr} & 189           & 3272             & 2275            & --             & \multicolumn{1}{c|}{1645}            & 492           & 1611             & 1221            & 364            & 905             \\
			\texttt{jas} & 251           & 128              & 1077            & --             & \multicolumn{1}{c|}{1919}            & 44            & 137              & 33              & 130            & 90              \\
			\texttt{mad} & 17            & 184              & 411             & --             & \multicolumn{1}{c|}{2512}            & 24            & 260              & 119             & 11             & 253             \\
			\texttt{phi} & 4             & 213              & 1781            & --             & \multicolumn{1}{c|}{4316}            & 13            & 309              & 214             & 229            & 194             \\
			\texttt{hiv} & 16            & 3234             & 2932            & --             & \multicolumn{1}{c|}{1210}            & 12            & 1618             & 6               & 189            & 651             \\
			\texttt{gui} & 1638          & 1161             & 6811            & --             & \multicolumn{1}{c|}{15850}           & 258           & 4281             & 1165            & 108            & 2796            \\
			\texttt{evi} & 12            & 8                & 13686           & --             & \multicolumn{1}{c|}{15635}           & 10            & 495              & 106             & 160            & 2027            \\
			\texttt{bio} & 546           & 489              & 1406            & --             & \multicolumn{1}{c|}{605}             & 8             & 1748             & 990             & 336            & 979             \\
			\texttt{sw1} & 235           & 2750             & 17201           & --             & \multicolumn{1}{c|}{7612}            & 481           & 2401             & 337             & 112            & 1298            \\
			\texttt{sw2} & 346           & 6766             & 9830            & --             & \multicolumn{1}{c|}{7182}            & 481           & 2401             & 1278            & 347            & 1133            \\
			\texttt{sw3} & 290           & 5849             & 10174           & --             & \multicolumn{1}{c|}{7156}            & 481           & 2401             & 1269            & 801            & 1155            \\
			\texttt{int} & 154          & 2623             & 2256            & --             & \multicolumn{1}{c|}{201}             & 312           & 1559             & 131             & 430            & 440             \\
			\texttt{qar} & 34            & 1350             & 1190            & --             & \multicolumn{1}{c|}{655}             & 52            & 1025             & 50              & 342            & 572             \\
			\texttt{qot} & 82            & 258              & 5860            & --             & \multicolumn{1}{c|}{2509}            & 93            & 740              & 232             & 777            & 539             \\ \hline
	\end{tabular}}{  Note: LISVM with $\ell_1$ regularizer does not provide a way to compute \texttt{nSV}}
\end{table}

%%%%%%%%%%%%%%%%%%%%%%%%%%%%%%%%%%%%%%%%%%%%5
\section{Conclusion} 

This paper aims to solve a nonsmooth and nonconvex SSVM-HM. We define a {\rm P}-stationary point to characterize its local minimizer. To find a {\rm P}-stationary point, we develop an inexact proximal augmented Lagrangian method (iPAL), which comprises a primal and multiplier step. Based on the {\rm P}-stationarity of the primal step, the inexactness measurement is carefully designed to ensure iPAL converges both globally and at a linear rate. To make the iPAL practically efficient, we design a projected gradient-Newton method (PGN) for computing the primal step with global and local quadratic rate. By the virtue of proximal operator of hard margin loss function and the projection of cardinality constraint, active samples and features can be identified to reduce the dimension of data matrix in PGN. In the extensive numerical comparison, iPAL shows effective reduction on active samples and features while ensuring high classification accuracy and fast computational speed.

This research brings new insights on nonconvex composite optimization with cardinality constraint. An interesting question is how to extend the convergence result to a more general model in which the quadratic term of SSVM-HM is replaced by a smooth function. 
In such an extension, the nice features of the strong convexity as well as the separable property
of the quadratic function would be lost.
Therefore, some proof techniques developed in this paper would not be applicable anymore. 
We leave the extension to future research.

%It seems that the Lyapunov function might not be sufficient descent along with iterate and new merit functions should be constructed in future research.

\section*{Acknowledgements} This work was supported by Fundamental Research Funds for the Central Universities (2022YJS099), the National Natural Science Foundation of China (12131004, 11971052), Beijing Natural Science Foundation (Z190002).

%%%%%%%%%%%%%%%%%%%%%%%%%%%%%%%%
\appendix

\section{Proof of Lemma \ref{Lemma-Prox}}

{\bf Proof}
It follows from (\ref{H-holding}) and (\ref{Prox-J}) that 
$\bfxi \in \Prox_{\beta \lambda} (\bfxi + \beta \bfv)$
if and only if one of the following cases occurs for each $i =1, \ldots, n$:
\[
 \left\{
 \begin{array}{lll}
 	(i) \ & \xi_i = 0, \ & 0 \le v_i < \sqrt{2 \lambda/\beta } \\
 	(ii)\ & \xi_i < 0 , \ & v_i = 0 \\
 	(iii)\ & \xi_i > \sqrt{2 \beta \lambda } , \ & v_i = 0 \\
 	(iv) \ & \xi_i = 0, \ & v_i = \sqrt{2 \lambda/\beta } \\
 	(v) \ & \xi_i = \sqrt{2 \beta \lambda }, \ & v_i = 0 .
 \end{array} 
 \right .
\]
Combing those cases leads to
\[
\left\{
\begin{array}{ll}
	v_i = 0 , & \mbox{if} \ \xi_i \in (-\infty, 0) \cup [\sqrt{2 \beta \lambda}, \; \infty ) \\  [0.2ex]
	v_i \in [0, \; \sqrt{2\lambda / \beta } ],   & \mbox{if} \ \xi_i = 0.
\end{array} 
\right .
\]
This means that $\xi_i$ must satisfy
\[
 \xi_i \in (-\infty, 0] \cup [\sqrt{2 \beta \lambda}, \; \infty )
\]
The characterization (\ref{Pxi}) implies that $\bfxi \in \Prox_{\beta\lambda J(\cdot)} (\bfxi)$. This proves
the necessity part of the lemma. The sufficiency part is by direct verification. 
\hfill $\blacksquare$

%%%%%%%%%%%%%%%%%%%%%%%%%%%%%%%%%%%%%%%%%%%%%%%%%%%%%%%%%%%%%%%%
\section{Proof of Theorem \ref{Thm-Stationarity}} \label{AP_opt}

The strong claim in Thm.~\ref{Thm-Stationarity} basically says that the
concept of the P-stationary point does not introduce any extra points 
other than those of local minimizers. 
Proof of this claim requires certain care and preparation.  
In particular, a smooth reformulation of \eqref{Constrained-HM} plays an important role in this process. We define this reformulation first.

Given a reference point $\bfu^* := (\bfw^*, \bfxi^*)$ belonging to the feasible region of \eqref{Constrained-HM}, let us define
\begin{eqnarray*}
	&&\mathcal{S}^*:= \{ i \in [n]: w^*_i \neq 0 \}, \quad \mathcal{I}^*_-:= \{ i \in [m]: \xi^*_i \leq 0 \}, \quad \mathbb{T}^*:= \{  T \subseteq [n]: T \supseteq \mathcal{S}^*, |T| = s \}.
\end{eqnarray*}
Taking $T^* \in \mathbb{T}^*$, we consider the following nonlinear programming associated with $T^*$ (abbreviated as NLP-$T^*$)
\begin{equation} \tag{NLP-$T^*$}
	\begin{aligned}
		\min_{\bfw,\bfxi} \quad  \frac{1}{2}\| \bfw \|^2, \quad
		\mbox{s.t.} \quad   \bfw_{\oT^*} = 0, \quad 
		\bfxi_{\mathcal{I}^*_-} \leq 0, \quad
		A\bfw + \bfone = \bfxi .
	\end{aligned}
\end{equation}
%Although this problem is also associated with $\mathcal{I}^*_-$, this set can be uniquely identified when $\bfu^*$ is given. 
The Lagrange function of NLP-$T^*$ is denoted by
\begin{align*}
	\mathcal{L}_{T^*} (\bfu,\bfq_w,\bfq_\xi,\bfz):= \frac{1}{2}\| \bfw \|^2 + \langle \bfq_w, \bfw_{\oT^*} \rangle + \langle \bfq_\xi, \bfxi_{\mathcal{I}^*_-} \rangle + \langle \bfz, A\bfw + \bfone - \bfxi \rangle,
\end{align*}
where $(\bfq_w,\bfq_\xi,\bfz) \in \mathbb{R}^{|\overline{T}^*|} \times \mathbb{R}^{|\mathcal{I}_-|} \times \mathbb{R}^m$ are multipliers associated with the three constraints in (NLP-$T^*$).
Thereby, the KKT system of (NLP-$T^*$) can be represented as 
\begin{equation}\label{KKT}
	\left\{  \begin{aligned}
		&(\bfw + A^\top \bfz)_{T^*} = 0, \ \bfw_{\oT^*} = 0, \\
		&\bfz_{\mathcal{I}^*_-} \geq 0, \ \bfxi_{\mI} \leq 0, \ \langle \bfz_{\mI}, \bfxi_{\mI} \rangle = 0, \ \bfz_{\moI} = 0, \\
		& A\bfw + \bfone - \bfxi = 0, \ \bfq_w = -( \bfw + A^\top \bfz )_{\oT^*}, \ \bfq_\xi = \bfz_{\mI}.
	\end{aligned} \right.
\end{equation}
We say $(\bfw, \bfxi)$ satisfying (\ref{KKT}) is a KKT point of (NLP-$T^*$)
with Lagrange multipliers $(\bfq_w,\bfq_\xi,\bfz)$.
We will prove Theorem \ref{Thm-Stationarity} based on the following
arguments.

\begin{itemize}
	\item[$\bullet$] A {\rm P}-stationary point $\bfu^*=(\bfw^*,\bfxi^*)$ is equivalent to a KKT point of (NLP-$T^*$) for any $T^* \in \mathbb{T}^*$ (Lemma \ref{P_stat_KKT}).
	
	\item[$\bullet$] We will prove a local minimizer $\bfu^*$ of \eqref{Constrained-HM} is also a local minimizer of (NLP-$T^*$). Since NLP-$T^*$ is a smooth programming with linear constraint, then $\bfu^*$ is also a KKT point of NLP-$T^*$ (see \cite{nocedal2006numerical}). By using Lemma \ref{P_stat_KKT},
	$\bfu^*$ is a {\rm P}-stationary point of \eqref{HM-J}. This is exactly the necessary optimality condition in Theorem \ref{Thm-Stationarity}.
	
	\item[$\bullet$] Second-order sufficient condition of (NLP-$T^*$) holds at any KKT point $\bfu^*$ (Lemma \ref{sosc-lemma}). Then following from \cite[Theorem 2.4]{nocedal2006numerical}, $\bfu^*$ is a strict local minimizer of (NLP-$T^*$) satisfying quadratic growth condition. This can further lead to the sufficient optimality condition in Theorem \ref{Thm-Stationarity} (ii). 
\end{itemize}

%The relationship between a $P$-stationary point of \eqref{Constrained-HM} and a KKT point of NLP-$T^*$ is given as follows.

\begin{lemma} \label{P_stat_KKT}
	Given a point $\bfu^*:= (\bfw^*,\bfxi^*)$, the following statements hold.
	
	\begin{itemize}
		\item[(i)] If $\bfu^*$ is a $P$-stationary point of \eqref{Constrained-HM}, then it is a KKT point of (NLP-$T^*$) for any $T^* \in \mathbb{T}^*$.
		
		\item[(ii)] If $\bfu^*$ is a KKT point with multiplier $(\bfq_w^*,\bfq_\xi^*,\bfz^*)$ of (NLP-$T^*$) for any $T^* \in \mathbb{T}^*$, then $(\bfq_w^*,\bfq_\xi^*)$ must be taken as
		\begin{align} \label{z-xy}
			\bfq_w^* = \left\{ \begin{aligned}
				&[ \bfw^* + A^\top \bfz^* ]_{\mathcal{S}^*}, &&\ \mbox{if} \ \| \bfw^* \|_0 = s, \\
				& 0, &&\ \mbox{if} \ \| \bfw^* \|_0 < s
			\end{aligned} \right. \ \mbox{and} \ \bfz^*_\xi = \bfz^*_{\mI}.
		\end{align}
		Moreover, $(\bfu^*,\bfz^*)$ is a {\rm P}-stationary pair with constant $\alpha \in (0,\alpha^*)$ and $\beta \in (0, \min\{ \beta_\xi^*,\beta_z^* \})$, where
		\begin{align} \notag
			&\alpha^* := \left\{ \begin{aligned}
				&\frac{|\bfw^*|_{(s)}}{\max_{i \in \overline{\mathcal{S}}^*} | (\bfw^* + A^\top \bfz^* )_i |}, \ &&\mbox{if} \ (\bfw^* + A^\top \bfz^*)_{\overline{\mathcal{S}}^*} \neq 0, \\
				& \infty, \ &&\mbox{otherwise}.
			\end{aligned} \right.\\
			&\beta_\xi^* := \left\{ \begin{aligned}
				& \infty, \ &&\mbox{if} \ \bfxi^* \leq 0, \\
				& \min \{ \frac{{\xi^*_i}^2}{2\lambda}: \xi^*_i > 0, i \in [m] \}, \ &&\mbox{otherwise}.
			\end{aligned} \right. \notag \\
			& \beta_z^* := \left\{ \begin{aligned}
				& \infty, \ &&\mbox{if} \ \bfz^* \leq 0, \\
				& \min \{ \frac{2\lambda}{{z^*_i}^2}: z^*_i > 0, i \in [m] \}, \ &&\mbox{otherwise}. \notag
			\end{aligned} \right.
		\end{align}
	\end{itemize}
\end{lemma}
{\bf Proof}
	(i) If $\bfu^*$ is a ${\rm P}$-stationary point of \eqref{Constrained-HM}, then there exists a {\rm P}-stationary multiplier $\bfz^*$ such that $(\bfu^*, \bfz^*)$ satisfies \eqref{P-stat}. Let us first prove $(\bfw^* + A^\top \bfz^*)_{T^*} = 0$ and $\bfw^*_{\oT^*} = 0$. 
	
	The claim $\bfw^*_{\oT^*} = 0$ directly follows from $T^* \supseteq \mathcal{S}^*$.
	% where $\mathcal{S}^* := \{ i \in [n]: w^*_i \neq 0 \}$. 
If $\| \bfw^* \|_0 = s$, then $T^* \in \mathbb{T}^* = \mathcal{S}^*$. By the definition of $\mathcal{S}^*$, \eqref{proj_explicit} implies $(\bfw^* + A^\top \bfz^*)_{T^*} = 0$. If $\| \bfw^* \| < s$, then $|\bfw^*|_{(s)} = 0$ and $\bfw^* + A^\top \bfz^* = 0$ holds from \eqref{proj_explicit}. The second line of \eqref{KKT} can be obtained from Lemma \ref{Lemma-Prox}. The third line of \eqref{KKT} holds when taking $\bfq_w^* = -( \bfw^* + A^\top \bfz^* )_{\oT^*}$ and $\bfq_\xi^* = \bfz_{\mI}^*$. 
	
	(ii) Let us first verify \eqref{z-xy}. $\bfz^*_\xi = \bfz^*_{\mI}$ directly follows from \eqref{KKT}. If $\| \bfw^* \|_0 = s$, then we can uniquely identify $T^* = \mathcal{S}^*$ for the reasons $| T^* | = | \mathcal{S}^* | = s$ and $T^* \supseteq \mathcal{S}^*$. If $\| \bfw^* \|_0 < s$, \eqref{proj_explicit} leads to $\bfq_w^* = - (\bfw^* + A^\top \bfz^*)_{\oT^*}$ = 0. Taking these two facts as well as \eqref{KKT} into consideration, $\bfq_w^*$ must be taken as \eqref{z-xy}.  
	
	We now prove $\alpha^* > 0$ so that $\alpha \in (0,\alpha^*)$ is well defined. Assume by contradiction that $\alpha^* = 0$, then $|\bfw^*|_{(s)} = 0$ implies $\| \bfw^* \|_0 < s$. Meanwhile, $(\nabla f(\bfu^*) + A^\top \bfz^*)_{\overline{S}^*} \neq 0$ holds. Since $\| \bfw^* \|_0 < s$, for any $i_0 \in \overline{\mathcal{S}}^*$, there exists $T^*_0 \in \mathbb{T}^*$ such that $i_0 \in T^*_0$ and thus $(\nabla f(\bfu^*) + A^\top \bfz^*)_{i_0} = 0$ from the first line of \eqref{KKT}. This contradicts to the fact $(\nabla f(\bfu^*) + A^\top \bfz^*)_{\overline{S}^*} \neq 0$. Therefore, we must have
	 $\alpha^* > 0$.
	
	Now we begin proving $\bfw^* \in \proj (\bfw^* - \alpha ( \bfw^* + A^\top \bfz^* ))$. Since $(\bfw^* + A^\top \bfz^*)_{T^*} = 0$ for any $T^* \in \mathbb{T}^*$ holds from \eqref{KKT} and $T^* \supseteq \mathcal{S}^*$, we have $(\bfw^* + A^\top \bfz^*)_{\mathcal{S}^*} = 0$.  $\mbox{If} \ (\bfw^* + A^\top \bfz^*)_{\overline{\mathcal{S}}^*} \neq 0$, then from the definition of $\alpha^*$, we obtain
	\begin{align} \notag
		| (\bfw^* + A^\top \bfz^* )_i | \leq \max_{i \in \overline{\mathcal{S}}^*} | (\bfw^* + A^\top \bfz^* )_i | \leq |\bfw^*|_{(s)}/\alpha, \ \forall i \in \overline{\mathcal{S}}^*.
	\end{align}
	$\mbox{If} \ (\bfw^* + A^\top \bfz^*)_{\overline{\mathcal{S}}^*} = 0$, then $| (\bfw^* + A^\top \bfz^* )_i | \leq |\bfw^*|_{(s)}/\alpha$ also holds for any $i \in \overline{\mathcal{S}}^*$. \eqref{proj_explicit} indicates that $\bfw^* \in \proj (\bfw^* - \alpha ( \bfw^* + A^\top \bfz^* ))$ is true.
	
	Next we will verify $\bfxi^* \in \prox(\bfxi^* + \beta \bfz^*)$. The definition of $\beta_\xi^*$ and $\beta_z^*$ indicates that for any $\xi^*_i > 0$ and $z^*_j > 0$, we have
	\begin{align*}
		&\beta < \beta_\xi^* \leq (\xi_i^*)^2/(2\lambda) \quad \ \mbox{and} \ \quad \beta < \beta_z^* \leq 2\lambda /(z^*_j)^2, 
	\end{align*}
	which means that $\xi^*_i > \sqrt{2\lambda\beta}$ if $y^*_i >0$, and $z^*_j < \sqrt{2\lambda/\beta}$ if $z^*_j > 0$. This result together with second line in \eqref{KKT} further leads to
	\begin{align*}
		\left\{ \begin{aligned}
			& z^*_i = 0, \ &&\mbox{if} \ \xi^*_i \in (-\infty, 0) \cup (\sqrt{2\lambda\beta}, \infty), \\
			& z^*_i \in [0,\sqrt{2\lambda/\beta}), \ &&\mbox{if} \ \xi^*_i = 0.
		\end{aligned} \right.
	\end{align*}
	Then Lemma \ref{Lemma-Prox} implies that $\bfxi^* \in \prox(\bfxi^* + \beta \bfz^*)$. Finally, $A\bfw^* + \bfone - \bfxi^* = 0$ directly follows from the third line of \eqref{KKT}. 
	This completes the proof.
\hfill $\blacksquare$

Another favorite property of (NLP-$T^*$) is that it naturally satisfies
the second-order necessary condition (SOSC, see e.g. \cite[Theorem 12.5]{nocedal2006numerical}) well defined for smooth optimization.

\begin{lemma} \label{sosc-lemma}
	Given a KKT pair $(\bfu^*, \bfq_w^*, \bfq_\xi^*, \bfz^*)$ of (NLP-$T^*$), the following SOSC naturally holds 
	\begin{align} \label{sosc}
		[ \bfd^w;\bfd^\xi ]^\top \nabla^2_{\bfu,\bfu} \mathcal{L}_{T^*} (\bfu^*, \bfq_w^*, \bfq_\xi^*, \bfz^*) [ \bfd^w;\bfd^\xi ] > 0, 
		\quad \forall\ \ [ \bfd^w;\bfd^\xi ] \in \mathcal{C}^* \backslash \{ 0 \}, 
	\end{align}
	where $\mathcal{C}^*:=\{ (\bfd^w,\bfd^\xi) \in \mathbb{R}^{m+n}: A \bfd^w = \bfd^\xi, \ \bfd^w_{\oT^*} = 0, \bfd^\xi_{\mathcal{I}^*_0} \leq 0, \ \bfd^\xi_{\mathcal{I}^*_+} = 0 \}$ is the critical cone of (NLP-$T^*$), $ \mathcal{I}^*_0:=  \{ i \in [n]: \xi^*_i = 0, z^*_i = 0 \}$ and $ \mathcal{I}^*_+:=  \{ i \in [n]: \xi^*_i = 0, z^*_i > 0 \}$.
\end{lemma}

{\bf Proof}
	The Hessian of the Lagrangian of (NLP-$T^*$) with respect to $\bfu$ can be written as 
	\begin{align*}
		\nabla^2_{\bfu,\bfu} \mathcal{L}_{T^*} (\bfu^*, \bfz_w^*, \bfz_\xi^*, \bfz^*) = \left[\begin{array}{cc}
			I &0 \\
			0 & 0
		\end{array}\right]
	\end{align*}
	Thus, \eqref{sosc} actually means 
	\begin{align*}
		\| \bfd^w \|^2 > 0 \quad \forall\ \ [ \bfd^w;\bfd^\xi ] \in \mathcal{C}^* \backslash \{ 0 \} . 
	\end{align*}
	Notice that $\bfd^w \neq 0$ must hold.
	Otherwise $A\bfd^w = \bfd^\xi$ would imply $0 = (\bfd^w,\bfd^\xi)$,
	contradicting with the assumption that it is not zero. 
	Therefore $\| \bfd^w \|^2 > 0$ and the SOSC naturally holds.
\hfill $\blacksquare$

%Now we are ready to prove Theorem~\ref{Thm-Stationarity}.

\noindent\textbf{Proof of Theorem \ref{Thm-Stationarity}.} 
We define the feasible regions of \eqref{Constrained-HM} and (NLP-$T^*$) 
\begin{align*}
	&\mathcal{F}:= \{ \bfu=(\bfw, \bfxi) : \| \bfw \|_0 \leq s, A \bfw + \bfone - \bfxi = 0 \},\\
&	\mathcal{F}_{T^*} := \{ \bfu=(\bfw, \bfxi) : \bfw_{\overline{T}^*} = 0, \bfxi_{\mathcal{I}_-^*} \leq 0, A \bfw + \bfone - \bfxi = 0 \}. 
\end{align*}

\textbf{(a) local minimizer $\Longrightarrow$ {\rm P}-stationary point.}
Given a local minimizer $\bfw^*$ of \eqref{HM-J}, $\bfu^*:= ( \bfw^*,\bfxi^* )$ with $\bfxi^* = A\bfw^* + \bfone$ is a local minimizer of \eqref{Constrained-HM}, there exists $\epsilon^* >0$ such that 
	\begin{align} \label{loc_min}
		\frac{1}{2}\| \bfw \|^2 + \lambda J( \bfxi) \geq \frac{1}{2}\|\bfw^*\|^2 + \lambda J(\bfxi^*),\ \mbox{for all} \ \bfu \in \mathcal{N}(\bfu^*, \epsilon^*) \cap \mathcal{F}.
	\end{align}
	Given $T^* \in \mathbb{T}^*$, we have $\mathcal{F}_{T^*} \subseteq \mathcal{F}$ and let us consider $\bfu \in \mathcal{N}(\bfu^*, \epsilon^*) \cap \mathcal{F}_{T^*}$. Since $\mathcal{F}_{T^*} \subseteq \mathcal{F}$ and $J(\bfxi^*) \geq J(\bfxi)$, from \eqref{loc_min}, we have
	\begin{align*} 
		\frac{1}{2}\| \bfw \|^2 \geq \frac{1}{2}\|\bfw^*\|^2,\ \mbox{for all} \ \bfu \in \mathcal{N}(\bfu^*, \epsilon^*) \cap \mathcal{F}_{T^*},
	\end{align*}
	which means that $\bfu^*$ is also a local minimizer of (NLP-$T^*$).  Noticing that for each $T^* \in \mathbb{T}$, (NLP-$T^*$) is a smooth nonlinear optimization problem with linear constraints, we can further deduce that for any $T^* \in \mathbb{T}^*$, $\bfu^*$ is a KKT point of (NLP-$T^*$) with corresponding multiplier $(\bfq_w^*,\bfq_\xi^*,\bfz^*)$ in \eqref{z-xy}. 
	Thus, we can prove the desired conclusion by Lemma \ref{P_stat_KKT} (ii).
	
	\textbf{(b) {\rm P}-stationary point $\Longrightarrow$ local minimizer.} Given a {\rm P}-stationary point $\bfw^*$ of \eqref{HM-J}, $\bfu^* = (\bfw^*, \bfxi^*)$ with $\bfxi^*:= A \bfw^* + \bfone$ is a KKT point of (NLP-$T^*$) for each $T^* \in \mathbb{T}^*$ from Lemma \ref{P_stat_KKT}. Meanwhile, noticing that the SOSC \eqref{sosc} holds, it follows from \cite[Theorem 12.6]{nocedal2006numerical} that there exists $\epsilon_{T^*} > 0$ and $ c_{T^*} > 0$ such that 
	\begin{align} \label{nlp-qua-gro}
		\frac{1}{2}\| \bfw \|^2 \geq \frac{1}{2}\|\bfw^*\|^2 +  c_{T^*} \| \bfu - \bfu^* \|^2, \ \forall \bfu \in \mathcal{N} (\bfu^*, \epsilon_{T^*}) \cap \mathcal{F}_{T^*} .
	\end{align}   
	Denote $ c^*: = \min_{T^* \in \mathbb{T}^*}  c_{T^*}$. Now we take a radius $\epsilon^*$ satisfying
	\begin{align} \label{eps_range}
		\epsilon^* < \min_{T^* \in \mathbb{T}^*} \epsilon_{T^*}
		\quad \mbox{and} \quad
		  c^* {\epsilon^*}^2 < \lambda/2.  
	\end{align}
	We also assume that $\epsilon^*$ is small enough such that for any $\bfu \in \mathcal{N}(\bfu^*, \epsilon^*)$, the following relationships hold
	\begin{gather} 
		 \mathcal{S}^* \subseteq \{ i \in [n]: w_i \neq 0 \} \ \mbox{and} \ \{ i \in [n]: \xi^*_i > 0\} \subseteq \{ i \in [n]: \xi_i > 0 \}, \label{subseteq_xy} \\
		| \| \bfw \|^2 - \| \bfw^* \|^2 | < \lambda, \label{cont-f}
	\end{gather}
	where the inequality follows from the continuity of $\| \cdot \|^2$. Particularly, \eqref{subseteq_xy} further leads to
	\begin{align} \label{indicate_+0}
		 J( \bfxi) \geq J(\bfxi^*).
	\end{align} 
	Denoting $\mathcal{F}^*:= \bigcup_{T^* \in \mathbb{T}^*} \mathcal{F}_{T^*} \subseteq \mathcal{F}$, then from \eqref{nlp-qua-gro} and \eqref{indicate_+0}, we can obtain
	\begin{align*} 
		\frac{1}{2}\| \bfw \|^2 + \lambda J( \bfxi) \geq \frac{1}{2}\|\bfw^*\|^2 + \lambda J(\bfxi^*) +  c^* \| \bfu - \bfu^* \|^2,\ \forall \bfu \in \mathcal{N}(\bfu^*, \epsilon^*) \cap \mathcal{F}^*.
	\end{align*}
	If we take $\bfu \in \mathcal{N}(\bfu^*, \epsilon^*) \cap 
	(\mathcal{F} \backslash \mathcal{F}^*)$, considering $\mathcal{S}^* \subseteq \{ i \in [n]: w_i \neq 0 \}$ in \eqref{cont-f} and $\| \bfw \|_0 \leq s$, there must exists $T^* \in \mathbb{T}$ such that $\bfw_{\oT^*} = 0$. This together with $\bfu \notin \F^*$ lead to $\bfxi_{\mI} \nleq 0$. There exists an index $i_0 \in \mI$ such that $\xi_{i_0} > 0$. Combining this with \eqref{subseteq_xy} leads to $J( \bfxi) \geq J(\bfxi^*) + 1$. Then taking \eqref{cont-f} and \eqref{indicate_+0} into consideration, we have
	\begin{align*}
		\frac{1}{2}\| \bfw \|^2 + \lambda J( \bfxi) &\geq \frac{1}{2}\|\bfw^*\|^2 + \lambda J(\bfxi^*) + \lambda/2  \mathop{\geq}\limits^{(\ref{eps_range})} \frac{1}{2}\|\bfw^*\|^2 + \lambda J(\bfxi^*) +   c^* \| \bfu - \bfu^* \|^2. 
	\end{align*} 
Overall, we have obtained
	\begin{align*}
		\frac{1}{2}\| \bfw \|^2 + \lambda J( \bfxi) \geq \frac{1}{2}\|\bfw^*\|^2 + \lambda J(\bfxi^*) +   c^* \| \bfu - \bfu^* \|^2,~ \forall \bfu \in \mathcal{N}(\bfu^*, \epsilon^*) \cap \mathcal{F}.
	\end{align*}
Finally, \eqref{Quadratic-Growth} follows from the definition of $\mathcal{F}$. \hfill $\blacksquare$
%	 from the definition of $\mathcal{F}^*$, one of the following two cases must hold.
%	
%	\textbf{Case I:} For any $T^* \in \mathbb{T}^*$, we have
%	$\bfw_{\oT^*} \neq 0$.  
%	This will lead to $\| \bfw \|_0 \geq s+1$. Indeed, if this is not true ($\| \bfw \|_0 \leq s$), we can take $T^* = \{ i \in [n]: w_i \neq 0 \}$. Then \eqref{subseteq_xy} and $\| \bfw \|_0 \leq s$ imply $T^* \in \mathbb{T}^*$ but $\bfw_{\oT^*} = 0$. This is a contraction and thus $\| \bfw \|_0 \geq s+1$ holds, which means that $\delta_{\mathbb{S}}(\bfw) = \infty$.
%	
%	\textbf{Case II:} There exists an index $i_0 \in \mI$ such that $\xi_{i_0} > 0$. Combining this with \eqref{subseteq_xy} leads to $J( \bfxi) \geq J(\bfxi^*) + 1$. Then taking \eqref{nlp-qua-gro}, \eqref{cont-f} and \eqref{indicate_+0} into consideration, we have
%	\begin{align*}
%		\frac{1}{2}\| \bfw \|^2 + \lambda J( \bfxi) + \delta_{\mathbb{S}}(\bfw) &\geq \frac{1}{2}\|\bfw^*\|^2 + \lambda J(\bfxi^*) + \delta_{\mathbb{S}}(\bfw^*) + \lambda/2  \\
%		&\mathop{\geq}\limits^{(\ref{eps_range})} \frac{1}{2}\|\bfw^*\|^2 + \lambda J(\bfxi^*) + \delta_{\mathbb{S}}(\bfw^*) +   c^* \| \bfu - \bfu^* \|^2. 
%	\end{align*}
%	Overall, we arrive at the desired conclusion. 

%%%%%%%%%%%%%%%%%%%%%%%%%%%%%%%%%%%%%%%%%%%%%%%%%%%%%
\section{Proofs on Global Convergence of iPAL} 

In this part, our ultimate goal is to prove Theorem \ref{Thm-Global}.
It is beneficial to briefly explain the main ideas behind our proofs.

\begin{itemize}
	\item[$\bullet$]  First, we will prove Proposition \ref{Prop-Sufficient-Decrease}, including the sufficient decrease of
	 Lyapunov function in \eqref{lya-des}, 
	 boundedness of the sequence $\{ (\bfu^k,\bfz^k) \}_{k \in \mathbb{N}}$,
	  and the convergence of difference of successive iterates \eqref{suc_change}.
	  
	\item[$\bullet$] The boundedness of sequence ensures that there must exist an accumulated point. The inexact criteria \eqref{error-metric} actually means that each iterate approximately satisfies a {\rm P}-stationary system and the degree of approximation can be measured
	 by $\| \bfw^{k+1} - \bfw^k \|$. 
	For such a sequence, each accumulated point is a {\rm P}-stationary point of \eqref{Constrained-HM} by using $\eqref{suc_change}$ and the proximal behavior \cite[Theorem 1.25]{rockafellar1976augmented}. This result is referred to as a subsequence convergence property (see Lemma \ref{sub-con}).
	
	\item[$\bullet$] We will mainly use \cite[Proposition 7]{kanzow1999qp} to prove that the whole sequence generated by iPAL is convergent. The requirements for using this proposition are \eqref{suc_change} and the isolatedness of accumulation points. The isolatedness property follows from Theorem \ref{Thm-Stationarity}.  
\end{itemize}

\noindent\textbf{Proof of Proposition \ref{Prop-Sufficient-Decrease}.}
 By the definition of $g_k$ and \eqref{Multiplier-update}, we have
	\begin{align}
		& \nabla_\bfw g_k (\bfu^{k+1}) = \bfw^{k+1} + \mu ( \bfw^{k+1} - \bfw^k ) + A^\top \bfz^{k+1} \label{nabla_x_g} \\
		& \nabla_\bfxi g_k (\bfu^{k+1}) = -\bfz^{k+1}. \label{nabla_y_g}
	\end{align}
	These facts will be frequently used in the following proofs.
	
	(i) First, we need to estimate an upper bound for $\| \bfz^{k+1} - \bfz^k \|$. If $| T_{k+1} \cap T_{k} | \geq r$, from \eqref{nabla_x_g}, we have
	\begin{align*}
		A^\top_{:, T_{k+1}\cap T_k} (\bfz^{k+1} - \bfz^{k+1}) = &  [\nwg_k(\bfu^{k+1}) - \nwg_{k-1} (\bfu^{k})  - (\bfw^{k+1} - \bfw^{k})  \\
		& - \mu (\bfw^{k+1} - \bfw^k ) + \mu (\bfw^k - \bfw^{k-1}) ]_{T_{k+1}\cap T_k}
	\end{align*}
	Using Assumption \ref{assum}, we can further estimate
	\begin{align}
		\gamma \| \bfz^{k+1} - \bfz^k \| \leq& \| A^\top_{:, T_{k+1}\cap T_k} (\bfz^{k+1} - \bfz^{k}) \| \leq \| \nabla_{T_{k+1}} g_k (\bfu^{k+1}) \| + \| \nabla_{T_{k}} g_{k-1}(\bfu^{k}) \| \notag \\
		&+ \| \bfw^{k+1} - \bfw^{k} \| + \mu \| \bfw^{k+1} - \bfw^k \| + \mu \| \bfw^k - \bfw^{k-1} \| \notag \\
		\mathop{\leq}\limits^{\eqref{error-metric}}& ( c_1  + \mu + 1) \| \bfw^{k+1} - \bfw^k \| + ( c_1 + \mu ) \| \bfw^k - \bfw^{k-1}  \|. \label{up-bound-z}
	\end{align} 
	If $| T_{k+1} \cap T_{k} | < r$, then taking $|T_{k+1}| = |T_k| = s$ into account, $| T_{k+1} \cap \oT_{k} | = | \oT_{k+1} \cap T_{k} | \geq r$ holds. By \eqref{nabla_x_g} and Assumption \ref{assum}, we can obtain
	\begin{align}
		\gamma\|\bfz^{k+1}\| \leq& \| A^\top_{:, T_{k+1}\cap \oT_k} \bfz^{k+1} \| \notag \\ \mathop{\leq}\limits^{\eqref{nabla_x_g}} & \| \nabla_{T_{k+1}\cap \oT_k} g_k (\bfu^{k+1}) \| + \| [\bfw^{k+1} + \mu (\bfw^{k+1} - \bfw^k)]_{T_{k+1}\cap \oT_k} \| \label{up-bound} \\
		\mathop{\leq}\limits^{\eqref{error-metric}}& c_1 \| \bfw^{k+1} - \bfw^k \| + \|  \bfw^{k+1}_{T_{k+1}\cap \oT_k} \| + \mu \| \bfw^{k+1} - \bfw^k \| \notag \\
		\leq & ( c_1+\mu) \| \bfw^{k+1} - \bfw^k \| + \| [\bfw^{k+1} - \bfw^k]_{T_{k+1}\cap \oT_k} \| + \| \bfw^k_{T_{k+1}\cap \oT_k} \| \notag \\
		\mathop{\leq}\limits^{\eqref{error-metric}} & ( c_1 + \mu + 1) \| \bfw^{k+1} - \bfw^k \| +  c_1 \| \bfw^{k} - \bfw^{k-1} \|. \notag 
	\end{align}
	\begin{align*}
		\gamma\|\bfz^{k}\| \leq& \| A^\top_{:, \oT_{k+1}\cap T_k} \bfz^{k} \| \mathop{\leq}\limits^{\eqref{nabla_x_g}} \| \nabla_{\oT_{k+1}\cap T_k} g_{k-1} (\bfu^{k}) \| + \| [\bfw^{k} + \mu (\bfw^{k} - \bfw^{k-1})]_{\oT_{k+1}\cap T_k} \| \\
		\mathop{\leq}\limits^{\eqref{error-metric}}& c_1 \| \bfw^{k} - \bfw^{k-1} \| +  \| \bfw^{k}_{\oT_{k+1}\cap T_k} \| + \mu \| \bfw^{k} - \bfw^{k-1} \| \\
		\leq & ( c_1+\mu) \| \bfw^{k} - \bfw^{k-1} \| + \| [\bfw^{k+1} - \bfw^k]_{\oT_{k+1}\cap T_k} \| +  \| \bfw^{k+1}_{\oT_{k+1}\cap T_k} \| \\
		\mathop{\leq}\limits^{\eqref{error-metric}} & ( c_1 + \mu ) \| \bfw^{k} - \bfw^{k-1} \| + (  c_1 + 1)  \| \bfw^{k+1} - \bfw^{k} \|.
	\end{align*}
	Adding the two inequalities above yields
	\begin{align*}
		\gamma \| \bfz^{k+1} - \bfz^k \| \leq & \gamma \| \bfz^{k+1} \| + \gamma \| \bfz^k \| \\
		\leq & ( 2c_1 + \mu + 2 ) \| \bfw^{k+1} - \bfw^{k} \| + ( 2c_1 + \mu) \| \bfw^{k} - \bfw^{k-1} \|.
	\end{align*}
	Combining this inequality and \eqref{up-bound-z} leads to
	\begin{align}
		\| \bfz^{k+1} - \bfz^k \| \leq  c_3 \| \bfw^{k+1} - \bfw^{k} \| +  c_4 \| \bfw^{k} - \bfw^{k-1} \|.  \label{up-bound-z1}
	\end{align}
	By using arithmetic mean and quadratic mean inequality, we can obtain
	\begin{align}\label{up-bound-z2}
		\| \bfz^{k+1} - \bfz^k \|^2 \leq 2 c_3^2 \| \bfw^{k+1} - \bfw^{k} \|^2 + 2 c_4^2 \| \bfw^{k} - \bfw^{k-1} \|^2.
	\end{align}
	From the definition of Lyapunov function and the first line of \eqref{error-metric}, we have the following chain of inequalities
	\begin{align*}
		\mathcal{L}_\rho ( \bfu^k, \bfz^k ) - \mathcal{L}_\rho ( \bfu^{k+1},  \bfz^{k +1}) =& \mathcal{L}_\rho ( \bfu^k, \bfz^k ) - \mathcal{L}_\rho ( \bfu^{k+1},  \bfz^{k}) + \mathcal{L}_\rho ( \bfu^{k+1},  \bfz^{k}) - \mathcal{L}_\rho ( \bfu^{k+1},  \bfz^{k +1}) \\
		\mathop{\geq}\limits^{(\ref{error-metric}, \ref{Multiplier-update})} & \frac{\mu}{2} \| \bfw^{k+1} - \bfw^k \|^2 - \frac{1}{\rho} \| \bfz^{k+1} - \bfz^k \|^2 \\
		\mathop{\geq}\limits^{\eqref{up-bound-z2}} & (\frac{\mu}{2} - \frac{2 c_3^2}{\rho}) \| \bfw^{k+1} - \bfw^k \|^2 - \frac{2 c_4^2}{\rho} \| \bfw^{k} - \bfw^{k-1} \|^2.
	\end{align*}
	Then we can further estimate
	\begin{align*}
		\M_k - \M_{k+1} =& \mathcal{L}_\rho ( \bfu^k, \bfz^k ) + \frac{\eta}{2} \| \bfw^{k} - \bfw^{k-1} \|  - \mathcal{L}_\rho ( \bfu^{k+1},  \bfz^{k +1}) - \frac{\eta}{2}  \| \bfw^{k+1} - \bfw^k \| \\
		\geq & (\frac{\mu}{2} - \frac{2 c_3^2}{\rho} - \frac{\eta}{2}) \| \bfw^{k+1} - \bfw^k \|^2 + ( \frac{\eta}{2} - \frac{2 c_4^2}{\rho}) \| \bfw^{k} - \bfw^{k-1} \|^2 \\
		\mathop{\geq}\limits^{\eqref{para_set}} & \frac{\mu}{4} \| \bfw^{k+1} - \bfw^k \|^2
	\end{align*}
	
(ii) From \eqref{nabla_x_g}, we can obtain
	\begin{align*}
		A^\top_{:T_{k+1}} \bfz^{k+1} = [\nwg_k ( \bfu^{k+1} ) - \bfw^{k+1} - \mu ( \bfw^{k+1} - \bfw^k )]_{T_{k+1}}
	\end{align*}
	Using Assumption \ref{assum} and \eqref{error-metric}, 
	we derive 
	%\eqref{up-bound}, we can estimate
	\begin{align*}
		\gamma\| \bfz^{k+1} \| \leq  c_1 \| \bfw^{k+1} - \bfw^k \| + \| \bfw^{k+1} \| + \mu \| \bfw^{k+1} - \bfw^k \| \leq \| \bfw^{k+1} \| + ( c_1 + \mu) \| \bfw^{k+1} - \bfw^k \|.
	\end{align*}
	By using arithmetic and quadratic mean inequality, we have
	\begin{align} 
		\| \bfz^{k+1} \|^2 \leq & \frac{2}{\gamma^2}  \| \bfw^{k+1} \|^2 + \frac{2( c_1 + \mu)^2}{\gamma^2} \| \bfw^{k+1} - \bfw^k \|^2   \label{z-bound}
	\end{align}
	%where the last inequality follows from $\bfw = D \bfw$ and its $\ell_f$-Lipschitz continuity.
	%Noticing that $f$ is a quadratic function, we also have
	%\begin{align} \label{ieq-f}
	%	\frac{1}{2}\| \bfw^{k+1} \|^2 - \frac{1}{\rho \gamma^2} \| \bfw^{k+1} \|^2  \geq (\frac{\sigma_f}{2} - \frac{\ell_f}{\rho\gamma^2}) \| \bfw^{k+1} \|^2
	%\end{align}
	The following chain of inequalities holds by \eqref{lya-des} and the definition of the Lyapunov function
	\begin{align}
		\M_1 \geq & \M_{k+1} = \frac{1}{2}\| \bfw^{k+1} \|^2 + \langle \bfz^{k+1}, A\bfw^{k+1} + \bfone - \bfxi^{k+1} \rangle + \frac{\rho}{2} \| A \bfw^{k+1} + \bfone - \bfxi^{k+1} \|^2 \notag \\
		&+ \frac{\eta}{2} \| \bfw^{k+1} - \bfw^{k} \|^2 + \delta_{\mathbb{S}}(\bfw^{k+1}) + \lambda J(\bfxi^{k+1}) \notag \\
		\geq & \frac{1}{2}\| \bfw^{k+1} \|^2 + \frac{\rho}{2} \| A \bfw^{k+1} + \bfone - \bfxi^{k+1} + \bfz^{k+1}/\rho \|^2 + \frac{\eta}{2} \| \bfw^{k+1} - \bfw^{k} \|^2 - \frac{1}{2 \rho} \| \bfz^{k+1} \|^2 \notag \\
		\mathop{\geq}\limits^{\eqref{z-bound}} & 
		\left(\frac{1}{2} - \frac{1}{\rho\gamma^2} \right) \| \bfw^{k+1} \|^2 
		+ \left( \frac{\eta}{2} - \frac{( c_1 + \mu)^2}{\rho\gamma^2} \right) 
		\| \bfw^{k+1} - \bfw^k \|^2  \notag \\
		& + \frac{\rho}{2} \| A \bfw^{k+1} + \bfone - \bfxi^{k+1} + \frac{1}{\rho} \bfz^{k+1} \|^2 . \label{v-bound}
	\end{align}
	%where the last inequality follows from $\frac{1}{2}\| \bfw \|^2:= \bfw^\top D \bfw/2$ and its $\sigma_f$-strong convexity. 
	Taking \eqref{para_set} into account, both 
	quantities $(1/2 - 1/(\rho\gamma^2))$ and $(\eta/2 - ( c_1+\mu)^2/(\rho\gamma^2))$ are positive. 
	Thus the sequences  $\{ \bfw^{k+1} \}_{k\in \mathbb{N}}$, $\{ \bfw^{k+1} - \bfw^k  \}_{k\in \mathbb{N}}$ and $\{ A \bfw^{k+1} + \bfone - \bfxi^{k+1} + \bfz^{k+1}/\rho \}_{k\in \mathbb{N}}$ are bounded. Then \eqref{z-bound} leads to the boundedness of
	 $\{ \bfz^{k+1} \}_{k\in \mathbb{N}}$. 
	The bound
	\begin{align*}
		\| \bfxi^{k+1} \| \leq \| A \bfw^{k+1} + \bfone - \bfxi^{k+1} + \bfz^{k+1}/\rho \| + \| A \| \| \bfw^{k+1} \| + \| \bfone \| + \| \bfz^{k+1} \|/\rho
	\end{align*}
implies the boundedness of
	$\{ \bfxi^{k+1} \}_{k\in \mathbb{N}}$. 
	Overall, the generated sequence $\{ (\bfu^k, \bfz^k) \}_{k \in \mathbb{N}}$ is bounded.
	
	Finally, let us prove the successive changes of the sequence converge to zero. Actually, \eqref{v-bound} implies $\M_{k+1} \geq 0$ for all $k \in \mathbb{N}$. Combining this and the nonincreasing property \eqref{lya-des}, it follows from the monotone convergence theorem that sequence $\{ \M_{k+1} \}_{k \in \mathbb{N}}$ must be convergent. Therefore, $\lim_{k \to \infty} \| \bfw^{k+1} - \bfw^k \| = 0$. Considering that \eqref{up-bound-z1} holds, we have $\lim_{k \to \infty} \| \bfz^{k+1} - \bfz^k \| = 0$. Finally, by using \eqref{Multiplier-update}, we can obtain
	\begin{align*}
		\| \bfxi^{k+1} - \bfxi^k \| \leq  ( \| \bfz^{k+1} - \bfz^k \| + \| \bfz^{k} - \bfz^{k-1} \| )/\rho + \| A \|\| \bfw^{k+1} - \bfw^k \|.
	\end{align*}
	which implies $\lim_{k \to \infty} \| \bfxi^{k+1} - \bfxi^k \| = 0$ by $\lim_{k \to \infty} \| \bfw^{k+1} - \bfw^k \| = 0$ and $\lim_{k \to \infty} \| \bfz^{k+1} - \bfz^k \| = 0$. \hfill $\blacksquare$
 
\begin{lemma} (Subsequence Convergence) \label{sub-con}  Suppose that Assumption \ref{assum} holds and parameters are chosen as \eqref{para_set}. Let $\{ (\bfu^{k} ; \bfz^k ) \}_{k \in \mathbb{N}}$ be a sequence generated by iPAL, then each of its accumulations points is a {\rm P}-stationary pair of \eqref{Constrained-HM}. Furthermore, $\bfu^*$ is a strict local minimizer of \eqref{Constrained-HM}.
\end{lemma}

{\bf Proof}
	Suppose that $(\bfu^*,\bfz^*)$ is an accumulation point of $\{ (\bfu^{k} ; \bfz^k ) \}_{k \in \mathbb{N}}$. Then there exists a subsequence $\{ (\bfu^{k} ; \bfz^k ) \}_{k \in \K}$ with $\lim_{k \in \K, \ k \to \infty} (\bfu^k, \bfz^k) = (\bfu^*, \bfz^*)$. 
	It follows from \eqref{suc_change} that
	$\{(\bfu^{k+1}, \bfz^{k+1})\}_{k \in K}$ also converges to
	$(\bfu^*, \bfz^*)$. Let us take
	\begin{align*}
		\overline{\bfw}^{k+1}:= \left[ \begin{array}{c}
			[\bfw^{k+1} - \alpha \nabla_\bfw g_k ( \bfu^{k+1} )]_{T_{k+1}} \\
			\bfzero
		\end{array} \right] \ \mbox{and} \ \overline{\bfxi}^{k+1}:= \left[ \begin{array}{c}
			[\bfxi^{k+1} - \beta \nabla_\bfxi g_k ( \bfu^{k+1} )]_{\Gamma_{k+1}} \\
			\bfzero
		\end{array} \right] .
	\end{align*}
	By the definition of $T_{k+1}$ and $\Gamma_{k+1}$, $\overline{\bfw}^{k+1}$ and $\overline{\bfxi}^{k+1}$ actually satisfy
	\begin{align} \label{prox_iterate}
		\overline{\bfw}^{k+1} \in \proj ( \bfw^{k+1} - \alpha \nabla_\bfw g_k ( \bfu^{k+1} ) ) \ \ \mbox{and} \ \
		\overline{\bfxi}^{k+1} \in \Prox_{\beta \lambda J(\cdot)  } ( \bfxi^{k+1} - \beta \nabla_\bfxi g_k ( \bfu^{k+1} ) ).
	\end{align}
	We can also estimate
	\begin{align*}
		&\| \overline{\bfw}^{k+1} - \bfw^{k+1} \| = \| [ \alpha \nabla_{T_{k+1}} g_k (\bfu^{k+1}); \bfw^{k+1}_{\oT_{k+1}} ] \| \leq  \max\{c_1\alpha, c_1\} \| \bfw^{k+1} - \bfw^k \| \\
		&\| \overline{\bfxi}^{k+1} - \bfxi^{k+1} \| = \| [ \beta \nabla_{\Gamma_{k+1}} g_k (\bfu^{k+1}); \bfxi^{k+1}_{\oG_{k+1}} ] \| \leq  \max\{ c_2\beta,c_2 \} \| \bfw^{k+1} - \bfw^k \|^2. 
	\end{align*}
	Considering \eqref{suc_change}, $\lim_{k \to \infty} \| \overline{\bfw}^{k+1} - \bfw^{k+1} \| = \lim_{k \to \infty} \| \overline{\bfxi}^{k+1} - \bfxi^{k+1} \| = 0$ hold. This together with $\lim_{k \to \infty, k\in \K} \|  \bfu^{k+1} - \bfu^* \| = 0$ leads to
	\begin{align} \label{lim-ow}
		\lim_{k \to \infty, k \in \K} \overline{\bfw}^{k+1} = \bfw^* \ \mbox{and} \ \lim_{k \to \infty, k \in \K} \overline{\bfxi}^{k+1} = \bfxi^*.
	\end{align}
	Besides, passing $k \to \infty$ for $k \in \K$ on both sides of \eqref{nabla_x_g}, \eqref{nabla_y_g} and \eqref{Multiplier-update} leads to
	\begin{equation} \label{lim-nab}
		\begin{aligned} 
			&\lim_{k \in, k \to \infty} \nwg_k (\bfu^{k+1}) = \bfw^* + A^\top \bfz^* \\ &\lim_{k \in \K, k \to \infty} \nyxi_k (\bfu^{k+1}) = - \bfz^*. \\
			& A \bfw^* + \bfone - \bfxi^* = 0.
		\end{aligned}
	\end{equation}
	Since \eqref{prox_iterate}, \eqref{lim-ow} and \eqref{lim-nab} hold, it follows from \cite[Theorem 1.25]{rockafellar1976augmented} that $(\bfu^*,\bfz^*)$ will be a {\rm P}-stationary pair satisfying \eqref{P-stat}. Finally, using Theorem \ref{Thm-Stationarity}, we can conclude that $\bfu^*$ is also a strict local minimizer of \eqref{Constrained-HM}.
\hfill $\blacksquare$
 
\noindent\textbf{Proof of Theorem \ref{Thm-Global}}
	Let us first prove that $\lim_{k \to \infty} \bfu^k = \bfu^*$ and $\bfu^*$ is a {\rm P}-stationary point of \eqref{Constrained-HM}. Lemma \ref{sub-con} and Theorem \ref{Thm-Stationarity} indicate that each accumulation point of $\{ \bfu^k \}_{k \in \mathbb{N}}$ is isolated. Moreover, taking \eqref{suc_change} into account, \cite[Proposition 7]{kanzow1999qp} 
	implies $\lim_{k \to \infty} \bfu^k = \bfu^*$.
	We further estimate
	\begin{align*}
		\| \bfw^*_{\oT_{k+1}} \| \leq \|  [\bfw^{k+1} - \bfw^*]_{\oT_{k+1}} \| + \| \bfw^{k+1}_{\oT_{k+1}} \| \mathop{\leq}\limits^{\eqref{error-metric}} \| \bfw^{k+1} - \bfw^*\| +  c_1 \| \bfw^{k+1} - \bfw^k \|
	\end{align*}
	Taking limit as $k \to \infty$ on both sides of above inequality leads to $ \lim_{k \to \infty} \| \bfw^*_{\oT_{k+1}} \| = 0$, which means that $\bfw^*_{\oT_{k+1}} = 0$ when $k$ is sufficiently large. 
	We then have 
	\begin{align} \label{T-S}
		T_{k+1} \supseteq \mathcal{S}^* : = \{ i \in [n]: w^*_i \neq 0 \}.
	\end{align}
	
	Suppose that $\bfz^*$ is a {\rm P}-stationary multiplier associated with $\bfu^*$. Now let us prove $\lim_{k \to \infty} \bfz^k = \bfz^*$. To achieve this goal, we need to give an upper bound for $\| \bfz^{k+1} - \bfz^* \|$. We claim that the following equation holds when $k$ is sufficiently large
	\begin{align*}
		(\bfw^* + A^\top \bfz^*)_{T_{k+1}} = 0.
	\end{align*}
Indeed, if $\|\bfw^*\|_0 = s$, then the $T_{k+1} = \mathcal{S}^*$ follows from \eqref{T-S} and $| T_{k+1} | = s$. This and \eqref{proj_explicit}   further leads to the above equation. If $\|\bfw^*\|_0 < s$, then we have $\bfw^* + A^\top \bfz^* = 0$ by \eqref{proj_explicit}. Moreover, considering that $| T_{k+1} | = s$ holds, we can use Assumption \ref{assum} and \eqref{nabla_x_g} to derive
	\begin{align}
		\gamma \| \bfz^{k+1} - \bfz^* \| \leq & \| [A^\top ( \bfz^{k+1} - \bfz^* )]_{T_{k+1}} \| \leq \| [\nabla_\bfw g_{k} (\bfu^{k+1}) - ( \bfw^{k+1} - \bfw^*) - \mu (\bfw^{k+1} - \bfw^{k})]_{T_{k+1}}  \| \notag \\
		\leq & \| \nabla_{T_{k+1}} g_{k} (\bfu^{k+1}) \| + \| \bfw^{k+1} - \bfw^* \| + \mu \| \bfw^{k+1} - \bfw^{k} \| \notag \\
		\mathop{\leq}\limits^{\eqref{error-metric}} & ( c_1+\mu) \| \bfw^{k+1} - \bfw^{k} \| +  \| \bfw^{k+1} - \bfw^* \| . \label{z-z*}
	\end{align}

	Considering that we have proved $\lim_{k \to \infty} \bfu^k = \bfu^*$, taking limit on both sides of the above inequality yields $\lim_{k \to \infty} \bfz^k = \bfz^*$. Overall, we have verified $\lim_{k \to \infty} (\bfu^k,\bfz^k) = (\bfu^*,\bfz^*)$. Using Lemma \ref{sub-con} and Theorem \ref{Thm-Stationarity}, we can arrive at the desired conclusion. 
	\hfill $\blacksquare$

%%%%%%%%%%%%%%%%%%%%%%%%%%%%%%%%%%%%%%%%%%%%%%%%%%%%%%%
\section{Corollary from global convergence}

\begin{corollary} \label{col_convergence}
		Under the premise of Theorem \ref{Thm-Global}, 
		the following holds.
		
		\begin{itemize}
			\item[(i)] For $k$ is sufficiently large, it holds
			\begin{align}
				&\| \bfw^*_{\oT_{k+1}} \| = 0, \ \left\{ \begin{aligned}
					& T_{k+1} \supseteq \mathcal{S}^*, \ \mbox{if} \ \| \bfw^* \|_0 < s, \\
					& T_{k+1} = \mathcal{S}^*, \ \mbox{if} \ \| \bfw^* \|_0 = s.
				\end{aligned} \right. \label{set-relation} \\
			& \| \bfxi^*_{\oG_{k+1}} \| = 0,~ \| \bfz^*_{\Gamma_{k+1}} \| = 0,~ J(\bfxi^{k+1}) = J(\bfxi^*) \label{y+0-iden} \\
				& 	\| \nabla_{\bfw} g_k (\bfu^{k+1}) \|	\leq   c_5 \| \bfw^{k+1} - \bfw^k \| +  c_6 \| \bfw^{k+1} - \bfw^* \|, \ \mbox{if} \ \| \bfw^* \|_0 < s \label{ng-up-bound}	
			\end{align}
			\item[(ii)] It holds 
			\begin{align*}
				\lim_{k \to \infty} \M_k = \M_*: = \M_{\rho,\eta} ( \bfu^*, \bfz^*, \bfu^* ) = \frac{1}{2}\|\bfw^*\|^2 + \lambda J(\bfxi^*).
			\end{align*}
		\end{itemize}
	\end{corollary}
	
	\bp \
	(i) Formulas \eqref{set-relation} has been proved in Theorem \ref{Thm-Global}. Moreover, $\| \bfxi^*_{\oG_{k+1}} \| = 0$ and $\| \bfz^*_{\Gamma_{k+1}} \| = 0$ can be derived from $\R_2( \bfu^{k+1}  ) \le  c_2 \| \bfw^{k+1} - \bfw^k \|^2$ by a similar procedure as that of $\| \bfw^*_{\oT_{k+1}} \| = 0$. We will first prove $J(\bfxi^{k+1}) = J(\bfxi^*)$ when $k$ is large enough. From the last line of \eqref{error-metric} and the definition of Moreau envelop, we have 
	\begin{align*}
		(\beta/2) \| \nabla_\bfxi g_k (\bfu^{k+1}) \|^2 + \lambda J(\bfxi^{k+1})  \leq & \Phi_{\beta\lambda J(\cdot)} ( \bfxi^{k+1} - \beta \nabla_\bfxi g_k(\bfu^{k+1}) ) + \vartheta_{k} \\
		\leq & \frac{1}{2\beta} \| \bfw^* - ( \bfw^{k+1} - \beta \nabla_\bfxi g_k ( \bfu^{k+1} ) ) \|^2 + \lambda J(\bfxi^*) + \vartheta_{k}.
	\end{align*}
	Taking the superior limits on both sides of the above inequality implies
	\begin{align*}
		\limsup_{k \to \infty} J(\bfxi^{k+1}) \leq J(\bfxi^*).
	\end{align*}
	Combining this with the lower semi-continuity of $J(\cdot)$ leads to \eqref{y+0-iden}. 
	
	Now we will prove \eqref{ng-up-bound}. If $\| \bfw^* \|_0 < s$, then from \eqref{T-S}, $ T_{k+1} \cap \overline{\mathcal{S}}^* \neq \varnothing$ holds. By the definition of $T_{k+1}$, we have the following chain of inequalities
	\begin{align}
		\| [\bfw^{k+1} - \alpha \nwg_k (\bfu^{k+1}) )]_{\oT^{k+1}} \| \leq & | \oT_{k+1} | | [\bfw^{k+1} - \alpha \nwg_k (\bfu^{k+1}) )]_i | \ \mbox{for any} \ i \in T_{k+1} \cap \overline{\mathcal{S}}^*  \notag \\
		\leq & (n-s) \| [\bfw^{k+1} - \alpha \nwg_k (\bfu^{k+1}) )]_{T_{k+1} \cap \overline{\mathcal{S}}^*} \| \notag \\
		\leq & (n-s)( \| [\bfw^{k+1} - \bfw^*]_{T_{k+1} \cap \overline{\mathcal{S}}^*} \| + \alpha \| [  \nwg_k (\bfu^{k+1}) ]_{T_{k+1} \cap \overline{\mathcal{S}}^*} \| ) \notag \\
		\mathop{\leq}\limits^{\eqref{error-metric}} & (n-s)( \| \bfw^{k+1} - \bfw^* \| + \alpha  c_1 \| \bfw^{k+1} - \bfw^k \| ), \label{x-nxg}
	\end{align}
	where the second inequality follows from the fact that $T_{k+1}$ contains the best $s$ largest elements of $\bfw^{k+1} - \alpha \nabla_\bfw g( \bfu^{k+1} )$ in absolute value. Then we can estimate
	\begin{align*}
		\| \nabla_{\overline{T}_{k+1}} g_k (\bfu^{k+1}) \| = & \| [\bfw^{k+1} - ( \bfw^{k+1} - \alpha \nwg_k (\bfu^{k+1}) )]_{\oT^{k+1}} \| / \alpha \\
		\leq & (\| \bfw^{k+1}_{\oT^{k+1}} \| + \| [\bfw^{k+1} - \alpha \nwg_k (\bfu^{k+1}) )]_{\oT^{k+1}} \|) / \alpha \\
		\mathop{\leq}\limits^{\eqref{x-nxg}} & (  c_1 \| \bfw^{k+1} - \bfw^k \| + (n-s) \| \bfw^{k+1} - \bfw^* \| + \alpha  c_1 (n-s) \| \bfw^{k+1} - \bfw^k \|
		)/\alpha.
	\end{align*}
	This result further leads to
	\begin{align} 
		\| \nwg_k (\bfu^{k+1}) \| \leq & \| \nabla_{T_{k+1}} g_k (\bfu^{k+1}) \| + \| \nabla_{\oT_{k+1}} g_k (\bfu^{k+1}) \| \mathop{\leq}\limits^{\eqref{error-metric}}   c_5 \| \bfw^{k+1} - \bfw^k \| +  c_6 \| \bfw^{k+1} - \bfw^* \|,  \notag
	\end{align}
	where $ c_5:=  c_1/\alpha +  c_1(n+1-s)$ and $ c_6:= (n-s)/\alpha$.
	
	(ii) Applying the fact $\lim_{k \to \infty} (\bfu^k,\bfz^{k}) = (\bfu^*,\bfz^*)$, $\delta_{\mathbb{S}}(\bfw^{k+1}) = \delta_{\mathbb{S}}(\bfw^*) =0 $ and \eqref{y+0-iden}, we can derive $\lim_{k \to \infty} \M_k = \M_*$.
	\hfill $\blacksquare$

%%%%%%%%%%%%%%%%%%%%%%%%%%%%%%%%%%%%%%%%%%%%%%%%%%%%%%%
\section{Proof of Theorem \ref{Thm-Convergence-Rate}  on Convergence Rate of iPAL} 

The main steps for convergence rate analysis is as follows.

\begin{itemize}
	\item[$\bullet$] To prove Theorem \ref{Thm-Convergence-Rate} (i), we will first estimate an upper bound  of $\M_{k+1} - \M_*$ (see \eqref{vk-v*}). This, together with the sufficient descent property \eqref{lya-des}, leads to a recursion formula \eqref{recursion}. This will give rise to the linear convergence rate of the Lyapunov function value sequence, see \eqref{r-linear-v}.
	
	\item[$\bullet$] For the linear convergence rate of iterate sequence, we will first investigate the relationship between $\| \bfw^{k+1} - \bfw^* \|$ and $\M_{k+1} - \M_*$ (see \eqref{w-w*<=v-v*}). Then we will use \eqref{r-linear-v} to derive linear convergence rate of $\| \bfw^{k+1} - \bfw^* \|$. Following a similar procedure, we can prove linear convergence rate of $\| \bfxi^{k+1} - \bfxi^* \|$ and $\| \bfz^{k+1} - \bfz^* \|$.
\end{itemize}

\noindent\textbf{Proof of Theorem \ref{Thm-Convergence-Rate}.} 
(i) We start with several inequalities. The first one is a direct computation
\begin{align} \label{up-bd-f}
	\frac{1}{2}\| \bfw^{k+1} \|^2 - \frac{1}{2}\|\bfw^*\|^2 + \langle  \bfw^{k+1}, \bfw^* - \bfw^{k+1} \rangle = -\frac{1}{2} \| \bfw^{k+1} - \bfw^* \|^2.
\end{align}
We now estimate an upper bound for $| \langle  \nyxi_k ( \bfu^{k+1} ), \bfxi^{k+1} - \bfxi^k \rangle |$. Since the sequence boundedness has been proved in Theorem \ref{Prop-Sufficient-Decrease} (ii), we can assume $ \| ( \bfu^{k+1};\bfz^{k+1} ) \| \leq \tau $ for $\tau > 0$. 
Using the fact $\bfxi^{*}_{\oG_{k+1}} = 0$ for sufficiently large $k$, we obtain
\begin{align*}
	| \langle \nabla_{\oG_{k+1}} g_k ( \bfu^{k+1} ), [\bfxi^{k+1} - \bfxi^*]_{\oG_{k+1}} \rangle | =&  | \langle \nabla_{\oG_{k+1}} g_k ( \bfu^{k+1} ), \bfxi^{k+1}_{\oG_{k+1}} \rangle | \leq \| \nabla_{\oG_{k+1}} g_k ( \bfu^{k+1} )\| \| \bfxi^{k+1}_{\oG_{k+1}} \| \\
	\leq & \tau  c_2 \| \bfw^{k+1} - \bfw^k \|^2   \\ 
	| \langle \nabla_{\Gamma_{k+1}} g_k ( \bfu^{k+1} ), [\bfxi^{k+1} - \bfxi^*]_{\Gamma_{k+1}} \rangle | \leq & \| \nabla_{\Gamma_{k+1}} g_k ( \bfu^{k+1} ) \| ( \| \bfxi^{k+1}_{\Gamma_{k+1}} \| + \| \bfxi^*_{\Gamma_{k+1}} \| ) \\
	\leq  & 2\tau c_2 \| \bfw^{k+1} - \bfw^k \|^2 .
\end{align*}
Adding the above inequalities implies
\begin{align} \label{up-bd-y}
	| \langle  \nyxi_k ( \bfu^{k+1} ), \bfxi^{k+1} - \bfxi^k \rangle | \leq 3\tau c_2 \| \bfw^{k+1} - \bfw^k \|^2 
\end{align}
%From the last line of inexact criteria \eqref{error-metric}, we can obtain
%\begin{align*}
%	(\beta/2) \| \nabla_\bfxi g_k (\bfu^{k+1}) \|^2 + \lambda J(\bfxi^{k+1})  \leq & \Phi_{\beta\lambda J(\cdot)} ( \bfxi^{k+1} - \beta \nabla_\bfxi g_k(\bfu^{k+1}) ) +  c_3 \| \bfu^{k+1} - \bfu^k \|^2 \\
%	\leq &  \| \bfw^* - ( \bfw^{k+1} - \beta \nabla_\bfxi g_k ( \bfu^{k+1} ) ) \|^2 /(2\beta) + \lambda J(\bfxi^*).
%\end{align*}
%This further leads to
%\begin{align}
%	\langle  \nabla_\bfxi g_k (\bfu^{k+1}), \bfxi^{k+1} - \bfxi^* \rangle + \lambda J(\bfxi^{k+1}) - \lambda J(\bfxi^*) \leq &  \| \bfxi^{k+1} - \bfxi^* \|^2 /(2\beta) +  c_3 \| \bfu^{k+1} - \bfu^k \|^2. \notag \\
%	\leq &  \| \bfu^{k+1} - \bfu^* \|^2/(2\beta) +  c_3 \| \bfu^{k+1} - \bfu^k \|^2.\label{up-bd-y}
%\end{align}
We shall also derive an upper bound for $| \langle \nabla_{\bfw} g_k (\bfu^{k+1}), \bfw^{k+1} - \bfw^* \rangle |$. If $\| \bfw^* \|_0 = s$,  using $T_{k+1} = \mathcal{S}^*$ and \eqref{error-metric}, we can derive
\begin{align*}
	&| \langle \nabla_{\bfw} g_k (\bfu^{k+1}), \bfw^{k+1} - \bfw^* \rangle |  \\
	= & | \langle \nabla_{T_{k+1}} g_k (\bfu^{k+1}), [\bfw^{k+1} - \bfw^*]_{T_{k+1}} \rangle | \\
	\leq & \| \nabla_{T_{k+1}} g_k (\bfu^{k+1}) \| \| [\bfw^{k+1} - \bfw^*]_{T_{k+1}} \| \leq   c_1 \| \bfw^{k+1} - \bfw^k \| \| \bfw^{k+1} - \bfw^* \|.
\end{align*}
If $\| \bfw^* \|_0 < s$, using $T_{k+1} \supseteq \mathcal{S}^*$, we have 
\begin{align*}
&	\langle \nwg_k ( \bfu^{k+1} ), \bfw^{k+1} - \bfw^* \rangle \\
 =& \langle \nabla_{ T_{k+1} } g_k ( \bfu^{k+1} ), [\bfw^{k+1} - \bfw^*]_{T_{k+1}} \rangle + \langle \nabla_{ \oT_{k+1} } g_k ( \bfu^{k+1} ), [\bfw^{k+1} - \bfw^*]_{\oT_{k+1}} \rangle \\
	=& \langle \nabla_{ T_{k+1} } g_k ( \bfu^{k+1} ), [\bfw^{k+1} - \bfw^*]_{T_{k+1}} \rangle + \langle \nabla_{ \oT_{k+1} } g_k ( \bfu^{k+1} ), \bfw^{k+1}_{\oT_{k+1}} \rangle .
\end{align*}
Then from \eqref{ng-up-bound} and \eqref{error-metric}, the following chain of inequalities holds
\begin{align*}
	&| \langle \nabla_{\bfw} g_k (\bfu^{k+1}), \bfw^{k+1} - \bfw^* \rangle | \\
	\leq&  \| \nabla_{ T_{k+1} } g_k ( \bfu^{k+1} ) \|  \| [\bfw^{k+1} - \bfw^*]_{T_{k+1}} \| +  \| \nabla_{ \oT_{k+1} } g_k ( \bfu^{k+1} ) \| \| \bfw^{k+1}_{\oT_{k+1}} \|  \\
	\leq &  c_1 \| \bfw^{k+1} - \bfw^k \| \| \bfw^{k+1} - \bfw^* \| +  c_1( c_5 \| \bfw^{k+1} - \bfw^k \| +  c_6 \| \bfw^{k+1} - \bfw^* \|) \| \bfw^{k+1} - \bfw^k \| \\
	\leq &  c_1 c_5 \| \bfw^{k+1} - \bfw^k \|^2 + (  c_1 c_6 +  c_1 ) \| \bfw^{k+1} - \bfw^k \| \| \bfw^{k+1} - \bfw^* \| .
\end{align*}
These two cases lead to
\begin{equation} \label{up-bd-x}
	| \langle \nabla_{\bfw} g_k (\bfu^{k+1}), \bfw^{k+1} - \bfw^* \rangle | \leq  c_1 c_5 \| \bfw^{k+1} - \bfw^k \|^2 + (  c_1 c_6 +  c_1 ) \| \bfw^{k+1} - \bfw^k \| \| \bfw^{k+1} - \bfw^* \|
\end{equation}
Now let us consider $\M_{k+1} - \M_*$. For sufficiently large $k$, using definition of Lyapunov function, $J(\bfxi^{k+1}) = J(\bfxi^*)$ and $A \bfw^* + \bfone - \bfxi^* = 0$, we have
\begin{align*}
	&\M_{k+1} - \M_*\\ =  &\frac{1}{2}\| \bfw^{k+1} \|^2 - \frac{1}{2}\|\bfw^*\|^2 + \langle \bfz^{k+1}, A\bfw^{k+1} + \bfone - \bfxi^{k+1} \rangle + \frac{1}{2\rho} \| \bfz^{k+1} - \bfz^k \|^2 + \frac{\eta}{2} \| \bfw^{k+1} - \bfw^k \|^2 \\
	=&\frac{1}{2}\| \bfw^{k+1} \|^2 - \frac{1}{2}\|\bfw^*\|^2 + \langle \bfz^{k+1}, A\bfw^{k+1} + \bfone - \bfxi^{k+1} - ( A\bfw^* + \bfone - \bfxi^* ) \rangle \\
	 &  + \frac{1}{2\rho} \| \bfz^{k+1} - \bfz^k \|^2 + \frac{\eta}{2} \| \bfw^{k+1} - \bfw^k \|^2 \\  
	= & \frac{1}{2}\| \bfw^{k+1} \|^2 - \frac{1}{2}\|\bfw^*\|^2 + \langle A^\top \bfz^{k+1}, \bfw^{k+1} - \bfw^* \rangle - \langle \bfz^{k+1}, \bfxi^{k+1} - \bfxi^* \rangle \\
	& + \frac{1}{2\rho} \| \bfz^{k+1} - \bfz^k \|^2 + \frac{\eta}{2} \| \bfw^{k+1} - \bfw^k \|^2 .
\end{align*}
Applying \eqref{nabla_x_g} and \eqref{nabla_y_g}, we can further derive
\begin{align*}
	&\M_k - \M_* \\
	&=  \frac{1}{2}\| \bfw^{k+1} \|^2 - \frac{1}{2}\|\bfw^*\|^2 - \langle \bfw^{k+1}, \bfw^{k+1} - \bfw^* \rangle +  \langle \nyxi_k (\bfu^{k+1}), \bfxi^{k+1} - \bfxi^* \rangle  \\
	& + \langle \nwg_k (\bfu^{k+1}), \bfw^{k+1} - \bfw^* \rangle - \mu \langle \bfw^{k+1} - \bfw^k, \bfw^{k+1} - \bfw^* \rangle +\frac{1}{2\rho} \| \bfz^{k+1} - \bfz^k \|^2 + \frac{\eta}{2} \| \bfw^{k+1} - \bfw^k \|^2 
\end{align*}
Then we can use the previous inequalities \eqref{up-bound-z2}, \eqref{up-bd-f}, \eqref{up-bd-y}, \eqref{up-bd-x} as well as the fact $- \mu \langle \bfw^{k+1} - \bfw^k, \bfw^{k+1} - \bfw^* \rangle \leq \mu \| \bfw^{k+1} - \bfw^k \|  \| \bfw^{k+1} - \bfw^* \|$ to obtain
\begin{align}
	&\M_k - \M_* \notag \\
	\leq & - \frac{1}{2} \| \bfw^{k+1} - \bfw^* \|^2 + \underbrace{(  c_1  c_6 +  c_1 + \mu )}_{:=  c_7} \| \bfw^{k+1} - \bfw^k \| \| \bfw^{k+1} - \bfw^* \|  \notag \\
	& + \underbrace{( 3\tau c_2 +  c_1  c_5 + \frac{\eta}{2} + \frac{ c_3^2}{\rho}  )}_{:= c_8} \| \bfw^{k+1} - \bfw^k \|^2 + \frac{ c_4^2}{\rho} \| \bfw^k - \bfw^{k-1} \|^2 \notag \\
	=& -\frac{1}{2} (\| \bfw^{k+1} - \bfw^* \| - c_7 \| \bfw^{k+1} - \bfw^k \| )^2 + (  c_8 + \frac{ c_7^2}{2} ) \| \bfw^{k+1} - \bfw^k \|^2 + \frac{ c_4^2}{\rho} \| \bfw^k - \bfw^{k-1} \|^2 \notag \\
	%	 \leq & \underbrace{(  c_3 + \frac{ c_5}{2} + \frac{ c_3^2}{\rho} + \frac{\mu^2}{4} )}_{:= c_7} \| \bfu^{k+1} - \bfu^k \|^2 + \frac{ c_4^2}{\rho} \| \bfu^{k} - \bfu^{k-1} \|^2   + \underbrace{( \frac{ c_5}{2} +  c_6 + \frac{1}{2\beta} + \frac{\eta}{2} + 1 )}_{:=  c_8} \| \bfu^{k+1} - \bfu^* \|^2 \notag \\
	\leq & \tau_1 ( \| \bfw^{k+1} - \bfw^k \|^2 + \| \bfw^{k} - \bfw^{k-1} \|^2 ), \label{vk-v*}
\end{align} 
where $\tau_1 : =  c_8 +  c_7^2/2$. Now taking the descent property \eqref{lya-des} into account, we can estimate
\begin{align*}
	(\M_{k-1} - \M_*) - (\M_{k+1} - \M_*) = & \M_{k-1} - \M_k + \M_k - \M_{k+1} \\
	 \geq & \frac{\mu}{4} ( \| \bfw^{k+1} - \bfw^k \|^2 + \| \bfw^k - \bfw^{k-1} \|^2 ).
\end{align*}
Combining this with \eqref{vk-v*} leads to 
\begin{align} \label{recursion}
	\M_{k+1} - \M_* \leq \frac{1}{1 + \mu/(4\tau_1)} ( \M_{k-1} - \M_* ).
\end{align}
This means that there exists a sufficiently large $k^*$ such that \eqref{r-linear-v} holds for constants
\begin{align*}
	q: = \sqrt{\frac{1}{1 + \mu/(4\tau_1)}} \quad \mbox{and} \quad
	  c_m := (1/q)^{k^*} ( \M_0 - \M_* ).
\end{align*} 

(ii) Suppose that index $k$ is sufficiently large. It follows from \eqref{vk-v*} that
\begin{align*}
	\sqrt{\M_{k+1} - \M_*} \leq \sqrt{\tau_1 (\| \bfw^{k+1} - \bfw^k \|^2 + \| \bfw^{k} - \bfw^{k-1} \|^2)} \leq \sqrt{\tau_1} (\| \bfw^{k+1} - \bfw^k \| + \| \bfw^{k} - \bfw^{k-1} \|) .
\end{align*}
Using this relationship and the concavity of $\sqrt{(\cdot)}$, we can obtain
\begin{align*}
	\varepsilon_k := \sqrt{\M_{k} - \M_*} - \sqrt{\M_{k+1} - \M_*} \geq \frac{\M_{k} - \M_{k+1}}{2\sqrt{\M_k - \M_*}} \geq \frac{\mu\| \bfw^{k+1} - \bfw^k \|^2}{8\sqrt{\tau_1}( \| \bfw^{k} - \bfw^{k-1} \| + \| \bfw^{k-1} - \bfw^{k} \| )},
\end{align*}
This further leads to
\begin{align*}
	\| \bfw^{k+1} - \bfw^k \| \leq & \left( \frac{8\sqrt{\tau_1}\varepsilon_k}{\mu} (\| \bfw^{k} - \bfw^{k-1} \| + \| \bfw^{k-1} - \bfw^{k-2} \|) \right)^{\frac{1}{2}} \\
	\leq & \frac{1}{4} ( \| \bfw^{k} - \bfw^{k-1} \| + \| \bfw^{k-1} - \bfw^{k-2} \| ) + \frac{8\sqrt{\tau_1}}{\mu} \varepsilon_k .
\end{align*}
Let us consider the sum of the above terms from $\ell = k+2$ to $\ell = \widetilde{k}$.
\begin{align*}
	\sum_{\ell = k+2}^{\widetilde{k}} \| \bfw^{\ell + 1} - \bfw^\ell \| \leq & \frac{1}{4} \sum_{\ell = k+2}^{\widetilde{k}} \| \bfw^\ell - \bfw^{\ell - 1} \| + \frac{1}{4} \sum_{\ell = k+2}^{\widetilde{k}} \| \bfw^{\ell - 1} - \bfw^{\ell - 2} \| + \frac{8\sqrt{\tau_1}}{\mu} \sum_{\ell = k+2}^{\widetilde{k}} \varepsilon_\ell \\
	\leq & \frac{1}{4} \sum_{\ell = k+2}^{\widetilde{k}} \| \bfw^{\ell+1} - \bfw^{\ell} \| + \frac{1}{4} \sum_{\ell = k+2}^{\widetilde{k}} \| \bfw^{\ell + 1} - \bfw^{\ell} \| + \frac{8\sqrt{\tau_1}}{\mu} \sum_{\ell = k+2}^{\widetilde{k}} \varepsilon_\ell \\
	& + \frac{1}{2} \| \bfw^{k+2} - \bfw^{k+1}  \|  + \frac{1}{4} \| \bfw^{k+1} - \bfw^{k}  \|
\end{align*}  
After some algebraic manipulation, we have
\begin{align*}
	\sum_{\ell = k}^{\widetilde{k}} \| \bfw^{\ell + 1} - \bfw^\ell \| \leq & \frac{3}{2} \| \bfw^{k+1} - \bfw^k \| + 2 \| \bfw^{k+2} - \bfw^{k+1} \| + \frac{16\sqrt{\tau_1}}{\mu} \sum_{\ell = k+2}^{\widetilde{k}} \varepsilon_\ell \\
	\leq & \frac{3}{\sqrt{\mu}} \sqrt{\M_k - \M_{k+1}} + \frac{4}{\sqrt{\mu}} \sqrt{\M_{k+1} - \M_{k+2}} + \frac{16\sqrt{\tau_1}}{\mu} \sqrt{\M_{k+2} - \M_*} \\
	\leq & (\frac{7}{\sqrt{\mu}} + \frac{16\sqrt{\tau_1}}{\mu})  \sqrt{\M_{k} - \M_*}.
\end{align*}
The taking $\widetilde{k} \to \infty$ for above inequality yields
\begin{align} \label{w-w*<=v-v*}
	\| \bfw^{k} - \bfw^*  \| \leq \sum_{\ell = k}^{\infty} \| \bfw^{\ell + 1} - \bfw^\ell \| \leq \underbrace{(7/\sqrt{\mu} + 16\sqrt{\tau_1}/\mu)}_{:= \tau_{2}}  \sqrt{\M_{k} - \M_*}.
\end{align}
Since \eqref{r-linear-v} holds, we can derive the R-linear convergence rate for $\{ \bfw^k \}_{k \in \mathbb{N}}$ in \eqref{r-linear-wz} with constant $ c_w:= \tau_2 \sqrt{ c_m}$. 

Next, we will prove the R-liner convergence rate of $\{ \bfz^{k} \}_{k \in \mathbb{N}}$. From \eqref{z-z*}, we can estimate
\begin{align*}
	\gamma \| \bfz^k - \bfz^* \| 
	\mathop{\leq} & ( c_1+\mu) \| \bfw^k - \bfw^{k-1} \| +  \| \bfw^k - \bfw^* \| \\ \mathop{\leq}\limits^{\eqref{lya-des},\eqref{w-w*<=v-v*}} & (2( c_1 + \mu)/\sqrt{\mu}) \sqrt{\M_{k-1} - \M_k} + \tau_2 \sqrt{\M_{k} - \M_*} \\
	\leq & \underbrace{( 2( c_1 + \mu)/\sqrt{\mu} + \tau_2  )}_{:=\tau_3} \sqrt{\M_{k-1} - \M_*}  .
\end{align*}
%By \eqref{w-w*<=v-v*} and descent property \eqref{lya-des}, we can obtain
By using \eqref{r-linear-v}, we can arrive at $\| \bfz^{k} - \bfz^* \| \leq  c_z \sqrt{q}^k$ with constant $ c_z := (\tau_3/\gamma)\sqrt{ c_m/q}$. Finally, we prove the linear convergence rate of $\{ \bfxi^k \}_{k \in \mathbb{N}}$. Using \eqref{Multiplier-update} and $A\bfw^* + \bfone - \bfxi^* = 0$, we can estimate
%\begin{align*}
%	\| \bfxi^{k} - \bfxi^* \| \leq & \| A \| \| \bfw^{k} - \bfw^* \| +  \| \bfz^{k} - \bfz^{k-1} \| / \rho \mathop{\leq}\limits^{\eqref{up-bound-z1}} \| A \| \| \bfw^{k} - \bfw^* \| + \frac{ c_3}{\rho} \| \bfw^{k} - \bfw^{k-1} \| + \frac{ c_4}{\rho} \| \bfw^{k-1} - \bfw^{k-2} \| \\
%	\leq & \|A\| \tau_2 \sqrt{ \M_k - \M_* } + \frac{2 c_3}{\rho\sqrt{\mu}} \sqrt{ \M_{k-1} - \M_k } + \frac{2 c_4}{\rho\sqrt{\mu}} \sqrt{ \M_{k-2} - \M_{k-1} } \\
%	 \leq & \underbrace{ \left( \| A \| \tau_2 + \frac{2( c_3 +  c_4)}{\rho\sqrt{\mu}} \right)}_{:= \tau_4} \sqrt{ \M_{k-2} - \M_* }
%\end{align*}
\begin{align*}
	\| \bfxi^{k} - \bfxi^* \| \leq & \| A \| \| \bfw^{k} - \bfw^* \| +  \| \bfz^{k} - \bfz^{k-1} \| / \rho \mathop{\leq}\limits^{\eqref{up-bound-z1}} \| A \| \| \bfw^{k} - \bfw^* \| \\
	& + ( c_3/\rho) \| \bfw^{k} - \bfw^{k-1} \| + ( c_4/\rho) \| \bfw^{k-1} - \bfw^{k-2} \| \\
	\mathop{\leq}\limits^{(\ref{w-w*<=v-v*}, \ref{lya-des})} & \|A\| \tau_2 \sqrt{ \M_k - \M_* } + \frac{2 c_3}{\rho\sqrt{\mu}} \sqrt{ \M_{k-1} - \M_k } + \frac{2 c_4}{\rho\sqrt{\mu}} \sqrt{ \M_{k-2} - \M_{k-1} } \\
	\leq & \underbrace{ \left( \| A \| \tau_2 + \frac{2( c_3 +  c_4)}{\rho\sqrt{\mu}} \right)}_{:= \tau_4} \sqrt{ \M_{k-2} - \M_* } .
\end{align*}
This means $\| \bfz^{k} - \bfz^* \| \leq  c_z \sqrt{q}^k$ can be verified with $ c_z = \tau_4 \sqrt{ c_m/q^2}$.
\hfill $\blacksquare$

\section{Proofs on Convergence Properties of PGN}

First we explain the general ideas for the proof of Theorem \ref{Thm-PGN-Global}.

\begin{itemize}
	\item To prove Theorem \ref{Thm-PGN-Global} (i), we first show the objective function $G$ enjoys sufficient descent \eqref{half_des} on the proximal gradient iterate $\bfu^{j+1/2}$. Then if Newton step is accepted, $G$ also enjoys the sufficient descent \eqref{Newton-Condition}. 
	These results lead to \eqref{des-G}. 
	We then show the convergence of $\{ G(\bfu^j) \}_{j \in \mathbb{N}}$,
	 which further implies \eqref{inn-suc-chan}.
	
	\item The procedure to prove (ii) is similar to the global convergence of iPAL. First, we show that the sequence $\{ \bfu^j \}_{j \in \mathbb{N}}$ is bounded. 
	Second, the boundedness ensures the existence of accumulated points and we will prove each of them is a {\rm P}-stationary point of subproblem \eqref{inner-sub}. Finally, we will utilize \cite[Proposition 7]{kanzow1999qp} to show the whole sequence is convergent. Again, this proposition requires \eqref{inn-suc-chan} and isolatedness of the {\rm P}-stationary points, which we will show in the following proof.
	
	\item[$\bullet$] For (iii), since we have proved the sequence $\{\bfu^j\}_{\in \mathbb{N}}$ converges to a {\rm P}-stationary point of \eqref{inner-sub} and the inexact criteria is just an approximation of the {\rm P}-stationary system. Then the iterate can satisfy the inexact criteria after finite steps.
\end{itemize}

\noindent\textbf{Proof of Theorem \ref{Thm-PGN-Global}}
	(i) Let us first prove the descent property of $G$. From \eqref{w-half}, and the definition of projection and proximal operator, we have
	\begin{align*}
		\left\{ \begin{aligned}
			& \frac{1}{2 \alpha} \| \bfw^{j+1/2} - ( \bfw^j - \alpha \nabla_\bfw g (\bfu^j) ) \|^2 \leq \frac{\alpha}{2} \| \nabla_\bfw g (\bfu^j) \|^2 \\
			& \frac{1}{2 \beta} \| \bfxi^{j+1/2} - ( \bfxi^j - \beta \nabla_\bfxi g (\bfu^j) ) \|^2 + \lambda J(\bfxi^{j+1/2}) \leq \frac{\beta}{2} \| \nabla_\bfxi g (\bfu^j) \|^2 + \lambda J(\bfxi^{j}) .
		\end{aligned}  \right.
	\end{align*}
	By some simple algebraic manipulation, the following inequalities can be deduced 
	\begin{align}\label{prox1}
		\left\{ \begin{aligned}
			& \langle \nabla_\bfw g( \bfu^j ), \bfw^{j+1/2} - \bfw^j \rangle \leq - \frac{1}{2\alpha} \| \bfw^{j+1/2} - \bfw^j \|^2 \\
			& \langle \nabla_\bfxi g( \bfu^j ), \bfxi^{j+1/2} - \bfxi^j \rangle + \lambda J(\bfxi^{j+1/2}) - \lambda J(\bfxi^j) \leq - \frac{1}{2\beta} \| \bfxi^{j+1/2} - \bfxi^j \|^2 .
		\end{aligned} \right.
	\end{align}
	Using the descent lemma \cite[Lemma 5.7]{beck2017first} on function $G$ yields
	\begin{align} \label{des-lem-xy}
		\begin{aligned}
			& g(\bfu^{j+1/2}) \leq g ( \bfu^j ) + \langle \nabla g( \bfu^j ), \bfu^{j+1/2} - \bfu^j \rangle + \frac{\ell_g}{2} \| \bfu^{j+1/2} - \bfu^j \|^2 .
		\end{aligned} 
	\end{align}
	%\begin{align} \label{des-lem-xy}
	%	\left\{ \begin{aligned}
	%		& g(\bfw^{j+1/2}, \bfxi^j) \leq g ( \bfu^j ) + \langle \nabla_\bfw g( \bfu^j ), \bfw^{j+1/2} - \bfw^j \rangle + \frac{\ell_g}{2} \| \bfw^{j+1/2} - \bfw^j \|^2  \\
	%		& g(\bfu^{j+1/2}) \leq g ( \bfw^{j+1/2}, \bfxi^j ) + \langle \nabla_\bfxi g( \bfw^{j+1/2}, \bfxi^j ), \bfxi^{j+1/2} - \bfxi^j \rangle + \frac{\ell_g}{2} \| \bfxi^{j+1/2} - \bfxi^j \|^2 
	%	\end{aligned} \right.
	%\end{align}
	Taking $\delta_{\mathbb{S}}(\bfw^{j+1/2}) = \delta_{\mathbb{S}}( \bfw^j )$ into account and adding \eqref{prox1} and \eqref{des-lem-xy}, 
	we obtain
	%\begin{align*} 
	%	\left\{ \begin{aligned}
	%		& G( \bfw^{j+1/2}, \bfxi^j ) \leq G( \bfu^j ) + ( \frac{\ell_g}{2} - \frac{1}{2\alpha} ) \| \bfw^{j+1/2} - \bfw^j \|^2  \\
	%		& G ( \bfu^{j+1/2} ) \leq G( \bfw^{j+1/2}, \bfxi^j ) + ( \frac{\ell_g}{2} - \frac{1}{2\beta} ) \| \bfxi^{j+1/2} - \bfxi^j \|^2
	%	\end{aligned} \right.
	%\end{align*}
	%The descent property for the gradient step can be derived by this result
	\begin{align} \label{half_des}
		G(\bfu^j) - G(\bfu^{j+1/2}) \geq \zeta \| \bfu^{j+1/2} - \bfu^j \|^2.
	\end{align}
	Then in each case of the update step \eqref{Newton-Condition}, we have
	\begin{align}\label{New-des}
		G ( \bfu^{j+1/2} ) - G ( \bfu^{j+1}) \geq (\sigma_g/4) \|  \bfu^{j+1/2} - \bfu^{j+1} \|^2
	\end{align}
	Adding the above two inequalities directly leads to \eqref{des-G}. 
	
	Since $g$ is strongly convex, and $\delta_{\mathbb{S}}(\cdot)$ and $J(\cdot)$ are lower bounded, we can conclude $G$ is also bounded below. Then $\{ G(\bfu^j) \}_{j \in \mathbb{N}}$ is a nonincreasing and bounded sequence, which implies the sequence is convergent. This result together with \eqref{des-G} yields \eqref{inn-suc-chan}.
	
	(ii) We will prove the global convergence of $\{ \bfu^j \}_{j\in \mathbb{N}}$ according to the three steps mentioned at the beginning of this section.
	
	\textbf{Step 1.} Since $g$ is strongly convex, it is also coercive, i.e. $\lim_{\| \bfu \| \to \infty} g( \bfu ) = \infty$. Combining this with the lower boundedness and lower semi-continuity of $\delta_{\mathbb{S}}(\cdot)$ and $\lambda J(\cdot)$ implies that $G$ is lower semi-continuous and coercive. It follows from \cite[Theorem 4.10]{mordukhovich2013easy} that $\{ 
	\bfu^{j} \}_{j \in \mathbb{N}}$ is bounded. 
	
	\textbf{Step 2.} Suppose that $\hbu$ is an accumulation point of $\{ \bfu^j \}_{j \in \mathbb{N}}$. Then there exists a subsequence $\{ \bfu^j \}_{j \in \mathcal{J}}$ converging to $\hbu$. 
	It follows from the continuous differentiability of $g$ that
	\begin{align*}
		\lim_{j \in \J} \bfw^{j} - \alpha \nabla_\bfw g_k (\bfu^j)  =  \hbu - \alpha \nabla_\bfw g_k (\hbu) \ \mbox{and} \ \lim_{j \in \J} \bfxi^{j} - \beta \nabla_\bfxi g_k (\bfu^j)  =  \hbxi - \alpha \nabla_\bfw g_k (\hbu).  
	\end{align*}
	Since $\| \bfu^{j+1/2} - \hbu \| \leq \| \bfu^{j+1/2} - \bfu^j \| + \| \bfu^j - \hbu \|$, by using $ \lim_{j \to \infty, j\in \J} \bfu^j = \hbu $ and \eqref{inn-suc-chan}, we have 
	\begin{align*}
		\lim_{j \to \infty, j\in \J} \bfu^{j+1/2} = \hbu
	\end{align*}
	Finally, it follows from \cite[Theorem 1.25]{rockafellar1976augmented} that $\hbu$ must satisfy system \eqref{P-stat-sub}.
	
	\textbf{Step 3.} Let $\hbu$ be an accumulation point of $\{ \bfu^j \}_{j \in \mathbb{N}}$. We define $\widehat{\mathcal{I}}_- := \{ i \in [m] : \widehat{\bfxi}_i \leq 0 \}$ and select $\widehat{T} \in \mathbb{T}:= \{ T : T \supseteq \widehat{\mathcal{S}}, \ | T | = s \}$. We consider the following convex programming.
	%\begin{align*}
	%	\min_{\bfu:=(\bfw,\bfxi)} g (\bfu) \quad s.t. \quad \bfxi_{\wmI_-} \leq 0, \quad \bfw_{\overline{\wT}} = 0
	%\end{align*} 
	\begin{align} \label{con-pro}
		\min_{\bfu:=(\bfw,\bfxi)} g (\bfu) \quad s.t. \quad \bfxi_{\hmI_-} \leq 0, \quad \bfw_{i} = 0, \  i \notin \hT.
	\end{align} 
	Since the objective function $g$ is strongly convex and the constraints are linear, if a point satisfies the following KKT system, then it must be the unique global minimizer of the above convex programming 
	\begin{align*}
		\left\{ \begin{aligned}
			& [\nabla_\bfw g (\bfu)]_{\widehat{T}} = 0, \quad \bfw_i = 0,\ i \notin \widehat{T} \\
			& \bfxi_{\hmI_-} \leq 0, \quad  \bfz_\xi \geq 0, \quad \langle \bfxi_{\hmI_-}, \bfz_\xi \rangle = 0 \\
			& \bfz_\xi = - [\nabla g(\bfu)]_{\hmI_-}, \quad  [\nabla_\bfxi g (\bfu)]_{i} = 0, \ i \notin \hmI_-
		\end{aligned} \right.
	\end{align*}
	Considering that the accumulation point $\hbu$ satisfies \eqref{P-stat-sub}, Lemmas \ref{Lemma-Proj} and \ref{Lemma-Prox} imply that
	 $\hbu$ must satisfy the above KKT system, and thus it is the unique global minimizer of \eqref{con-pro}. Since the numbers of the choices of $\hmI$ and $\hT$ are finite, the number of accumulation for $\{ \bfu^j \}_{j \in \mathbb{N}}$ is also finite. Therefore, each accumulation point must be isolated. Finally, taking \eqref{inn-suc-chan} into account, it follows from \cite[Proposition 7]{kanzow1999qp} that the whole sequence $\{ \bfu^j \}_{j \in \mathbb{N}}$ must converge to $\hbu$.
	
	(iii) By the assertion (ii), $\{ \bfu^{j} \}_{j \in \mathbb{N}}$ must converge to a {\rm P}-stationary point $\hbu$ of \eqref{inner-sub}. Let us first prove 
	\begin{align} \label{y+-iden}
		\lim_{j \to \infty} J(\bfxi^{j+1/2}) = \lim_{j \to \infty} J(\bfxi^{j+1}) = J(\hbxi) .
	\end{align}
	From \eqref{w-half} and the definition of proximal operator, we have
	\begin{align*}
		\frac{1}{2\beta} \| \bfxi^{j+1/2} - \bfxi^j + \beta \nabla_\bfxi g( \bfu^j ) \|^2 + \lambda J(\bfxi^{j+1/2}) \leq  \frac{1}{2\beta} \| \hbxi - \bfxi^j + \beta \nabla_\bfxi g( \bfu^j ) \|^2 + \lambda J(\hbxi) .
	\end{align*}
	Taking the superior limits on both sides of above inequality implies $\limsup_{j \to \infty} J(\bfxi^{j+1/2}) \leq J(\hbxi).$ Combining this with the lower semi-continuity of $J(\cdot)$ leads to $\lim_{j \to \infty} J(\bfxi^{j+1/2}) = J(\hbxi) $. 
	
	From \eqref{New-des} and the fact $\delta_{\mathbb{S}} (\bfw^{j+1}) = \delta_{\mathbb{S}} (\bfw^{j+1/2}) = 0$, we have
	\begin{align*}
		g(\bfu^{j+1}) + \lambda J(\bfxi^{j+1}) + (\sigma_g/4) \| \bfu^{j+1/2} - \bfu^{j+1} \| \leq g(\bfu^{j+1/2}) + \lambda J(\bfxi^{j+1/2}) .
	\end{align*} 
	Taking the superior limits on both sides of the above inequality, we have 
	\begin{align*}
		\limsup_{j \to \infty} J(\bfxi^{j+1}) \leq \limsup_{j \to \infty} J(\bfxi^{j+1/2}) = J(\hbxi) .
	\end{align*}
	This together with lower semi-continuity of $J(\cdot)$ leads to $\lim_{j \to \infty} J(\bfxi^{j+1}) = J(\hbxi)$.
	
We now show $ \lim_{j \to \infty} \mathcal{R}_i(\bfu^j) = 0$ for $i = 1,2,3$. The first line of \eqref{error-metric} directly follows from \eqref{des-G} and \eqref{w-half}.  Furthermore, we have
	\begin{align} \label{lim1}
		\begin{aligned}
			&  \mathcal{R}_1(\bfu^j) =  \| [\alpha\nabla_{T_{j}} g_k (\bfu^{j}); \bfw^{j}_{\overline{T}_{j}}]  \| \leq \max\{ \alpha,1 \} \| \bfw^{j+1/2} - \bfw^{j} \| \\
			& \mathcal{R}_2(\bfu^j) = \| [\beta\nabla_{\Gamma_{j}} g_k (\bfu^{j}); \bfxi^{j}_{\overline{\Gamma}_{j}}] \| \leq \max\{ \beta,1 \} \| \bfxi^{j+1/2} - \bfxi^{j} \| ,
		\end{aligned}
	\end{align}
	where $T_{j}$ and $\Gamma_{j}$ are corresponding index sets for the $j$-th identification step. Then we derive $\lim_{j\to \infty} \R_1 (\bfu^j) = \lim_{j\to \infty} \R_2 (\bfu^j) = 0$ by using \eqref{inn-suc-chan}. Applying the definition of Moreau envelop and \eqref{w-half} yields
	\begin{align} 
		& \lim_{j \to \infty} \mathcal{R}_3(\bfu^j) =	\lim_{j \to \infty} \frac{\beta}{2} \| \nabla_\bfxi g_k (\bfu^j) \|^2 + \lambda J(\bfxi^j) - \Phi_{\beta\lambda J(\cdot)} ( \bfxi^j - \beta \nabla_\bfxi g_k(\bfu^j) ) \notag \\
		=& \lim_{j \to \infty} \frac{\beta}{2} \| \nabla_\bfxi g_k (\bfu^j) \|^2 + \lambda J(\bfxi^j) -\frac{1}{2\beta} \| \bfxi^{j+1/2} - \bfxi^j + \beta \nabla_{\bfxi} g_k (\bfu^j) \|^2  - \lambda J(\bfxi^{j+1/2}) \notag \\
		= & \lim_{j \to \infty} - \frac{1}{2\beta} \| \bfxi^{j+1/2} - \bfxi^j \|^2 - \langle \nabla_\bfxi g_k ( \bfu^j ), \bfxi^{j+1/2} - \bfxi^j \rangle + \lambda J(\bfxi^j)   - \lambda J(\bfxi^{j+1/2}) \mathop{=}\limits^{\eqref{y+-iden}} 0 \label{lim2} .
	\end{align}
	Meanwhile, we can derive $\lim_{j \to \infty} \| \bfw^{k,j} - \bfw^k \| = \| \hbw - \bfw^k \| = \| \hbw - \bfw^{k,0} \| \neq 0 $. Combing this with \eqref{lim1} and \eqref{lim2}, we arrive at the desired conclusion.
\hfill $\blacksquare$

Next we will prove the local quadratic convergence rate of PGN. The main ideas for this proof are presented as follows.
\begin{itemize}
	\item[$\bullet$] The changeable index sets $T_j$ and $\Gamma_j$ in \eqref{subspace-newton} brings difficulties for convergence rate analysis. We will show that they exactly contain nonzero elements of the solution (active sets) after finite iterate 
	(see Lemma~\ref{finite_iden} below).
	
	\item[$\bullet$] Newton step plays a crucial role to ensure quadratic convergence rate. We will prove the update condition \eqref{Newton-Condition} always holds after finite iterations.
	Thus Newton step is accepted (see the first part of Theorem~\ref{qua_convergence}).
	
	\item[$\bullet$] Once the above two points prove to be true, the Newton iteration will be always performed on a fixed subspace. Then considering the strong convexity of $g$, the local quadratic convergence of PGN just follows from classical theory. It is also noteworthy that there is a gradient step before Newton step and both of them are performed on the same subspace 
	\begin{equation*}
		\bfu^j \ \longrightarrow \bfu^{j+1/2} \ (\mbox{gradient step}) \ \longrightarrow \bfu^{j+1} \ (\mbox{Newton step}) .
	\end{equation*}
We shall show the gradient iteration will not influence the quadratic convergence rate (see, the second part of Theorem \ref{qua_convergence}). 	 
\end{itemize}

\begin{lemma} (Finite Identification) \label{finite_iden}
	Let $\{ \bfu^j \}_{j \in \mathbb{N}}$ be a sequence converging to a {\rm P}-stationary point $\hbu$ of \eqref{inner-sub}. Suppose that $\hbxi$ and $\nabla_{\bfxi} g (\hbu)$ satisfy strictly complementary condition \eqref{scc}, then there exists sufficiently large integer $\widehat{j}$ such that
	
	\begin{align} \label{T-G-S1}
		\Gamma_j = \mathcal{S}(\bfxi^{j+1/2}) = \mathcal{S}(\hbxi), \quad	\left\{ \begin{aligned}
			& T_j = \mathcal{S}(\bfw^{j+1/2}) = \mathcal{S}(\hbw), &&\ \mbox{if} \ \| \hbw \|_0 = s, \\
			& T_j \supseteq \mathcal{S}(\bfw^{j+1/2}) \supseteq \mathcal{S}(\hbw), &&\ \mbox{if} \ \| \hbw \|_0 < s, 
		\end{aligned}  \right. \quad  \forall j \geq \widehat{j},
	\end{align}
where $\mathcal{S}(\cdot)$ includes indices of nonzero elements for a given vector.
	
\end{lemma}

{\bf Proof}
Let us first prove the relationship involving $T_j$. Since $\lim_{j \to \infty} \bfu^{j} = \widehat{\bfu}$ and \eqref{inn-suc-chan}, we know that $\lim_{j \to \infty} \bfu^{j+1/2} = \widehat{\bfu}$ also holds. Then when $j$ is sufficiently large, we have $\mathcal{S}(\bfw^{j+1/2}) \supseteq \mathcal{S}(\hbw)$. The relationship $T_j \supseteq \mathcal{S}(\bfw^{j+1/2})$ directly follows from \eqref{gradient_xy}, and thus $T_j \supseteq \mathcal{S}(\bfw^{j+1/2}) \supseteq \mathcal{S}(\hbw)$ always holds when $j$ is sufficiently large. When $\| \hbw \|_0 = s$, $|T_j| = s$ indicates $T_j = \mathcal{S}(\bfw^{j+1/2}) = \mathcal{S}(\hbw)$.
	
We now prove $\Gamma_j = \mathcal{S}(\hbxi)$. From $\lim_{j \to \infty } \bfu^j = \hbu$, $\mathcal{S}(\bfxi^{j+1/2}) \supseteq \mathcal{S}(\hbxi)$ holds when $j$ is large enough. Noticing that \eqref{w-half} and \eqref{Tj} lead to $\mathcal{S}(\bfxi^{j+1/2}) = \Gamma_j$, we can obtain $\Gamma_j \supseteq \mathcal{S}(\hbxi)$. 
	
Finally, we need to prove $\Gamma_j \subseteq \mathcal{S}(\hbxi)$. 
Suppose for the contradiction that there exists a infinite index $\widehat{\J}$ and $\Gamma_j \not\subseteq \mathcal{S}(\hbxi)$ for any $j \in \widehat{\J}$. Then considering $| \Gamma_j | \subseteq [m]$ is finite, without loss of generality, we can assume that there exists a fixed index $\widehat{i} \in \Gamma_j$ but $\widehat{i} \notin \mathcal{S}(\hbxi)$ for all $j \in \widehat{\J}$. From \eqref{gradient_xy}, we have $[\nabla_\bfxi g( \bfu^j )]_{\widehat{i}} = -[\bfxi^{j+1/2} - \bfxi^j]_{\widehat{i}}/\beta$. Passing limit $j \to \infty$ on both sides of this equality leads to $[\nabla_\bfxi g( \hbu )]_{\widehat{i}} = 0$. Since $ \widehat{i} \notin \mathcal{S}(\widehat{\bfxi})$, $\hbxi_{\widehat{i}} = 0$ must hold. This contradicts to the strictly complementary assumption. Therefore, $\Gamma_j = \mathcal{S}(\hbxi)$ holds.
\hfill $\blacksquare$ 
	
\noindent\textbf{Proof of the first part in Theorem \ref{qua_convergence}.} 
 Let us first prove $[\nabla g ( \hbu )]_{\Upsilon_j} = 0$. Indeed, $[\nabla_\bfw g (\hbu)]_{T_j} = 0$ follows from \eqref{T-G-S1} and \eqref{proj_explicit}. $[\nabla_\bfxi g (\hbu)]_{\Gamma_j} = 0$ can be verified by \eqref{T-G-S1} and Lemma \ref{Lemma-Prox}. 
	
	Denote $H(t):= [\nabla^2 g( \hbu + t ( \bfu^{j+1/2} - \hbu ) )]_{\Upsilon_j, \Upsilon_j}$. 
	We analyze the relationship between $\| \widetilde{\bfu}^{j+1} - \hbu \| $ and $\| \bfu^{j+1/2} - \hbu \|$ below.
	\begin{align}
		\| \widetilde{\bfu}^{j+1} - \hbu \|  \mathop{=}\limits^{\eqref{T-G-S1}}  & \| [\widetilde{\bfu}^{j+1} - \hbu]_{\Upsilon_j} \| = \| [\bfu^{j+1/2} - \hbu]_{\Upsilon_j} - (H^{j+1/2})^{-1} [\nabla g( \bfu^{j+1/2} )]_{\Upsilon_j} \| \notag	\\
		\leq & \frac{1}{\sigma_g} \| H^{j+1/2} [\bfu^{j+1/2} - \hbu]_{\Upsilon_j} - [\nabla g( \bfu^{j+1/2} )]_{\Upsilon_j} \| \notag \\
		\leq & \frac{1}{\sigma_g} \|  H^{j+1/2} [\bfu^{j+1/2} - \hbu]_{\Upsilon_j} - [\nabla g( \bfu^{j+1/2} ) - \nabla g( \hbu )]_{\Upsilon_j} \| \notag \\
		\leq & \frac{1}{\sigma_g} \| \int_{0}^{1} ( H^{j+1/2} - H(t) ) [\bfu^{j+1/2} - \hbu]_{\Upsilon_j}  {\rm d}t \| \leq \frac{1}{\sigma_g}  \int_{0}^{1} L_g (1-t) \| \bfu^{j+1/2} - \hbu \|^2 {\rm d}t \notag \\
		\leq & \frac{L_g}{2\sigma_g}  \| \bfu^{j+1/2} - \hbu \|^2, \label{q_shrink}
	\end{align}
	where the first inequality is derived by using $\sigma_g$-strong convexity of $g$ and the fourth inequality follows from the Lipschitz continuity of $\nabla^2 g$. Then $\lim_{j \to \infty}  \tbu^{j+1} = \hbu$ directly follows from $\lim_{j \to \infty}  \bfu^{j+1/2} = \hbu$. We also have the following equations by using \eqref{T-G-S1} and \eqref{Newton-Condition} when $j$ is sufficiently large
	\begin{align*}
		&J(\widetilde{\bfxi}^{j+1}) = J(\widetilde{\bfxi}^{j+1}_{\Gamma_j}) = J(\hbxi_{\Gamma_j}) ,~~ J(\bfxi^{j+1/2}) = J(\bfxi^{j+1/2}_{\Gamma_j}) = J(\hbxi_{\Gamma_j}) .
	\end{align*}
	Finally, we prove that the descent property in \eqref{Newton-Condition} 
	 when $j$ is sufficiently large, and thereby the Newton step will always be adopted.
	\begin{align*}
		&G(\tbu^{j+1}) - G(\bfu^{j+1/2}) \\
		= & g( \tbu^{j+1} ) - g(\bfu^{j+1/2}) + \lambda J(\widetilde{\bfxi}^{j+1}) - \lambda J(\bfxi^{j+1/2}) \\
		%	= & \langle \nabla g( \bfu^{j+1/2} ), \tbu^{j+1} - \bfu^{j+1/2} \rangle + \frac{1}{2} ( \tbu^{j+1} - \bfu^{j+1/2} )^\top H^{j+1/2} ( \tbu^{j+1} - \bfu^{j+1/2} ) + o( \| \tbu^{j+1} - \bfu^{j+1/2} \|^2 ) \\
		\mathop{=}\limits^{\eqref{T-G-S1}} & \langle [\nabla g( \bfu^{j+1/2} )]_{\Upsilon_j}, [\tbu^{j+1} - \bfu^{j+1/2}]_{\Upsilon_j} \rangle + \frac{1}{2} [\tbu^{j+1} - \bfu^{j+1/2}]_{\Upsilon_j}^\top H^{j+1/2} [\tbu^{j+1} - \bfu^{j+1/2}]_{\Upsilon_j} \\
		& + o( \| \tbu^{j+1} - \bfu^{j+1/2} \|^2 ) \\
		\mathop{\leq}\limits^{\eqref{Newton-Condition}} & - \frac{1}{2} [\tbu^{j+1} - \bfu^{j+1/2}]_{\Upsilon_j}^\top H^{j+1/2} [\tbu^{j+1} - \bfu^{j+1/2}]_{\Upsilon_j} + o( \| \tbu^{j+1} - \bfu^{j+1/2} \|^2 ) \\
		\leq & -\frac{\sigma_g}{2} \| [\tbu^{j+1} - \bfu^{j+1/2}]_{\Upsilon_j} \|^2 + o( \| \tbu^{j+1} - \bfu^{j+1/2} \|^2 ) \\ \mathop{=}\limits^{\eqref{T-G-S1}} & -\frac{\sigma_g}{2} \| \tbu^{j+1} - \bfu^{j+1/2} \|^2 + o( \| \tbu^{j+1} - \bfu^{j+1/2} \|^2 ) \leq - \frac{\sigma_g}{4}  \| \tbu^{j+1} - \bfu^{j+1/2} \|^2,
	\end{align*} 
	where the second inequality follows from the $\sigma_g$-strong convexity of $g$ and the last inequality is derived form $\lim_{j \to \infty} \| \tbu^{j+1} - \bfu^{j+1/2} \| = 0$. \\

\noindent\textbf{Proof of the second part in Theorem \ref{qua_convergence}.} 
Notice that \eqref{q_shrink} has indicated the relationship between  $\| \bfu^{j+1}  - \hbu\|$ and $\| \bfu^{j+1/2}  - \hbu\|$. To prove the quadratic convergence, we just need to analyze the relationship between $\| \bfu^{j+1/2}  - \hbu\|$ and $\| \bfu^{j}  - \hbu\|$. By using \eqref{gradient_xy}, \eqref{T-G-S1} and $[\nabla g( \hbu )]_{\Upsilon_j} = 0$, we have the following estimation:
\begin{align*}
	& \begin{aligned}
		\| \bfw^{j+1/2} - \hbw \| = & \| [\bfw^{k+1/2} - \hbw]_{T_j} \| = \| [ \bfw^j - \hbw - \alpha \nwg (\bfu^j) ]_{T_j} \| \\
		= & \| [\bfw^j - \hbw - \alpha ( \nwg(\bfu^j) - \nwg(\hbu) )]_{T_j} \| \\
		\leq & \| \bfw^j - \hbw \| + \alpha \ell_g \| \bfu^j - \hbu \| \leq \underbrace{( 1 + \alpha \ell_g )}_{:= \zeta_2} \| \bfu^j - \hbu \| ,
	\end{aligned} \\
	& \begin{aligned}
		\| \bfxi^{j+1/2} - \hbxi \| = & \| [\bfxi^{j+1/2} - \hbxi]_{\Gamma_j} \| = \| [ \bfxi^j - \widehat{\bfxi} - \beta \nyxi ( \bfu^j ) ]_{\Gamma_j} \| \\
		= & \| [ \bfxi^j - \hbxi ]_{\Gamma_j} - \beta [ \nyxi ( \bfu^j) - \nyxi (\hbu)  ]_{\Gamma_j} \| \\
		\leq & \| \bfxi^j - \hbxi \| + \beta \ell_g   \| \bfu^{j+1/2} - \hbu \|  \leq  \underbrace{(1+ \beta \ell_g)}_{:=\zeta_3} \| \bfu^j - \hbu \| .
	\end{aligned}
\end{align*} 
These two results lead to
\begin{align*}
	\| \bfu^{j+1/2} - \hbu \|^2 = \| \bfw^{j+1/2} - \hbw \|^2 + \| \bfxi^{j+1/2} - \hbxi \|^2 \leq (\zeta_2^2 + \zeta_3^2) \| \bfu^j - \hbu \|^2,
\end{align*}
which combining with \eqref{q_shrink} implies local quadratic rate.   \hfill $\blacksquare$

%%%%%%%%%%%%%%%%%%%%%%%%%%%%%%%%%%%%%%%%%%%%%%%%%%%%%%%%%%%%%%
\bibliographystyle{siamplain}
\bibliography{SSVM_Refs} % if more than one, comma separated

%%%%%%%%%%%%%%%%%%%%%%%%%%%%%%%%%%%%%%%%%%%%%%%%%%%%%%%%

\end{document}